\documentclass{article}

\usepackage{arxiv}

\usepackage[utf8]{inputenc} %
\usepackage[T1]{fontenc}    %
\usepackage{pifont}%
\usepackage{hyperref}       %
\usepackage{lmodern}
\usepackage{url}            %
\usepackage{booktabs}       %
\usepackage{nicefrac}       %
\usepackage[dvipsnames]{xcolor}         %
\usepackage[pdftex]{graphicx}
\usepackage{amsmath}
\usepackage{amsthm}
\usepackage{mathtools}
\usepackage{bm}
\usepackage{bbm}
\usepackage{enumitem}
\usepackage{cleveref}
\usepackage{amssymb}
\usepackage[mathscr]{eucal}
\usepackage[backend=bibtex, 
style=authoryear,
hyperref = true, 
backref=true,
doi=false, isbn=false, url=false, eprint=false, dashed=false,
maxcitenames=2]{biblatex}  %
\usepackage{crossreftools}
\hypersetup{%
    bookmarksnumbered, bookmarksopen=true, bookmarksopenlevel=1,%
}

\definecolor{burgundy}{rgb}{0.5, 0.0, 0.13}
\definecolor{dark_burgundy}{rgb}{0.6, 0.0, 0.05}
\definecolor{camel}{rgb}{0.76, 0.6, 0.42}
\definecolor{chamoisee}{rgb}{0.63, 0.47, 0.35}
\definecolor{grey1}{RGB}{128,128,128}
\definecolor{navy}{rgb}{0.0, 0.0, 0.5}
\definecolor{teal}{rgb}{0.0, 0.5, 0.5}
\definecolor{forest}{rgb}{0.13, 0.55, 0.13}
\definecolor{crimson}{rgb}{0.86, 0.08, 0.24}
\hypersetup{colorlinks=true,linkcolor=burgundy,citecolor=chamoisee,urlcolor=burgundy,linktoc=page}

\newtheorem{assumption}{Assumption}%
\newtheorem{remark}{Remark}[section]

\newtheorem{definition}{Definition}%
\newtheorem{proposition}{Proposition}[section]
\newtheorem{lemma}{Lemma}[section]
\newtheorem{theorem}{Theorem}

\newtheorem*{remark*}{Remark}

\numberwithin{equation}{section}

\newcommand*{\defeq}{\mathrel{\vcenter{\baselineskip0.5ex \lineskiplimit0pt
			\hbox{\scriptsize.}\hbox{\scriptsize.}}}%
	=}
\newcommand{\<}{\langle}
\renewcommand{\>}{\rangle}
\DeclareMathOperator{\Supp}{Supp}

\DeclareMathOperator{\Expectation}{\mathbb{E}}

\DeclareMathOperator{\Id}{Id}
\DeclareMathOperator*{\essentialsup}{ess\,sup}
\DeclareMathOperator*{\essentialinf}{ess\,inf}
\DeclareMathOperator{\LPC}{LPC}

\newcommand{\K}{{s_0}}
\def\target{\mu^0}
\def\targetdensity{f^0}
\def\estimator{\widehat\mu}
\newcommand{\empiricaldist}{\hat{f}_n}

\DeclareMathOperator*{\argmin}{arg\,min}
\newcommand*\diff{\mathop{}\!\mathrm{d}} %
\newcommand{\Measures}{\mathscr{M}}

\newcommand{\nearregion}{\mathcal{N}^{\textnormal{reg}}}
\newcommand{\farregion}{\mathcal{F}^{\textnormal{reg}}}

\newcommand{\Certificate}{\eta^0}

\newcommand{\Hilbert}{\mathcal{H}}
\newcommand{\model}{\textnormal{mod}}

\newcommand{\pivot}{\textnormal{pivot}}
\newcommand{\switch}{\textnormal{switch}}

\newcommand{\kernel}{\lambda}
\newcommand{\bandwidth}{\tau}

\DeclareMathOperator{\sinc}{sinc}

\DeclareMathOperator{\fourier}{\mathfrak{F}}

\newcommand{\metric}{\mathfrak{g}}

\bmdefine{\bUpsilon}{\Upsilon}

\newcommand{\bI}{\bm{I}}

\newcommand{\R}{\mathbb{R}}

\newcommand{\rawobs}{z}
\bmdefine{\obs}{\rawobs}
\bmdefine{\Obs}{Z}

\newcommand{\distance}{{\mathfrak{d}_{\mathfrak g}}}
\newcommand{\distancegeneric}{{\mathfrak{d}}}
\newcommand{\rawparam}{x}
\bmdefine{\param}{\rawparam}
\newcommand{\trueparam}{\param^0} %
\newcommand{\Param}{\mathcal{X}} %

\newcommand{\templatedensity}{\phi} %

\newcommand{\templatefourier}{\fourier[\templatedensity]} %
\newcommand{\nsketch}{m}

\newcommand{\Forward}{F}

\bmdefine{\bsketch}{\bm{y}_{\textnormal{sketch}}}
\newcommand{\SketchDist}{\Lambda}

\bmdefine{\ba}{a}
\bmdefine{\bb}{b}
\bmdefine{\bc}{c}

\bmdefine{\bu}{u}
\bmdefine{\bt}{t}
\bmdefine{\bv}{v}
\bmdefine{\bs}{s}
\bmdefine{\bq}{q}
\bmdefine{\by}{y}
\bmdefine{\bz}{z}
\bmdefine{\bx}{x}
\bmdefine{\bw}{w}
\bmdefine{\bS}{S}
\bmdefine{\bY}{Y}
\bmdefine{\bZ}{Z}
\bmdefine{\bF}{F}
\bmdefine{\bI}{I}
\bmdefine{\bomega}{\omega}
\bmdefine{\bBeta}{\beta}
\bmdefine{\balpha}{\alpha}
\bmdefine{\bgamma}{\gamma}
\bmdefine{\br}{r}
\bmdefine{\bxi}{\xi}

\bmdefine{\btau}{\tau}
\bmdefine{\btheta}{\theta}
\bmdefine{\bTheta}{\Theta}
\bmdefine{\bmu}{\mu}
\bmdefine{\bSigma}{\Sigma}

\newcommand{\indicator}{\mathbbm{1}}

\newcommand{\Model}{\mathcal{M}} %

\Crefname{pluralequation}{Equations}{Equations}
\crefformat{pluralequation}{equations~(#2#1#3)}
\Crefformat{pluralequation}{Equations~(#2#1#3)}

\makeatletter
\def\arrowfill@@#1#2#3#4{%
  $\m@th\thickmuskip0mu\medmuskip\thickmuskip\thinmuskip\thickmuskip
   \relax#4#1
   \xleaders\hbox{$#4#2$}\hfill
   #3$%
}
\def\leftarrowfill@@{\arrowfill@@\leftarrow\relbar\relax}
\newcommand{\xdashleftarrow}[2][]{\ext@arrow 3095\leftarrowfill@@{#1}{#2}}
\makeatother

\newcommand{\nocontentsline}[3]{}
\let\origcontentsline\addcontentsline
\newcommand\stoptoc{\let\addcontentsline\nocontentsline}
\newcommand\resumetoc{\let\addcontentsline\origcontentsline}

\graphicspath{{./figures/}}

\author{%
 Yohann De Castro\\
 Institut Camille Jordan - IUF - Centrale Lyon%
  \And
  Rémi Gribonval \\
  Inria, CNRS, ENS de Lyon, Université Claude Bernard Lyon 1, LIP, UMR 5668, 69342, Lyon cedex 07, France
 \And
  Nicolas Jouvin \\
  INRAE, AgroParisTech, Université Paris-Saclay, MIA-Paris, 91120, Palaiseau
}

\title{Effective regions and kernels in continuous sparse regularisation, with application to sketched~mixtures}

\begin{document}

\maketitle
\begin{abstract}
    This paper advances the general theory of continuous sparse regularisation on measures with the Beurling-LASSO (BLASSO). This TV-regularised convex program on the space of measures allows to recover a sparse measure using a noisy observation from a measurement operator. While previous works have uncovered the central role played by this operator and its associated kernel in order to get estimation error bounds, the latter requires a technical local positive curvature (LPC) assumption to be verified on a case-by-case basis. In practice, this yields only few LPC-kernels for which this condition is proved. In this paper, we prove that the ``sinc-4'' kernel, used for signal recovery and mixture problems, does satisfy the LPC assumption. Furthermore, we introduce the \textit{kernel switch} analysis, which allows to leverage on a known LPC-kernel as a \emph{pivot} kernel to prove error bounds. Together, these results provide easy-to-check conditions to get error bounds for a large family of translation-invariant model kernels.   
    Besides, we also show that known BLASSO guarantees can be made adaptive to the noise level. This improves on known results where this error is fixed with some parameters depending on the model kernel. We illustrate the interest of our results in the case of mixture model estimation, using band-limiting smoothing and sketching techniques to reduce the computational burden of BLASSO. 
\end{abstract}

\section{Introduction}
\label{sec:intro}
Sparse regularisation, particularly methods based on $\ell^1$ minimisation, has emerged as a cornerstone technique to solve high-dimensional linear inverse problems across numerous domains. These approaches have revolutionized signal processing, statistical learning, and computational imaging by enabling reconstruction from severely limited observations~\parencite{foucart2013mathematical, candes2006robust}. Despite their remarkable success, traditional sparse methods rely on discretisation, introducing inherent limitations in both computational efficiency and theoretical precision. Consequently, there has been growing interest in continuous counterparts to $\ell^1$ minimisation that operate in a grid-free manner. These "off-the-grid" approaches not only mitigate discretisation errors but also potentially offer more efficient solvers, sharper theoretical guarantees, and increased robustness to sampling limitations~\parencite[see \textit{e.g.}][]{candes2014towards, tang2013compressed,duval2015exact,duval2017sparse}. Such continuous formulations have proved particularly valuable when the underlying target structure inherently exists in a continuous domain, such as in the case of mixtures of distributions or when the parameters of interest are not naturally discretized.

A promising extension of $\ell^1$ minimisation to continuous domains involves encoding both positions and weights as a discrete Radon measure, replacing the $\ell^1$ norm with the total variation norm. This approach treats a measure as "sparse" when it consists of a finite sum of few Dirac masses. The resulting infinite-dimensional optimisation problem is cast over the space of measures and known as Beurling-LASSO (BLASSO). The latter has emerged as a powerful framework across multiple disciplines. In signal processing, BLASSO has demonstrated effectiveness in sparse deconvolution and super-resolution \parencite{candes2014towards,duval2015exact,boyer2017adapting}. Its applications extend to statistics for continuous sparse regression and off-the-grid compressed sensing \parencite{decastro2012exact,poon2023geometry}, and to machine learning for mixture model estimation \parencite{decastro2019sparse}. Practical implementations have brought significant advances in diverse fields: single-molecule fluorescent microscopy \parencite{boyd2017alternating}, neuronal spike sorting \parencite{ekanadham2014unified}, statistical mixture modelling \parencite{gribonval2021compressive,decastro2019sparse}, and training of shallow neural networks \parencite{bach2017breaking}. In this paper, we motivate our analysis and illustrate our results in the case of mixture model estimation \parencite{decastro2019sparse}, although our theoretical contributions apply broadly to all BLASSO applications.

The BLASSO objective is a TV-regularised data-fitting term, the latter inducing a \textit{model} kernel which is tied to the measurement process--or model--at hand, \textit{e.g.} compressed sensing, deconvolution, mixture model estimation etc. One of the main practical limitation of BLASSO is the computation cost of large kernel matrices, see for instance \textcite[Section 1.3]{decastro2023fastpart}. This well-known issue in kernel methods has motivated the use of sketching techniques~\parencite{rahimi2007random,gribonval2020sketching}, whose principle is to compress a dataset into a single vector, the \textit{sketch}, by computing averages (over the samples of the dataset) of~$m$ non-linear random features. This compression is designed to approximate the kernel computations while retaining statistical guarantees for the underlying learning task. For continuous sparse regression, BLASSO analyses for \textit{sketched} model kernels defined via random features have been investigated in~\textcite{poon2019support,poon2023geometry}, leading to the development of a corresponding \textit{sketched} BLASSO formulation. 
Realted albeit related recent studies \parencite{gribonval2021compressive,ayoub2024}have proven that the model sketching operator, mapping from a dataset to its sketch, acts as a quasi-isometry when restricted to the set of sparse measures whose support points are sufficiently separated.

\subsection{Main contributions}
\label{sec:main_contribution_intro}

We begin by presenting a schematic view of BLASSO theoretical analysis presented in this paper. For the sake of presentation, most of of the objects will not be properly defined as we wish to focus on the articulation between the state-of-the-art and our contributions. 

The BLASSO is a convex optimisation program of the space of measures, denoted by $\mathcal{P}_K$ hereafter, which depends on a reproducing kernel $K$. The goal of BLASSO theoretical analysis is to obtain error bounds on the recovery of target sparse measures $\target$ belonging to a model set $\Model_{\Delta, K}$. This model set involves a minimal separation $\Delta>0$ between atoms of $\target$, measured with respect to the Fisher-Rao distance $\distancegeneric_{K}$ of the kernel $K$. The crux of every theoretical work is to prove the existence of so-called \emph{non-degenerate dual certificates} for the BLASSO program~$\mathcal{P}_K$ and the model set~$\Model_{\Delta, K}$. These objects guarantees that the support of the BLASSO estimator ``localizes'' around the support of the target sparse measure~$\target$. This localisation error is stated in term of \emph{near regions} with a radius~$r>0$ around the true atoms, the lower the better the recovery.  In their comprehensive treatment, \textcite{poon2023geometry} formalized the local positive curvature assumption, $\LPC(K, \Delta)$, a technical condition on kernel~$K$ to ensure the existence of such certificates. Their $\LPC$ analysis proves the existence of certificates for the program~$\mathcal{P}_K$, as well as the \textit{sketched} program $\mathcal{P}_{\tilde K}$ when using sketching techniques to approximate the kernel $K \approx \tilde K$ with random features~\parencite{rahimi2007random, gribonval2020sketching}. Notably, in the sketching case, the existence of certificates holds with high-probability with respect to the sketch size, however the model set~$\Model_{\Delta, K}$ still depends on the kernel $K$ and not on its sketched version $\tilde K$.

Although appealing, this analysis requires to prove the $\LPC(K, \Delta)$ and to compute the Fisher-Rao distance $\distancegeneric_{K}$ which can be cumbersome and has to be done on a case-by-case basis. To date, only a limited set of kernels have been rigorously proved to enjoy the $\LPC$ property. Furthermore, the radius of the localization result is fixed and cannot adapt to the noise level. 

 We address these issues by introducing both the notion of \emph{effective} near regions and the \emph{kernel switch} technique. Additionally, we prove the $\LPC$ for the sinc-4 kernel, thereby incrementing the list of known $\LPC$-kernel. We also demonstrate the interest of these results in the case of mixture model estimation.

\paragraph{Effective near-regions and near-minimizers} We show that the radius $r$ of near region is {\em adaptive to the noise level}, denoted by $\gamma$, and that BLASSO can localize \textit{effective} near-regions up to~$r\sim\sqrt{\gamma\sqrt{\K}}$, where $\K$ is the sparsity level of the target measure $\target$. This is particularly relevant in statistical applications where the noise level decreases with sample size, denoted by $n$. As an example, in the mixture problem, the noise~$\gamma=\mathcal O(n^{-1/2})$, and one gets effective near regions of size essentially $\mathcal O({n^{-1/4}})$, narrowing around the true parameters as $n$ increases. Additionally, we also uncover an interesting feature of our convex programming problem. Namely, we establish that near minimizers of the BLASSO problem—defined as solutions with an objective value no greater than that of the target—share the same guarantees as the exact minimizer itself. In practice, this means one does not need to solve the BLASSO exactly, any sufficiently \textit{near solution} will do.

\begin{figure}[!b]
    \centering
    \includegraphics[width=0.8\linewidth]{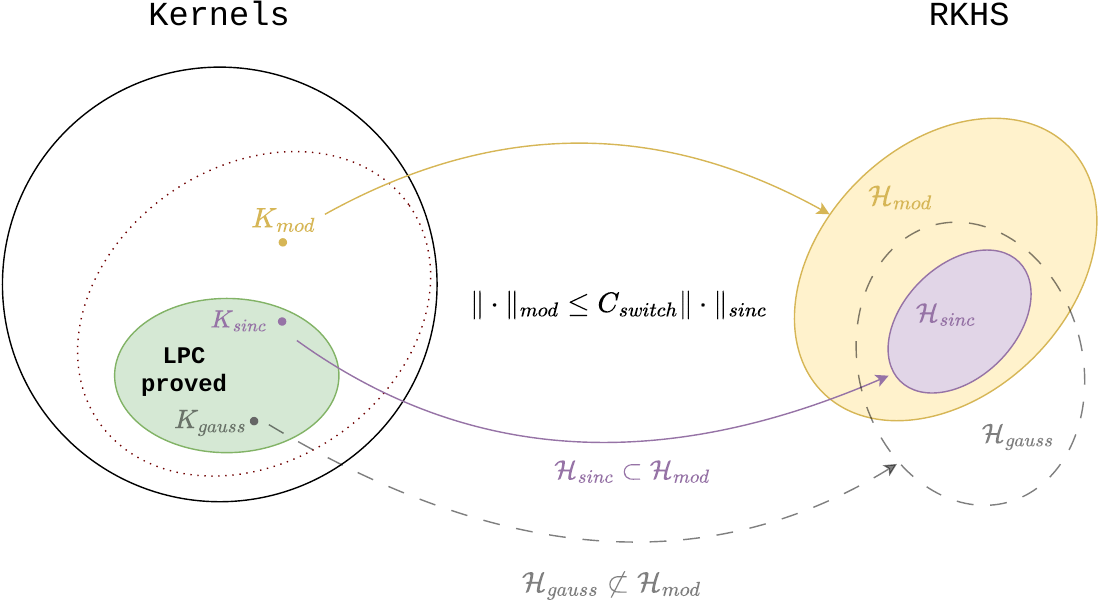}
    \caption[Illustration of the kernel switch principle.]{
    Illustration of the kernel switch principle. Left: A small set of pivot kernels known to satisfy the local positive curvature assumption ($\LPC$) enables theoretical guarantees for a much larger class of model kernels. Right: A model kernel $K_{\model}$ is admissible for our statistical guarantees if there exists a pivot LPC kernel (here~$K_{\textnormal{sinc}}$) whose RKHS is continuously embedded into the model RKHS. The error bounds remain valid with an additional scaling factor $C_{\textnormal{switch}}(K_\model, K_{\textnormal{sinc}})$. In this example, using the sinc kernel as pivot is valid because its RKHS is included in the model RKHS, while the Gaussian kernel cannot serve as a pivot since its RKHS contains functions with faster-decaying Fourier transforms than the model allows.
    }
    \label{fig:kernel_switch}
\end{figure}

\paragraph{Kernel switch}
This paper introduces a significant advancement by demonstrating that an existing kernel $K_0$ satisfying $\LPC(K_0,\Delta)$--termed a \textit{pivot kernel}--can be used to derive guarantees for BLASSO problems of interest~$\mathcal{P}_K$ involving a \textit{different} model kernel $K$, including the case of a sketched model kernel $\tilde K$ defined via sketching. This \textit{kernel switch} principle, illustrated in Figure~\ref{fig:kernel_switch}, is possible provided that the RKHS of $K_0$ is continuously embedded into the RKHS of $K$, which is ensured by a finite embedding constant: $C_\switch(K, K_0) < + \infty$. Under this condition, we prove the existence of certificates both for the program $\mathcal{P}_{K}$ and the sketched program~$\mathcal{P}_{\tilde K}$, with the model set $\Model_{\Delta, K_0}$ defined via the pivot kernel $K_0$ and its $\LPC(K_0, \Delta)$ condition. Our result yields the same localization error bounds as in the direct approach of proving $\LPC$, only factoring in the embedding constant $C_\switch$. This substantially broadens the scope of continuous sparse regression by allowing practitioners to leverage the established theoretical guarantees of well-characterized pivot kernels, while employing model kernels that might be more computationally advantageous or better suited to specific problem structures. 

\paragraph{The ``sinc-4'' pivot} We rigorously show for the first time that a specific kernel $K_0$, the sinc-4 defined as the fourth power of the sinus cardinal kernel~\parencite{decastro2019sparse}, satisfies the~LPC assumption and is amenable to sketching with random features. This is particularly helpful as: \textit{a) its RKHS is included in the well-known space of band-limited functions}; and \textit{b) its Fisher-Rao distance is a simple rescaled euclidean distance}, leading to a simple model set $\Model_{\Delta, \textnormal{sinc-4}}$ defined in terms of Euclidean separation between atoms. Coupled with our kernel switch analysis, this result implies that \emph{any BLASSO problem $\mathcal{P}_K$ with $K$ a translation-invariant kernel whose RKHS contains band-limited functions verifies the theoretical guarantees} resulting from our analysis. This also extends to sketched translation-invariant kernels under mild condition on the random features.

\paragraph{Contributions for Mixture Modelling} 
Our contribution on sketching sparse regularisation of translation-invariant mixtures is detailed in~\Cref{sec:Mixtures} with a special focus on the sinc-4 kernel. This lead to S2Mix, a more practical algorithm for the \textit{sketched} Supermix problem of~\textcite{decastro2019sparse}. 
The consequences of our results are as follows: S2Mix is a computationally efficient algorithm to estimate the parameters of a mixture model, even for large datasets. For a mixture with~$\K$ component, it achieves a sketch size complexity of~$\mathcal O(\K\log^2(\K))$, 
which can be significantly lower than the original sample complexity $n$. Furthermore, S2Mix covers a wide range of mixture models, including non-Gaussian template distributions, since it only requires that its characteristic function does not vanish in low-frequencies, a friendly condition to be verified. This makes it a versatile tool for various applications in machine learning and statistics. The algorithm is robust to noise and model misspecification, as it can handle a wide range of noise levels and model assumptions. This robustness is particularly important in real-world applications where data may be noisy or incomplete. Finally, S2Mix provides theoretical guarantees for parameter recovery, ensuring that the estimated parameters are close to the true model parameters with high probability.

\subsection{Related works}

In their pioneering work, \textcite{candes2014towards} studied the super-resolution deconvolution problem on the torus. Their proof technique can be reread under the viewpoint developed in this paper as using a kernel switch with the Dirichlet kernel as model kernel and the Jackson kernel (Dirichlet kernel raised to the power of~4) as pivot kernel. Building directly on this methodology,~\textcite{decastro2019sparse} utilized a similar proof technique in their detailed analysis of mixture model estimation with BLASSO, relying on the sinc-4 kernel as pivot. Notably, both papers implicitly leverage a kernel switch analysis without formally characterizing the general principles of this phenomenon. Furthermore, \textcite{decastro2019sparse} did not provide a proof that the sinc-4 kernel satisfies the $\LPC$ assumption, making it impossible to derive sketching statistical error bounds using the framework later developed by~\textcite{poon2023geometry}.

For Gaussian mixture models with fixed covariance, \textcite[Section~2.3 and Appendix~D.5.1]{poon2023geometry} showed that the Gaussian kernel with the same covariance satisfies the $\LPC$ assumption. Their BLASSO analysis, however, is essentially restricted to Gaussian mixtures where the empirical distribution is convolved with a prescribed Gaussian (referred to as smoothing kernel in this paper, see Section~\ref{sec:smoothing}). In contrast, our kernel switch analysis (under the mild assumptions in Section~\ref{sec:Mixtures}) covers a much broader class of mixture densities and smoothing kernels. Recently, building on our effective near-region analysis, \textcite{giard2025gaussian} further established BLASSO error bounds for Gaussian mixtures with \emph{unknown} diagonal covariances by leveraging the Fisher-Rao geometry of the model kernel via a novel semi-distance tailored to the $\LPC$ condition.

Although the near-minimizer property is implicit in the proofs of previous works on Bregman divergence \parencite{candes2014towards, duval2015exact,duval2017sparse, azais2015spike,poon2023geometry}, it has not been explicitly stated. A similar concept exists in the Compressed Sensing (CS) literature, where guarantees for exact minimizers of LASSO-related problems are often proven via more general results that also apply to near minimizers (see e.g., \cite[Theorems 4.18 \& 4.19, 4.21 \& 4.23]{foucart2013mathematical}).

The near-minimizer property is of great practical importance, as it justifies the use of approximate solvers. This is especially relevant for large-scale or high-dimensional problems where finding the exact minimizer is computationally infeasible. Recent algorithms to solve the BLASSO problem, such as Sliding Frank-Wolfe \parencite{denoyelle2019sliding} and stochastic versions of Conic Particle Gradient Descent (CPGD) \parencite{chizat2019sparse,decastro2023fastpart}, are designed precisely to find such near minimizers with convergence guarantees. However, while these CPGD-based methods can find a near minimizer in a finite number of steps, the required number of steps scales exponentially with the dimension $d$, currently limiting their application to moderate-dimensional regimes where $d$ is less than approximately $100$.

We illustrate our contribution on mixture models. The latter represent a fundamental framework in modern statistics, with diverse applications spanning machine learning, biology, genetics, and astronomy \parencite{mclachlan2004finite,fruhwirth2019handbook}. Traditional inference for these models typically relies on the non-convex maximum-likelihood paradigm solved via the Expectation-Maximisation (EM) algorithm \parencite{mclachlan2007algorithm}. A critical challenge in mixture modelling—determining the optimal number of components~$\K$—is conventionally addressed post-inference as a model selection problem \parencite{celeux2018model}. To handle large-scale datasets, the computational statistics community has developed sophisticated methodological tools including variational inference \parencite{blei2017variational} and amortized inference techniques \parencite{kingma2014auto,tomczak2024deep}, though theoretical guarantees for parameter recovery remain relatively limited. The sparse continuous regression perspective on mixture modelling, recently introduced by \textcite{decastro2019sparse, poon2023geometry}, offers compelling theoretical guarantees for parameter recovery while demonstrating robustness to noise and model misspecification. This approach complements practical algorithms developed to solve the BLASSO optimisation problem, such as sliding Frank-Wolfe \parencite{denoyelle2019sliding} and particle-based methods \parencite{chizat2019sparse, decastro2023fastpart}, providing both theoretical foundations and computational efficiency.

\subsection*{Outline}
The remainder of the paper is organized as follows:
\begin{itemize}
    \item Section~\ref{sec:sota_blasso} introduces the BLASSO along with state-of-the-art theoretical guarantees that can be obtained. The notions of model kernel, model set and sketching are properly defined; we also review kernels known to satisfy the $\LPC$ assumption. 
    \item Section~\ref{sec:Mixtures} introduces the sketching BLASSO problem for mixture models. We first present a resolution-based approach to mixture modelling and show how it can be formulated as a continuous sparse regularisation problem. We then develop a sketched version of the framework, called S2Mix, which achieves significant computational efficiency through random Fourier features. We derive statistical guarantees for this approach, demonstrating how the sinc-4 kernel serves as an effective pivot kernel that enables parameter estimation with near-optimal sample complexity.
    \item Section~\ref{sec:BlassoResults} presents our three main theoretical contributions. First, our effective near regions analysis in Theorem~\ref{thm:main_thm} demonstrates that the size of near regions can adapt to the noise level, yielding localisation guarantees of order $\mathcal{O}(\sqrt{\gamma\sqrt{\K}})$ that improve as noise decreases. Second, our kernel switch technique in Theorem~\ref{thm:main_thm} establishes that statistical error bounds hold when using a different kernel for analysis than for measurement, provided the RKHS embedding condition (Assumption~\ref{hyp:pivot_kernel}) is satisfied. Theorem~\ref{thm:main_thm_sketch} extends the kernel switch to the case of sketched kernels under mild assumption on their random features. Third, we prove in both theorems that near-optimal solutions—\textit{i.e.,}~measures with objective value not exceeding that of the target measure—enjoy the same statistical guarantees as the exact minimizer. This section also includes a rigorous analysis proving that the sinc-4 kernel satisfies the $\LPC$ assumption, and can therefore be used as a pivot in a kernel switch analysis.
\end{itemize}

Technical proofs, additional mathematical background, and supporting material are provided in the appendix. A comprehensive list of notation precedes the references.

\section{State-of-the-art on continuous sparse regularisation with BLASSO}
\label{sec:sota_blasso}
Let $\Param\subset\R^d$ be a compact set, referred to as the \emph{parameter space}, and consider $( \Measures(\Param),\|\cdot\|_{\mathrm TV} )$ the space of Radon measures, defined as the topological dual space of $(\mathcal C(\Param),\|\cdot\|_{\infty})$, the continuous functions endowed with the infinity norm. Let~$\mathcal{F}$ be a separable Hilbert space and let $F\,:\,\Measures(\Param)\to \mathcal{F}$ be a linear map, referred to as the forward measurement operator. We define the BLASSO problem as
\begin{align}
\label{eq:blasso}
\addtocounter{equation}{1}
\estimator \in  \argmin_{\mu \in \Measures(\Param)} J_\kappa(\mu)
\quad\textnormal{where}\quad 
J_\kappa(\mu)\defeq\frac{1}{2} \big\Vert \by - \Forward \mu \big\Vert_{\mathcal{F}}^2 + \kappa \Vert \mu \Vert_{\mathrm TV}\,, 
\tag{\theequation--$J_\kappa$}
\end{align}
with $\kappa>0$ a regularisation parameter and $\by\in\mathcal{F}$ some observation. One can prove that a solution $\estimator\in\Measures(\Param)$ to BLASSO~\eqref{eq:blasso} exists under the forthcoming hypothesis~\eqref{hyp:continuous_model_kernel}, see for instance \cite[Theorem~1.1]{decastro2023fastpart}.

The criterion to be minimized decomposes as the sum of two terms: a \emph{data-fitting} term, comparing some observation $\by$ with a candidate~$\Forward \mu$ in the separable Hilbert space $\mathcal{F}$; and \emph{TV-regularisation} term, which acts as a sparsity inducing penalty analogous to $\ell^1$-regularisation in finite dimension \parencite[Theorem 1.1]{foucart2013mathematical,boyer2019representer}. 

\subsection{Model kernel, class of sparse models and inverse problem under consideration}

The continuous regression problem under consideration offers a powerful generalisation of classical sparse regression, encoding both weights and positions within a Radon measure framework. Central to our analysis is the \emph{model kernel}~$K_\model$, which characterizes the interaction between point masses in the parameter space $\Param$ through the measurement operator $F$. Although it is traditionally referred to simply as "the kernel" in the BLASSO literature~\parencite{poon2023geometry}, we choose to deliberately distinguish it from what we will call the \emph{pivot kernel} later in our analysis. 

The model kernel~$K_\model$ is formally defined as the $\mathcal{F}$-inner product of the operator $F$ applied to pairs of Dirac masses at different locations in the parameter space~$\Param$. Its continuity is a mild technical condition which ensures the existence of solutions to the BLASSO optimisation problem \parencite{decastro2023fastpart}. This distinction between \emph{model kernel} and \emph{pivot kernel} represents one of our key contributions, allowing us to extend theoretical guarantees to a broader class of inverse problems than previously possible.

\medskip

\begin{definition}[Model kernel]
\label{def:model_kernel}
    Define the \textit{model kernel} by
        \begin{equation}
        \label{def:kmod_general}
            K_\model(\bs, \bt)\defeq\langle F\delta_\bs,F\delta_\bt \rangle_{\mathcal{F}}\quad \textnormal{for all }\bs,\bt\in\Param\,,
        \end{equation}
where $\delta_\bu$ is the Dirac mass at point $\bu\in\Param$. We assume that 
\begin{equation}
    \label{hyp:continuous_model_kernel}
    \tag{\theequation --$\mathbf{H}_{\model}$}
    K_\model(\cdot,\cdot)\text{ is continuous,}
    \addtocounter{equation}{1}
\end{equation}
where continuity is understood with respect to the pair of variables. It is a positive definite kernel and we denote as~$\Hilbert_{\model}$ its reproducing kernel Hilbert space (RKHS), which is continuously included in $\mathcal{F}$ as discussed in~\Cref{sec:BlassoResults}.
\end{definition}

In the following, the class of models $\Model_{s, \Delta,\distancegeneric}$ is defined as the set of measures that can be written as a sum of $s$ Dirac masses with a \emph{minimal separation} of $\Delta$ between two spikes for some distance $\distancegeneric$; the noise term $\Gamma$ is defined as the difference between the observation $\by$ and the operator~$F$ applied to the \emph{target} measure $\target$; and the \emph{noise level} $\gamma$ is defined as the norm of the noise term. 

\medskip

\begin{definition}[Sparse models]
	Let $\distancegeneric(\cdot,\cdot)$ be a distance defined on $\Param$. Let $\Delta$ be positive and let $s$ be greater than~$1$. We define the class of $s$-sparse measures with minimal separation $\Delta$ as
	\begin{equation}
    \label{eq:class_model}
	\Model_{s, \Delta, \distancegeneric} 
        \defeq 
        \Big\{ \mu\,:\,
            \mu = \sum_{k=1}^{s} a_k \delta_{\bx_k}
               \textnormal{ and }
            \min_{k,\ell} 
            \distancegeneric(\bx_k,\bx_\ell)
            \geq \Delta
        \Big\}\,,
	\end{equation}
 where $\Delta$ is referred to as the minimal separation, and $s$ as the sparsity.
\end{definition}

\begin{remark}
    The class of sparse models is usually defined with respect to the so-called \emph{Fisher-Rao} distance associated to some positive type kernel~$K(\cdot,\cdot)\,:\,\Param\times\Param\to\R$. While previous works choose this kernel to be the model kernel~\parencite{poon2019support, poon2023geometry}, this paper shows error bounds choosing a \underline{possibly different}~kernel. 
\end{remark}

The inverse problem under consideration is defined as follows. We consider a target measure $\target\in\Model_{\K, \Delta_0, \distancegeneric}$, for some sparsity $\K$, distance $\distancegeneric$ and minimal separation $\Delta_0$. Let a noisy observation $\by$ from $\target$ through the measurement operator $F$, we define the noise term and noise  respectively~as
\begin{subequations}
\label[pluralequation]{eqs:noise_level_term}
\begin{align}
    \Gamma &\defeq \by-F\target\label{eq:noise}\,,\\
    \text{resp. }\gamma_0 &:=\|\Gamma\|_{\mathcal{F}}\,.\label[pluralequation]{eq:noise_level}
\end{align}
\end{subequations}
Given a \emph{noisy} observation $\by$, the goal is to recover the true locations $\trueparam_k$ and weights $a_k^0$ using a solution $\estimator$ of BLASSO \eqref{eq:blasso}. The noise term $\Gamma$ in \eqref{eqs:noise_level_term} accounts for model misspecification, outliers, and sampling or measurement noise \parencite{decastro2019sparse,poon2023geometry}. Throughout this paper, we call \emph{noise level} some upper bound $\gamma$ on the noise $\gamma_0$, namely $\gamma\geq\gamma_0$. Error bounds are stated in terms of the noise level~$\gamma$ and hold for any $\kappa$; %
setting $\kappa \simeq \gamma/\sqrt{\K}$ yields sharper constants, as presented in this paper. In statistical settings, standard concentration gives $\gamma = \mathcal O_{\mathbb P}(1/\sqrt{n})$ where $n$ is the sample size (we will uncover such result for mixture models later in this paper), while in deterministic inverse problems $\gamma$ is often specified by the measurement device.

\subsection{Estimation error bounds via dual certificates}

The analysis of BLASSO estimation error bounds \eqref{eqs:error_bounds} below hinges on the concept of the so-called \emph{far} and \emph{near} regions, which define spatial partitioning around the support points of the target measure $\target$. These regions serve a fundamental role in characterizing the behaviour of any solution $\estimator$ to \eqref{eq:blasso}. Specifically, the BLASSO exhibits strong reconstruction properties in near regions—where $\estimator$ closely approximates $\target$—while maintaining controlled behaviour in far regions. This spatial characterisation provides precise quantitative guarantees on the support recovery capabilities of the estimator, referred to as the localisation property of BLASSO. We formally define these regions as follows.

\medskip

\begin{subequations}\label[pluralequation]{eqs:near_far_def}
\begin{definition}[Far and Near regions]
\label{def:Far_and_Near_regions}
	Let $\target \in \Model_{\K, \Delta_0, \distancegeneric}$, let $r > 0$, we define the near region of $\trueparam_k$ as the closed ball centered at $\trueparam_k$ for the distance $\distancegeneric$
	\begin{equation}
        \label{eq:near_regions}
	\nearregion_k(r) \defeq \left\{ \param \in \Param, \quad 
        \distancegeneric(\param,\trueparam_k)
        \leq r \right\},
	\end{equation}
	and the far region as the complement
	\begin{equation}
         \label{eq:far_regions}
	\farregion(r) \defeq \Param \setminus \nearregion(r), \quad \textnormal{with: } \nearregion(r) = \bigcup_{k=1}^{\K} \nearregion_k(r).
	\end{equation}
\end{definition}
\end{subequations}

\begin{figure}[!ht]
    \centering
    \includegraphics[width=0.65\linewidth]{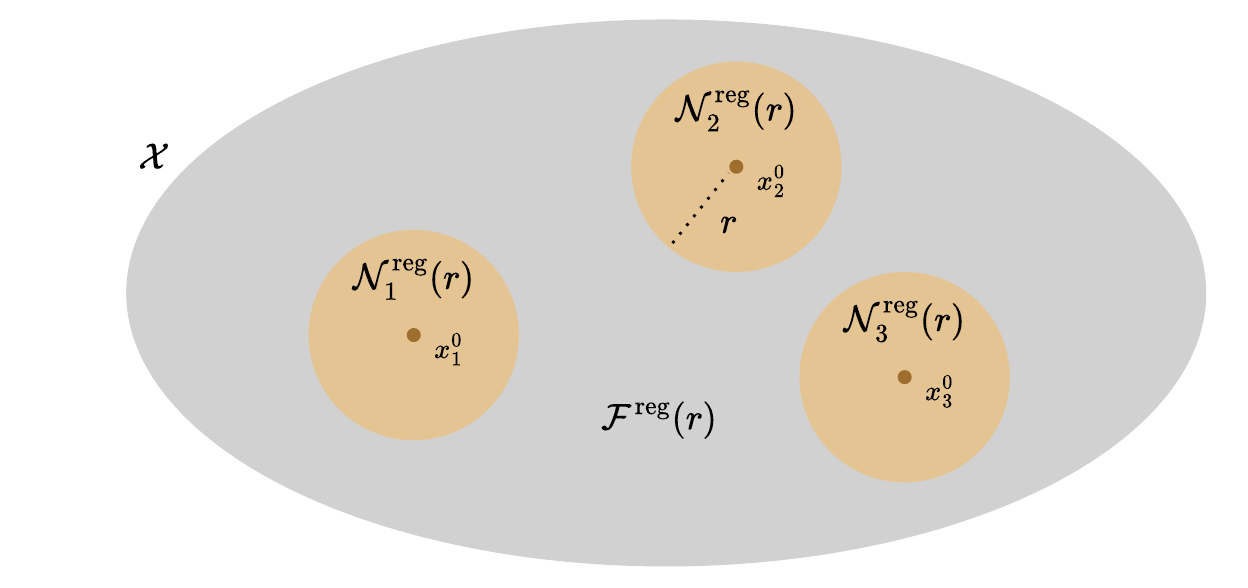}
    \caption[Two-dimensional illustration of near and far regions]{Two-dimensional illustration of near (light marroon) and far (light gray) regions with fixed radius $r$ and a distance $\distancegeneric(\bs,\bt) \propto \Vert \bs - \bt \Vert^2$.}
    \label{fig:near_far}
\end{figure}
\medskip

Existing bounds on estimation errors proved in the literature \parencite{azais2015spike,candes2014towards,poon2023geometry} depend on some fixed $r>0$, that may depend on $\Delta_0$, and are expressed as follows:
\begin{subequations}\label[pluralequation]{eqs:error_bounds}
\begin{itemize}
    \item \emph{Control of the far region:} 
    \begin{equation}
    \label{eq:control_far}
        |\estimator|\big(\farregion(r)\big)\lesssim_d \gamma \sqrt{\K}\,,
    \end{equation}
    \item \emph{Control of all the near regions:}
    \begin{equation}
    \label{eq:control_near}
        \big|\estimator\big(\nearregion_k(r)\big)-a^0_k \big|\lesssim_d \gamma \sqrt{\K}\,,
    \end{equation}
    \item \emph{Detection level: }For all borelian $A\subset\Param$ such that $|\estimator|(A)\gtrsim_d \gamma \sqrt{\K}$, there exists $\trueparam_k$ such that
    \begin{equation}
    \label{eq:detection_near}
    \distancegeneric(A,\trueparam_k)\defeq \min_{\bt\in A}\distancegeneric(t,\trueparam_k)\lesssim_d r\,,
    \end{equation}    
\end{itemize}
\end{subequations}
where $\lesssim_d$ denotes the inequality up to a multiplicative constant that may depend on the dimension $d$, $\estimator$ is a solution to \eqref{eq:blasso} with regularisation parameter $\kappa\simeq \gamma /\sqrt{\K}$ and $|\estimator|$ denotes the absolute part of $\estimator$. 

\medskip

\begin{remark}[On regularisation parameter calibration]
    \label{rem:calibration}
    The regularisation parameter $\kappa$ is crucial for the performance of BLASSO. In practice, it is often calibrated using cross-validation techniques or other model selection criteria. The choice of $\kappa$ can significantly impact the balance between bias and variance in the estimation process. In our analysis, we assume that $\kappa$ is chosen to be proportional to the noise level $\gamma$ and the sparsity $\K$, specifically $\kappa\simeq \gamma/\sqrt{\K}$. This choice ensures that the regularisation term effectively controls the trade-off between fitting the data and promoting sparsity in the estimated measure.

    When $\K$ is not known, other choices of the regularisation parameter can be considered, such as $\kappa\simeq\gamma$ which is a more realistic choice in practice, yet less optimal in terms of the theoretical guarantees. In this case, the error bounds \eqref{eqs:error_bounds} scale as $\K$ instead of $\sqrt{\K}$, which is a less favorable scaling, yet still provides useful information about the localisation of the solution, see Remark~\ref{rem:tuning_kappa_main_thm} below.
\end{remark}

\medskip

This error analysis of BLASSO \eqref{eq:blasso} relies on controlling both the noise level $\gamma$ and the Bregman divergence of the total variation norm
\begin{equation}
    \label{eq:def_bregman}
    \mathcal D_{\eta^0}(\mu\,||\,\target)\defeq\|\mu\|_{\mathrm{TV}}-\|\target\|_{\mathrm{TV}}-\langle\eta^0,\mu-\target\rangle_{\mathcal C(\Param)\times\mathcal M(\Param)}
\end{equation}
between the underlying $\target\in\Model_{\K, \Delta_0, \distancegeneric}$ 
and a solution $\mu=\estimator$ \parencite{decastro2012exact,azais2015spike} or, as we will see, any \textit{approximate} solution to \eqref{eq:blasso}. 

\medskip
\pagebreak[3]

\noindent
This error analysis requires that
\begin{itemize}
    \item $\eta^0$ is an element of the RKHS $\Hilbert_\model$,
    \item $\eta^0$ is a sub-gradient of the total variation norm at point~$\target$ with specific properties on the near regions~$\nearregion_k(r)$ and far region $\farregion(r)$,
\end{itemize}
Such $\eta^0$ is referred to as a \emph{non-degenerate dual certificate of radius $r$}~\parencite{decastro2012exact,candes2014towards,duval2015exact}. When such a certificate exists, the control on Bregman divergence directly translates into the estimation error bounds of~Equations \eqref{eqs:error_bounds}, yielding precise guarantees on both support recovery and parameter estimation accuracy \parencite{duval2017sparse}. 

\subsection{Construction of dual certificates under local positive curvature assumption (LPC)}

Crucially, a sufficient condition for BLASSO error analysis to hold is that the model kernel verifies a so-called \textit{local positive curvature} assumption~\parencite[$\LPC$,][Assumption 1]{poon2023geometry}, given in Section~\ref{sec:assumptions}. The $\LPC$ ensures the existence of a non-degenerate dual certificate $\eta^0$ of radius $r$ for any target measure $\target$ in the sparse model class~$\Model_{\K, \Delta_0, \mathfrak{d}_{\textnormal{FR}}}$ where $\Delta_0$ is taken sufficiently large with respect to $r$. There, the distance $\mathfrak{d}_{\textnormal{FR}}$ is the so-called \emph{Fisher-Rao distance} induced by the model kernel, which we will detail later on in \Cref{sec:assumptions}. 

To the best of our knowledge, the literature \parencite{azais2015spike, candes2014towards, poon2019support, poon2023geometry} has proved that the estimation error bounds of~\Cref{eqs:error_bounds} hold, but only for {\em fixed} near regions of size $r=r_0$ {\em irrespective of the noise level}. However, in view of~Equations \eqref{eqs:error_bounds}, the smaller the radius~$r$, the closer the BLASSO localizes around the true support, hence the sharper the error bounds. One of our key contributions is showing that \emph{these bounds extend to near regions whose size can adapt to the noise level~$\gamma$, decreasing as $\mathcal O(\sqrt{\gamma\sqrt{\K}})$}, which we call \emph{effective} near regions. This is a substantial improvement for small noise levels since the $\LPC$ is a property \textit{intrinsic} to the model kernel, with parameters independent of $\target$ and $\by$, and thus of the noise level $\gamma$. 

\subsection{Previous studies on kernels for BLASSO}

Fortunately, the existence of dual certificates has already been proved for some kernels. 
Some of these have also been shown to verify the $\LPC$ assumption. %
Pioneering work \parencite{candes2014towards} has proved that the non-degenerate dual certificate exist for the Jackson kernel:
\begin{subequations}
\begin{itemize}
    \item \textbf{Jackson kernel} \parencite{candes2014towards} on $\Param=[0,1)^d$ identified to the $d$-dimensional torus $(\R\setminus\mathbb{Z})^d$. The kernel is the fourth power of the Dirichlet kernel, i.e.
    \begin{equation}
    \label{eq:def_Jackson}
        K(\bs, \bt) = \prod_{i=1}^d \kappa(s_i - t_i), \quad\textnormal{with } \kappa(t) = \bigg(\frac{\sin\big(\big(\frac N2+1\big)\pi t\big) }{\big(\frac N2+1\big)\sin(\pi t)}\bigg)^4\textnormal{ and } N\geq 1\,.
    \end{equation}
    This kernel satisfies the $\LPC$ \parencite[Appendix C]{poon2023geometry}.
\end{itemize}

Using a similar technique as in \cite{candes2014towards}, a low-pass spherical harmonic kernel has been investigated in \cite{bendory2015exact} for \emph{spike detection} on the sphere. 
\begin{itemize}
    \item \textbf{Low-pass spherical harmonic} \parencite{bendory2015exact} on the $2$-sphere $\Param=\mathbb S^2\subset\R^3$, given by
    \begin{equation}
        K(\bs, \bt) = C(N) \sum_{k=0}^N \zeta(k/N) P_{k,3}(\langle\bs,\bt\rangle)\,,\notag
    \end{equation}
    where $N\geq 1$ is some degree, $C(N)$ is a normalizing constant, $P_{k,3}$ is the univariate ultraspherical Gegenbauer polynomial of order $3$ and degree $k$, and $\zeta(t)$ is a smooth non-negative univariate function equal to~$1$ for $t$ between $0$ and $1/2$, $0$ for $t$ greater than $1$ and less than $1$ otherwise. \\
    The $\LPC$ has not been proved for this kernel although the existence of a dual certificate has been proven via alternative techniques.
\end{itemize}

Recently, two others kernels have been proved to satisfy the $\LPC$, namely

\begin{itemize}
    \item \textbf{Gaussian kernel} \parencite[Appendix D]{poon2023geometry} given by
        \begin{equation}
            \label{eq:sketch_kernel_Gaussian}
            K(\bs, \bt) =\Theta_{\bm{\Omega}}(\bs-\bt)
                \quad
            \textnormal{where}
                \quad
                \Theta_{\bm{\Omega}}\defeq\exp{\big(-\frac12\big\|\cdot\big\|^2_{\Omega}\big)}
            \,,
        \end{equation}
and $\|\bt\|^2_{\bm{\Omega}}=\langle \bt,\bm{\Omega}^{-1} \bt\rangle$ for some positive definite matrix $\bm{\Omega}$;
\pagebreak[3]
\item \textbf{Laplace transform kernel} \parencite[Appendix E]{poon2023geometry}. This kernel is defined on~${\mathcal{X} = \mathbb{R}_+^d}$ and arises when studying continuous sampling of the Laplace transform of a positive signal, defined as~${\mathcal{L}[\mu](\bomega) = \int_\mathcal{X} e^{-\bomega^\top \bt} \diff \mu(\bt)}$. It was studied in~\textcite{denoyelle2019sliding, poon2019support} and is not to be confused with the \textit{Laplace kernel}. In particular, it is not translation-invariant and is given by
\begin{equation}
    \label{eq:sketch_kernel_laplace}
    K_\alpha(\bs, \bt) = \prod_{i=1}^d \kappa(s_i + \alpha_i, t_i+\alpha_i), \quad\textnormal{with } \kappa(a,b) = \frac{2 \sqrt{ab}}{a+b} \textnormal{ and } \alpha_i > 0.
\end{equation}
\end{itemize}

Importantly, these kernels (except the low-pass spherical harmonic) will serve as "pivot kernels" in our kernel switch analysis—a central concept to our contribution. As we will show, {\em a pivot kernel that satisfies the $\LPC$ enables the construction of non-degenerate dual certificates} 
with {\em other} model kernels. 
\begin{table}[!ht]
    \centering
    \renewcommand{\arraystretch}{1.5} %
    \resizebox{\textwidth}{!}{
    \begin{tabular}{@{}lllll@{}}
        \toprule
        Kernel name & Domain $\mathcal{X}$ & $K(\bs, \bt)$  &  \Cref{hyp:local_positive_curvature} ($\LPC$) & $\LPC$-parameters  \\ %
        \midrule
        Jackson & $(\mathbb{R} \backslash \mathbb{Z})^d$ & \eqref{eq:def_Jackson}  & \textcite[Appendix C]{poon2023geometry} &  $r_0 = \mathcal{O}\big(1\big)$, $\Delta_0 = \mathcal{O}\big(d^{\frac12} s_0^{\frac14}\big)$ \\
        Gaussian & $\mathbb{R}^d$ & \eqref{eq:sketch_kernel_Gaussian}  & \textcite[Appendix D]{poon2023geometry} & $r_0 = \mathcal{O}\big(1\big)$, $\Delta_0 = \mathcal{O}\big( \log(\K) \big)$ \\
        Laplace & $\mathbb{R}^d$ & \eqref{eq:sketch_kernel_laplace}  & \textcite[Appendix E]{poon2023geometry} & $r_0 = \mathcal{O}\big(1\big)$, $\Delta_0 = \mathcal{O}\big(d + \log(d^{\frac{3}{2}} s_0) \big)$  \\
        {\em Sinc-4} & $\mathbb{R}^d$ & \eqref{eq:sketch_kernel_sinc}  & {\em This work} - \Cref{thm:sinc4_lca} &  $r_0 = \mathcal{O}\big(\frac{1}{d}\big)$, $\Delta_0 = \mathcal{O}\big(d^{\frac74} s_0^{\frac14}\big)$ \\
        \bottomrule
    \end{tabular}
    }
    \vspace{5pt}
    \caption[Table of known LPC-kernels.]{Table of %
    known kernels verifying the $\LPC$ assumption, with associated references for proofs. The parameters $s_0$ and $r_0$ are the sparsity and radius of the near regions, respectively. The parameter $\Delta_0$ is the minimal separation between two spikes in the target measure.}
    \label{tab:known_pivot_kernels}
\end{table}

\subsection{The sinc-4 kernel}
Another important kernel introduced by \textcite[Section 4.1]{decastro2019sparse} to derive recovery guarantees for mixture models is the fourth power of the sinus cardinal kernel, coined sinc-4.

\begin{itemize}
    \item  \textbf{Sinc-4 kernel} where the $\sinc^4$ function is defined on $\mathbb{R}^d$ as the product over coordinates of the fourth power of the univariate sinus cardinal $\sin(u) /u$, and 
    \begin{equation}
    \label{eq:sketch_kernel_sinc}
    K(\bs, \bt) =\Psi_\bandwidth(\bt-\bs)
    \quad
    \textnormal{where}
    \quad
    \Psi_\bandwidth(\bu) \defeq\sinc^4\Big(\frac\bu{4\bandwidth}\Big)
    \,,\ \bandwidth>0\,.
    \end{equation}
     
\end{itemize}
\end{subequations}
This kernel is particularly interesting as pivot, since it is bandlimited so that the kernel switch analysis is valid with any translation invariant kernel with large enough spectral support.
However, the analysis of~\textcite[Section 4.1]{decastro2019sparse} does not allow to directly prove the~$\LPC$ for this kernel. One of the important contributions of our paper is~\Cref{thm:sinc4_lca} which proves that this kernel indeed satisfies the $\LPC$, yielding BLASSO guarantees for a wide range of translation-invariant kernels .

\subsection{Continuous sparse regression and sketching}
\label{sec:subsec_sketched_BLASSO}
A key motivation for the introduction of pivot kernels is to address model kernels arising from sketched problems, which we now detail.

The general goal is to estimate a sparse target measure $\target\in\Model_{s, \Delta, \distancegeneric}$ 
from a sketch $\bsketch\in\mathbb{C}^m$, by solving the corresponding BLASSO program~\eqref{eq:blasso}, with the sketched forward operator $F_{\textnormal{sketch}}$, which definition involves the sketching law~$\SketchDist$. Thus the corresponding setting is as follows:
\begin{subequations}
\label[pluralequation]{eqs:def_sketch_blasso}
\begin{itemize}
    \item \textit{Hilbert space} $\mathcal{F}_{\textnormal{sketch}}$: consider the standard dot product on $\mathbb{C}^m$, $\mathcal{F}_{\textnormal{sketch}}=\mathbb{C}^m$,
    \item \textit{Forward measurement} $F_{\textnormal{sketch}}: \,\mathcal{M}(\Param)\to \mathbb{C}^m$ given by
    \begin{equation}
    \label{eq:sketch_forward}
        (F_{\textnormal{sketch}}\mu)_i= \frac{1}{\sqrt{m}}\int_\Param {\varphi}_{\bomega_i}(\bt) \mathrm{d}\mu(\bt)\,,\quad i=1,\ldots,m\,,
    \end{equation}
    where $\varphi_\bomega\,:\,\Param\to\mathbb C$ is a complex-valued \underline{continuous} function, referred to as the {\em sketching function}.
    \item \textit{Sketching law} $\SketchDist$: the sequence of \textit{i.i.d.\!\!} random vectors $(\bomega_i)_i$ is drawn with respect to $\SketchDist$,
        \begin{equation}
        \label{eq:sketch_distribution}
        \bomega_i\underset{\mathrm{i.i.d.}}\sim\SketchDist\,,\quad i=1,\ldots,m\,.
    \end{equation}
\end{itemize}
\end{subequations}

\pagebreak[3]

\begin{subequations}
Thus, in line with \Cref{def:model_kernel}, we define the \textit{sketched} model kernel by
\begin{equation}
    \label{eq:sketch_model_kernel}
    {K}_{\textnormal{sketch},\model}(\bs, \bt)
    \defeq \Big\< F_{\textnormal{sketch}} \delta_\bs, F_{\textnormal{sketch}}\delta_\bt\Big\>_{\mathbb{C}^m}  =
    \frac{1}{m}\sum_{i=1}^m \varphi_{\bomega_i}(\bs)\overline{\varphi}_{\bomega_i}(\bt)\,,
\end{equation}
where $\overline{\varphi}_\bomega$ denotes complex conjugation. This kernel depends on the draw of the sketch, and taking its expectation naturally leads to the \textit{population} model kernel
\begin{equation}
    \label{eq:population_model_kernel}
    {K}_{\model}(\bs, \bt) \defeq \Expectation_{\bomega\sim \SketchDist}\big[\varphi_\bomega(\bs)\overline{\varphi}_\bomega(\bt)\big] \,.
\end{equation}

\end{subequations}
Note that every draw of the sketch model kernel meets Definition~\ref{def:model_kernel} and satisfies the continuity assumption~\eqref{hyp:continuous_model_kernel} with the construction~\eqref{eqs:def_sketch_blasso}. Furthermore, its limit as~$m$ tends to infinity is almost surely the population model kernel, for fixed $\bs$ and $\bt$. 

An interesting feature of the sketched problem is to reduce the computational cost induced by the model kernel by considering an empirical version of rank $m$. Moreover, \textcite{poon2023geometry} showed that the error bound analysis of the BLASSO also holds for its sketched counterpart, provided that $K_\model$ satisfies the $\LPC$ assumption, under additional probabilistic assumptions %
on the tails of the feature function $\varphi_\bomega$ and its derivatives, which are detailed later in~\Cref{hyp:sketch_function_bounded}.

However, proving the $\LPC$ assumption has to be done on a case-by-case basis, possibly leading to tedious calculations. Instead, we propose in this paper to \textit{switch} to another kernel (known to satisfy the $\LPC$) using an alternative (continuous) sketching function $\psi_\bomega\,:\,\Param\to\mathbb C$. We define the associated kernel by
\begin{equation}
    \label{eq:pivot_sketch_kernel}
    K(\bs, \bt)\defeq \Expectation_{\bomega\sim \SketchDist}\big[\psi_\bomega(\bs)\overline{\psi}_\bomega(\bt)\big]\,,
    \addtocounter{equation}{1}
\tag{\theequation--$\mathbf{H}_{\psi_\bomega,\SketchDist}$}
\end{equation}
where we assume that $K(\cdot,\cdot)$ is real valued and the expectation is taken with respect to the same distribution~$\SketchDist$. Given a model kernel $K_\model$, we characterize a broader class of kernels $K$ that can be used as \textit{pivot} in the BLASSO analysis, provided that $\Hilbert_\model$ is continuously included in their RKHS (see~\Cref{fig:kernel_switch}). Our result offers an alternative route to obtain the controls of~\Cref{eqs:error_bounds} when considering a specific BLASSO problem.

\section{Sketching mixture models via continuous sparse regression}
\label{sec:Mixtures}
This section focuses on the application of our general theoretical framework to the problem of mixture model estimation. We present specific consequences of the broader results that will be detailed in Section~\ref{sec:BlassoResults}. The key novelties for mixture modelling include the ability of our proposed algorithm, S2Mix, to handle a wide variety of template densities beyond the commonly assumed Gaussian. Furthermore, we establish new localisation error bounds for parameter recovery, demonstrating improved scaling as the sample size increases. These advancements significantly enhance the practical utility and theoretical understanding of continuous sparse regularisation for mixture models.
\subsection{A resolution-based approach of mixture modelling}
\label{sec:smoothing}
We demonstrate our theoretical contributions through their application to the important problem of mixture model estimation. Mixture models represent fundamental statistical tools used across machine learning, signal processing, and computational biology to model heterogeneous populations. Before introducing the technical details, let us sketch the key ideas: we reformulate mixture estimation as a continuous sparse regression problem where the unknown mixture components—both their weights and locations—are encoded as a sparse discrete measure. This measure can then be recovered using BLASSO, even when working with compressed data through sketching techniques.

In this section, we first present a resolution-based approach to mixture modelling via continuous sparse regression. We then develop an efficient sketching algorithm, coined S2Mix, that enables parameter estimation from large datasets with near-optimal sample complexity. Crucially, we establish novel statistical guarantees for this approach by leveraging the sinc-4 kernel $\Psi_\bandwidth$ as a pivot kernel within our theoretical framework. This demonstrates how our effective near regions analysis and kernel switch technique deliver practical improvements for an important class of statistical models.

We are interested in estimating the mixture positive weights $(a_k^0)_k$ and parameters $(\trueparam_k)_k$ from a $n$-sample $\Obs:=\{\obs_1,\ldots,\obs_n\}$ of a translation invariant mixture model with $\K$-components:
\begin{align}
\label{eq:MixtureDensity}
\targetdensity \defeq \sum_{k=1}^{\K} a_k^0 \templatedensity(\cdot - \trueparam_k) 
\quad \textnormal{with}\quad
\sum_{k} a_k^0 = 1
\textnormal{ and }
\trueparam_k\in\Param\,,
\end{align}
where $\Param$ is a  compact %
subset of $\R^d$. 

\pagebreak[3]

Here, $\templatedensity$ is a probability density function on $\R^d$ called the \textit{template distribution}. Popular instances include mixtures of Gaussians with fixed and known covariance, or mixtures of Laplace distributions. More generally, our proposed framework will also encompass mixtures of any distribution with known characteristic function, such as mixtures of stable distributions for instance.

\medskip

From a non-parametric statistics perspective, estimating a mixture model amounts to estimate the density~$\targetdensity$, which can be done via kernel density estimation ({\it i.e.,} smoothing the empirical sample distribution), namely
\begin{equation*}
	L_\bandwidth \hat{f}_n \defeq \kernel_\bandwidth \star \hat{f}_n
	\quad\textnormal{where}\quad
	\hat{f}_n \defeq \frac1n \sum_{j=1}^n \delta_{\obs_j}\,.
\end{equation*} 
The \textit{smoothing} function $\kernel_\bandwidth(\cdot) = \bandwidth^{-d} \lambda(\cdot / \bandwidth)$ is defined by the smoothing kernel $\lambda$ and its \emph{bandwidth} parameter~$\bandwidth$. In this work, we focus on the sinus cardinal smoothing kernel which is a standard choice in the literature \parencite{decastro2019sparse}, given by  
\begin{equation}
    \lambda(\bu) \defeq \sinc(\bu) = \prod_{i=1}^d \frac{\sin(u_i)}{u_i} \quad \textnormal{so that} \quad \lambda_\bandwidth(\bx) = \bandwidth^{-d} \sinc\Big(\frac{\bx}{\bandwidth}\Big). \label{eq:DefSinc}
\end{equation} 
This function acts as a low-pass filter in the Fourier domain, other filters are possible but we will consider the sinus cardinal function in this section. The crux of non-parametric methods is then to optimise the mean quadratic loss of the estimator with respect to the bandwidth parameter~$\bandwidth_n$, usually depending on $n$, the latter translating into a bias-variance trade-off.

In this work, we follow a different path considering that the bandwidth parameter $\bandwidth$ will be fixed by the practitioner and will not depend on the sample size~$n$ but rather on the \emph{minimal separation} between parameters $(\trueparam_k)_k$, also referred to as \textit{Rayleigh limit}. Precisely, we assume that the bandwidth $\bandwidth>0$ is such that 
\begin{equation}
        \addtocounter{equation}{1}
	\tag{\theequation--$\mathbf{H}_\bandwidth$}
    \bandwidth \leq \bandwidth_{\max}(\target) \defeq \frac{\min_{1\leq k\neq \ell\leq \K} \big\Vert \trueparam_k - \trueparam_\ell \big\Vert_2}{147.77\,\K^{1/4}d^{7/4}}\,.
	\label{hyp:Hbandwidth}
\end{equation}
As we will see, this hypothesis ensures that the bandwidth is small enough to resolve the mixture components.

\medskip

Now, consider a putative parameter $\theta$ given by 
\[
	\theta\defeq\big\{s;\, a_1,\ldots, a_s;\, x_1,\ldots x_s\big\}\quad
	\textnormal{and}\quad
	\mu_\theta\defeq\sum_{k=1}^s a_k\delta_{\bx_k}\,,
\]
with $s\geq 1$ the number of mixing components and $\mu_\theta$ is the lifting of the putative parameter on the space of measures, where~$\delta_\bt$ denotes the Dirac mass at location $t\in\Param$. 

The mixture model density defined in \Cref{eq:MixtureDensity} can be seen as an image of a discrete measure through the convolution operator: 
\begin{equation}
\begin{aligned}
\Phi : & & \mu \in \Measures(\Param) & & \mapsto & & & \Phi \mu \defeq \templatedensity \star \mu = \int \templatedensity(\cdot - \param) \diff \mu(\param)
=\int\templatedensity_{\param}\diff \mu\,,
\end{aligned}
\end{equation}
where $\templatedensity_{\param}\defeq\templatedensity(\cdot - \param)$. This mapping is reminiscent of the so-called Kernel Mean Embedding of $\mu$ in Machine Learning, and they exactly corresponds when the template $\templatedensity$ is of positive type~\parencite{muandet2017kernel}.

\pagebreak[3]

The mixture problem seeks to fit a parameter $\theta$ from the observation of $\hat{f}_n$, namely
\[
    \theta \xrightarrow{\textnormal{Non-linear}} \Bigg[f_{\theta} \defeq \Phi\mu_\theta  = \sum_{k=1}^{s} a_k \templatedensity(\cdot - t_k) \Bigg] \xrightarrow{\textnormal{Sampling}} \textcolor{NavyBlue}{\hat{f}_n}\,.
\]
where the step $\theta \to f_{\theta}$ is a non-linear mapping of the parameters to the mixture density. Directly solving the corresponding inverse problem is hard, both because of the non-linearity as well as the fact that those two objects do not belong to the same functional space. On the other hand, lifting the problem on the space of measures yields the following linear problem

\begin{equation}
    \mu_{\theta}\xrightarrow{\textnormal{Linear}} \Big[f_{\theta}=\Phi\mu_\theta \Big]\xrightarrow{\textnormal{Embedding}} \Big[L_\bandwidth f_{\theta}=(L_\bandwidth \circ\Phi)\mu_\theta \Big] \approx \textcolor{NavyBlue}{L_\bandwidth \hat{f}_n \xdashleftarrow{\text{Embedding}} \hat{f}_n},
\end{equation}
where we have embedded the candidate $L_\bandwidth f_{\theta}$ and the observation~$L_\bandwidth \hat{f}_n$ into the same Hilbert space $\mathcal{F}_\bandwidth$, the RKHS of $\lambda_\bandwidth$.

\subsection{Continuous sparse regularisation for mixture estimation -- Supermix}

Recent advances have established a powerful connection between continuous sparse regularisation techniques developed for super-resolution problems and the statistical challenge of mixture estimation \parencite{decastro2019sparse, poon2023geometry}. 
Indeed, as we have just seen, mixture modelling can be reformulated as the problem of recovering an unknown sparse discrete measure $\target$ from observations obtained through a linear measurement operator~$\Forward = L_\bandwidth \circ\Phi$ mapping from the space of finite Radon measures $\Measures(\Param)$ to a Hilbert space~$\mathcal{F}$, a possible approach to address it is to seek a (near) optimum of the BLASSO problem \eqref{eq:blasso}. 

This formulation offers the significant advantage that the number of components $\K$ need not be specified in advance. Speciminimalfically, this approach involves
\begin{subequations}
\label[pluralequation]{eqs:SketchedSupermix}
\begin{itemize}
    \item \textit{Hilbert space}: $\mathcal{F} \defeq \mathcal{F}_\bandwidth$, the RKHS associated to the smoothing kernel $\lambda_\bandwidth(t-s)$ defined in~\eqref{eq:DefSinc}. As a translation invariant reproducing kernel, it admits an explicit characterisation in the Fourier domain (see \Cref{app:RKHS}).
    \item \textit{Forward measurement operator} $F \defeq L_\bandwidth \circ \Phi: \mathcal{M}(\Param)\to \mathcal{F}_\bandwidth$.
    \item \textit{Observation} $\by \defeq L_\bandwidth \hat{f}_n$ which is of the form 
    \[
    	\by = F \target + \Gamma_n \quad \textnormal{with} \quad \Gamma_n:=L_\bandwidth \hat{f}_n-\Expectation[L_\bandwidth \hat{f}_n],
    \]
    the latter being a centered noise term due to $n$-sampling with $\Expectation[L_\bandwidth \hat{f}_n]=L_\bandwidth f^0=F\target $. 
\end{itemize}
\end{subequations}
Hence, one can use a sparse regression estimation procedure to fit the mixture model, which was coined \textit{Supermix} by \textcite{decastro2019sparse}. %
For the specific problem at hand, the model kernel defined in \eqref{def:kmod_general} is translation invariant and can be described in the Fourier domain (see~\Cref{app:RKHS}). Indeed, one has 
\begin{align*}
    K_\model(\bs,\bt) 
    &\defeq\langle F\delta_\bs,F\delta_\bt \rangle_{\mathcal{F}} = \langle \kernel_\bandwidth \star \templatedensity_\bs, \kernel_\bandwidth \star \templatedensity_\bt  \rangle_{\mathcal{F}_\bandwidth}, \\
    &= \frac{1}{(2 \pi)^d} \int \fourier[\lambda_\bandwidth](\bomega) \cdot \lvert \templatefourier \rvert^2 (\bomega) \cdot e^{+ \imath \bomega^\top (\bt - \bs)} \diff \bomega, \\
    &= \fourier^{-1} [\fourier[\lambda_\bandwidth] \cdot \vert \templatefourier \vert^2](\bt - \bs), \\
    &= (\lambda_\bandwidth \star \templatedensity \star \check{\templatedensity})(\bt - \bs)\,,
\end{align*}
where $\check{\templatedensity}(\bx)=\templatedensity(-\bx)$. At first glance, the natural road to derive estimation error bounds is trying to prove $\LPC$ (Assumption \ref{hyp:local_positive_curvature}) for this translation invariant kernel. However, this %
would require controlling non-trivial quantities
involving the model kernel and its associated Fisher-Rao metric. %
Alternatively, \textcite{decastro2019sparse} proposed to use the sinc-4 kernel $K_\pivot(\bs, \bt) = \Psi_\bandwidth(\bt - \bs)$ defined in~\Cref{eq:sketch_kernel_sinc} as a proxy to build non-degenerate certificates, implicitly relying on a \textit{kernel switch} analysis. The latter holds under the assumption that the RKHS $\Hilbert_{\pivot}$, defined by $K_{\pivot}$, is continuously included in $\Hilbert_{\model}$, defined by $K_{\model}$. In the special case where both $K_{\pivot}$ and $K_{\model}$ are translation invariant, this translates into an assumption on the boundedness of the ratio of their Fourier transforms~(see \textit{e.g.}~\textcite[Corrollary~3.2]{zhang2013inclusion} and \Cref{prop:sinc4_cswitch}): 
\begin{equation}
       C_{\textnormal{switch}} = C_{\textnormal{switch}}(\bandwidth, \templatedensity) \defeq \sup_{\bomega \in [-1/\bandwidth, 1/\bandwidth]^d}\sqrt{\frac{\fourier[\Psi_\bandwidth]}{ \fourier[\lambda_\bandwidth] \big\lvert \fourier[\templatedensity] \big\rvert^2} (\bomega)} < + \infty,
        \tag{\theequation--$\mathbf{H}_\templatedensity$}
        \label{hyp:Hphi}
\end{equation}
with the convention $\frac{0}{0} = 0$. This is a special instance of the general results discussed in \Cref{sec:BlassoResults}. The scaling of this constant $C_{\switch}$ with respect to the bandwidth parameter $\bandwidth$ depends on the smoothness of $\fourier[\templatedensity]$ as discussed in \textcite[][,Section 5.2]{decastro2019sparse}. In particular, the case of \emph{supersmooth} densities defined hereafter leads to a scaling in
\[
\exists p \in [1, +\infty], \; \alpha,\beta >0, \quad \fourier[\templatedensity](\bomega) \propto e^{-\alpha \Vert \bomega \Vert_p^\beta} \implies C_\switch = \mathcal{O}_d\Big(\tau^{d/2} e^{\alpha \big(\tfrac{d^{1/p}}{\bandwidth}\big)^\beta}\Big)\,,
\]   
where the constant in $\mathcal{O}_d$ may depend exponentially in the dimension $d$, but it does not depend on $\tau$. This case encompasses \textit{e.g.} Gaussian, multivariate Cauchy, product of univariate Cauchy, or more generally of centered stable distributions with known scale parameter and zero skewness. As an example, when $\templatedensity$ is a Gaussian density, the scaling is given by~$C_\switch = \mathcal{O}_d(\tau^{d/2} e^{d /2\bandwidth^2)})$. 

In addition, the scaling of the noise level $\gamma_n \geq \|\Gamma_n\|_{\mathcal{F}_\bandwidth}$ defined in \Cref{eqs:noise_level_term} can be controlled with high probability over the draw of the sample $\Obs$ given by the following Lemma
\begin{lemma}[Control of the noise level $\gamma_n$]
    \label{prop:control_population_noise_level_whp}
    Let $\alpha > 0$. Then with probability at least $1 - \alpha$, it holds that 
    \begin{equation}
        \|\Gamma_n\|_{\mathcal{F}_\bandwidth}\leq\gamma_n := C_\alpha %
        {\frac{{\bandwidth}^{-\frac{d}2}}{\sqrt n}}\,,
    \end{equation}
    where $C_\alpha \defeq 2\sqrt{1 + C_1 \log({C_2}/{\alpha})}$ is a constant only depending on $\alpha$, while $C_1$ and $C_2$ are universal constants.
\end{lemma}
\begin{proof}
    The result is proved in \textcite[][, Lemma 16]{decastro2019sparse} using concentration of U-process and \textcite[][, Proposition 2.3]{arcones1993limit}. A closely related proof is given for the sketched case in Lemma~\ref{lem:noise_level_sketched}, see also Remark~\ref{rem:link_noise_level_concentration_population_sketched}.
\end{proof}
Thus, we can tune the %
regularisation $\kappa$ of the BLASSO problem \eqref{eq:blasso}, which is crucial for the statistical recovery result. This leads to the following proposition, which is an instance of the upcoming \Cref{thm:main_thm} for the Supermix problem, using the sinc-4 kernel $\Psi_\bandwidth$ as pivot. It can also be seen as a refined version of \textcite[][, Theorem 10]{decastro2019sparse} as commented below.
\begin{proposition}[Statistical guarantees for the Supermix problem]
    \label{prop:supermix_thm10}
    Let $\Param\subset\R^d$ be a compact set, and the target $\target\in \Measures(\Param)$ be a sparse measure of $s_0$ spikes $\{\trueparam_1, \ldots, \trueparam_{s_0}\}$, with minimal separation~${\Delta = \min_{k\neq l} \Vert \trueparam_k - \trueparam_l \Vert^2}$. Fix the bandwidth parameter $\bandwidth \leq \bandwidth_{\max}(\target)$ as in \Cref{hyp:Hbandwidth}. Take any template distribution $\templatedensity$ such that it achieves a finite \emph{kernel switch} constant $C_\switch(\bandwidth, \templatedensity) < + \infty$ in~\Cref{hyp:Hphi}. Consider the forward measurement operator~${F : \mu \in \Measures(\Param) \to  \lambda_\bandwidth \star \templatedensity \star \mu \in \mathcal{F}_\bandwidth}$, where $\lambda_\bandwidth$ is the smoothing kernel defined in \Cref{eq:DefSinc}. Let $\Obs = \{\obs_1, \ldots, \obs_n\}$ be an i.i.d. sample from $\templatedensity \star \target$, and $\empiricaldist$ its associated empirical distribution. Define the following BLASSO problem, named Supermix, with observation $\by = \lambda_\bandwidth \star \empiricaldist$
    \begin{equation}
        \min_{\mu \in \Measures(\Param)}\left\{ J_\kappa(\mu) = \frac{1}{2} \Vert \by - F\mu \Vert_{\mathcal{F}_\bandwidth}^2 + \kappa \Vert \mu \Vert_{\mathrm{TV}} \right\}. 
        \addtocounter{equation}{1}
	\tag{\theequation--Supermix}
    \end{equation} 
    Take any diverging sequence $\delta_n$ and set the near and far regions' \emph{effective radius} $r_n = \delta_n n^{-1/4}$. Fix the regularisation~${\kappa = \frac{1}{\sqrt{2}C_\switch}  \frac{C_\alpha {\bandwidth}^{-\frac{d}2}}{\sqrt{n}}} \frac{1}{\sqrt\K}$. Then, for any near optimal $\mu \in \Measures(\Param)$ such that $J_\kappa(\mu) \leq J_\kappa(\target)$, $n \geq c_d \delta_n^4$ where $c_d$ depends polynomially on $d$, and $\alpha >0$, we have with probability $1 - \alpha$ over the drawing of~$\Obs$:
    \begin{subequations} 
    \begin{itemize}
        \item Control of the far region: 
        \begin{equation}
            |\mu|\big(\farregion(r)\big)\leq \frac{256 \sqrt{2}\,C_\alpha C_\switch {\bandwidth}^{-\frac{d}2}}{23\, \delta_n^2} \sqrt\K,
        \end{equation}
        \item Control of all the near regions: for all $k\in[\K]$,
        \begin{equation}
        \label{eq:control_near_mixture}
            \Big|\mu\big(\nearregion_k(r)\big)-a^0_k\Big|\leq \frac{1536\sqrt{2}\,C_\alpha C_\switch {\bandwidth}^{-\frac{d}2}}{23\, \delta_n^2} \sqrt\K
        \end{equation}
        \item Detection level: for all Borelian $A\subset\Param$ such that $|\mu|(A)> \frac{256\sqrt{2}C_\alpha C_\switch {\bandwidth}^{-\frac{d}2}}{23 \delta_n^2} \sqrt\K$, there exists $\trueparam_k$ such~that
        \begin{equation}
        \label{eq:detection_near_mixture}
            \min_{\bt\in A} \frac{1}{2 \sqrt{3} \bandwidth} \Vert\bt - \trueparam_k \Vert_2  \leq  r_n\,,
        \end{equation} 
    \end{itemize}

    \end{subequations}
\end{proposition}
The proof is given in \Cref{proof:supermix_thm10}, and is an application of the main theorem of \Cref{sec:BlassoResults},~\Cref{thm:main_thm}. This can be adjusted with regularisation $\kappa$ independent of $s_0$, with the price of slightly degraded error bounds as discussed in~\Cref{sec:BlassoResults} after~\Cref{thm:main_thm}.  Moreover, this result is a digest and improvement of~\textcite[][, Theorem~10]{decastro2019sparse} as made clear by the following remarks. 

\begin{remark}[Notational comparison with SuperMix]
    When comparing with \textcite{decastro2019sparse}, the following notational equivalences are relevant:
    \begin{itemize}
        \item The expression $1/(4\bandwidth)$, related to inverse bandwidth (or spatial resolution) in this paper, corresponds to their notation $m$.
        \item Our noise level, $\gamma_n \geq \| \Gamma_n \|$, corresponds to their notation $\rho_n$.
        \item Our constant $C_{\switch}(\bandwidth, \templatedensity)$ from \Cref{hyp:Hphi} corresponds to their notation $C_m(\varphi, \lambda)$.
        \item Their radius is denoted as $\epsilon$ and does not depend on $n$ (see \Cref{rem:effective_near_region} below).
    \end{itemize}
\end{remark}

\pagebreak[3]

\begin{remark}[Model set with sinc-4 pivot]
    The sinc-4 pivot kernel $\Psi_\bandwidth$ induces the so-called Fisher-Rao distance, properly defined in \Cref{eq:Fisher_Rao_distance}, which is here just a rescaled version of the Euclidean distance by a factor~${2\sqrt{3}\bandwidth}$.
    Consequently, for a measure to belong to the model class $\Model_{\K, \Delta_0, \distance}$, its support points $\{\trueparam_k\}$ must satisfy~${\|\trueparam_k - \trueparam_l\|_2 \geq 2\sqrt{3}\bandwidth \Delta_0}$.
    The parameter $\Delta_0$ is the minimal separation inherent to the Local Positive Curvature $(\LPC)$ property of $\Psi_\bandwidth$, which is established for this kernel in Theorem~\ref{thm:sinc4_lca}.
    Assumption \eqref{hyp:Hbandwidth} ensures that the chosen bandwidth $\bandwidth$ aligns the actual separations in the target measure $\target$ with this $\LPC$-derived requirement for $\Psi_\bandwidth$. Since the bandwidth $\bandwidth$ is chosen by the practitioner, this assumption means that they must chose $\bandwidth$ small enough to resolve the mixture components, which is a standard assumption in the literature \parencite{decastro2019sparse}.
\end{remark}

\medskip

\begin{remark}[Effective near regions]
    \label{rem:effective_near_region}
    One of the novelties of \Cref{prop:supermix_thm10} is that the radius of the near regions~${r_n = \delta_n n^{-1/4}}$ now depends on the number of samples. Moreover, any sequence $\delta_n$ diverging to infinity is acceptable provided that $r_n \to 0$. The latter controls a trade-off between how fast effective near regions narrow around the true support of $\target$ at a rate $r_n = \delta_n n^{-1/4}$, and how slow the estimator $\estimator$ detects the effective near and far regions at a rate $v_n = \delta_n^{-2}$.

    Thus, one have different \textit{scalings} of the near regions given by the choice of the diverging sequence $\delta_n$
    \begin{itemize}
        \item Polynomial: $\delta_n = n^{\alpha}$ where $0 < \alpha < 1/4$. The effective near regions scales in $n^{-1/4 + \alpha}$ and $v_n = n^{-2 \alpha}$.
        \item Logarithmic: $\delta_n = \sqrt{\log n}$. In this case the effective near regions narrows at the rate $(\log(n) / n)^{1/4}$ while~$v_n = \log(n)$.
        \item Log-Polynomial: $\delta_n = n^{\alpha}\log^\beta n $ where $0 < \alpha < 1/4$ and $\beta\geq0$. The effective near regions scales in~$n^{-1/4 + \alpha}\log^\beta n$ and~$v_n = n^{-2 \alpha}\log^{-2\beta}n$.
    \end{itemize}
\end{remark}

\subsection{Statistical error bounds for {\em sketched} mixture models}
While the preceding results are interesting from a theoretical perspective, a major practical difficulty arise when trying to solve the Supermix as the objective function $J_\kappa$, or its gradient, involves $d$-dimensional integrals. In this work, we focus on a \textit{sketched} version of the Supermix BLASSO problem leading to computational gains. 

\paragraph{Random Fourier features} The main idea revolves on an approximation of the translation invariant model kernel $K_{\model}$ using weighted random Fourier features \parencite{rahimi2007random, gribonval2020statistical}. The latter stems from the representation of a translation invariant $K_\model$ as an expectation over a free \textit{sketch} distribution, which is then approximated by $m$ Monte-Carlo samples.

Akin to importance sampling in the Fourier domain, this \textit{sketching distribution} $\SketchDist$ may be chosen different from (a properly normalized version of) the spectral measure~$U_\bandwidth \defeq \fourier[\lambda_\bandwidth] / (2 \pi)^d$ of the smoothing kernel. Formally, let $\SketchDist$ be \emph{any p.d.f. with support containing %
the hypercube~${\Vert \bomega \Vert_\infty \leq 1 / \bandwidth}$}, which is the support of $U_\bandwidth$. Then, the model kernel is amenable to sketching as it can be written in the form of~\Cref{eq:population_model_kernel} as follows:
\begin{subequations}
\label[pluralequation]{eqs:RFF_mixture}
\begin{align}
    K_\model(\bt,\bs) &=  \int \overset{U_\bandwidth}{\overbrace{\frac{\fourier[\lambda_\bandwidth]}{(2 \pi)^d}}}(\bomega) \lvert \templatefourier \rvert^2 (\bomega) e^{+ \imath \bomega^\top (\bt - \bs)} \diff \bomega,\notag \\
    &= \int \SketchDist(\bomega) \frac{U_\bandwidth}{\SketchDist}(\bomega) \lvert \templatefourier \rvert^2 (\bomega) e^{+ \imath \bomega^\top (\bt - \bs)} \diff \bomega,\notag \\
    &=\mathbb E_{\bomega \sim \SketchDist} \left[ \varphi_\bomega(\bs) \overline{\varphi_\bomega(\bt)} \right],\label{eq:kmod_as_expected_sketch_functions}\\
\text{where }\varphi_\bomega(\bt) &:= W(\bomega) \,\fourier[\templatedensity](\bomega)\, e^{-\imath \bomega^\top \bt}, \quad \textnormal{with: } W(\bomega) = \sqrt{\frac{U_\bandwidth}{\SketchDist}(\bomega)}.
\label{eq:sketch_functions_mixture}
\end{align}
Then, this leads to a Monte-Carlo approximation of the kernel, defining the \textit{sketched} model kernel as in \eqref{eq:sketch_model_kernel}
\begin{equation}
    \label{eq:mixture_sketch_model_kernel}
    {K}_{\textnormal{sketch},\model}(\bs,\bt)
    \defeq \frac{1}{m}\sum_{i=1}^m \varphi_{\bomega_i}(\bs)\overline{\varphi}_{\bomega_i}(\bt)\,,
\end{equation}
\end{subequations}
where the sequence of \textit{i.i.d.\!\!} random vectors $(\bomega_i)_i$ is drawn with respect to $\SketchDist$.

\paragraph{Sketching the supermix problem} 
Recalling the definition of \Cref{eqs:def_sketch_blasso}, we are in presence of a sketched version of the Supermix BLASSO problem, coined S2Mix hereafter, with the following characteristics
\begin{subequations}
    \label[pluralequation]{eqs:sketch_mixture_framework}
\begin{itemize}
    \item \textit{Finite dimensional} Hilbert space: $\mathcal{F}_{\textnormal{sketch}} = \mathbb C^m$.
    \item \textit{Sketching law}: the sequence of \textit{i.i.d.\!\!} random vectors $(\bomega_i)_i$ is drawn with respect to $\SketchDist$,
    \begin{equation}
        \bomega_i\underset{\mathrm{i.i.d.}}\sim\SketchDist\,,\quad i=1,\ldots,m\,.
    \end{equation}
    \item \textit{Sketched\footnote{Related to the \textit{sketching operator} $\mathcal{A}$ \parencite{gribonval2021compressive} which verifies~${\bsketch \defeq \mathcal{A}\hat{f}_n \approx \mathcal{A}(\templatedensity \star \target) = F_{\textnormal{sketch}} \target}$.} forward operator} $F_{\textnormal{sketch}} : \Measures(\Param) \to \mathbb C^m$ defined as
    \begin{equation}
        \label{eq:F_sketched}
        F_{\textnormal{sketch}} \mu = \frac{1}{\sqrt{\nsketch}} \left( \int \varphi_{\bomega_i}(\bt) \diff \mu(\bt) \right)_{i=1}^m.
    \end{equation}
    This forward sketched operator yields the sketched model kernel given by Equation~\eqref{eq:mixture_sketch_model_kernel}, which can in turn be viewed as a Monte-Carlo approximation of the model kernel given by Equation~\eqref{eq:kmod_as_expected_sketch_functions}.
    \item \textit{Sketch vector}: $\bsketch$, a summary of the dataset easily computed as a weighted version the empirical distribution's Fourier transform taken at the random frequencies $\bomega_i$
\begin{equation}    
\label{eq:sketch_definition}
\bsketch \defeq \frac{1}{\sqrt{\nsketch}} \left( W(\bomega_i) \frac{1}{n}\sum_{j=1}^n e^{- \imath \bomega_i^\top \obs_j } \right)_{i=1}^m\,. 
\end{equation}
The latter is of the form: 
\[ \bsketch = F_{\textnormal{sketch}} \target + \Gamma_{\textnormal{sketch}} \quad \textnormal{with: } \Gamma_{\textnormal{sketch}} = \bsketch - F_{\textnormal{sketch}} \target.  \]
The noise term is centered as $\Expectation_{\Obs}\left[\bsketch \right]= \frac{1}{\sqrt{\nsketch}} \left( W(\bomega_i) \mathbb{E}_{\obs \sim \templatedensity \star \target}[e^{- \imath \bomega_i^\top \obs}] \right)_{i=1}^m = F_{\textnormal{sketch}} \target$.
Hence the data sketch is an empirical version of $\Forward_{\textnormal{sketch}}$ at the target point $\target$. 
\end{itemize}
\end{subequations}

\medskip

\begin{remark}[Sketching and compression] Crucially, computing $\bsketch$ acts as a compression of the original $n \times d$ dataset $\Obs$ into an $m$-dimensional sketch in $\mathbb{C}^m$, and is done \emph{only once} prior to learning. This the essence of sketching and compressive learning methods~\parencite{gribonval2021compressive}, seeking statistical guarantees with respect to the sketch size $m$.
\end{remark}

As in the population case, we can use concentration arguments to control the sketched noise level $\gamma_{\textnormal{sketch}}$.
\begin{lemma}
\label{lem:noise_level_sketched}
    Let $\alpha \in (0,1)$, $\Obs = \{\obs_j\}_{j=1}^n$ be an i.i.d. $n$-sample from $f^0$ and $\{\bomega_i\}_{i=1}^m$ be an $m$-sample from $\SketchDist$. Then, with probability at least $1- \alpha$ over the joint draw of the sample and the sketch, it holds that 
    \[
    \Vert \Gamma_{\textnormal{sketch}} \Vert_{\mathbb{C}^m} \leq C_{\alpha,m } \frac{1}{\sqrt{n}},
    \]
    with 
    \[
    C_{\alpha,m} = 2\sqrt{\bigg[ \Big(\frac1\bandwidth\Big)^{d} + \frac{1}{2 \sqrt{m}}\left\Vert \frac{U_\bandwidth}{\SketchDist} \right\Vert_\infty\! \log \Big( \frac{2}{\alpha} \Big) \bigg]  \bigg[ 1 + C_1 \log \Big( \frac{2C_2}{\alpha} \Big) \bigg]}\,,
    \]
    and $C_1, C_2>0$ the same universal constants as in \Cref{prop:control_population_noise_level_whp}.
\end{lemma}

\begin{remark}
\label{rem:link_noise_level_concentration_population_sketched}
    Note that, as $m$ goes to infinity, we uncover that $C_{\alpha,m}$ tends to $C_{\alpha/2} \left(\frac{4}{\bandwidth}\right)^{\frac d2}$ and we recover~\Cref{prop:control_population_noise_level_whp}, which bounds the noise level in the population case.
\end{remark}
The proof is given in~\Cref{proof:lemma_noise_level_sketch} and, again, we may tune the optimal regularisation $\kappa$. Applying the upcoming \Cref{thm:main_thm_sketch} to this specific sketched BLASSO, using the sinc-4 pivot kernel $\Psi_\bandwidth$, one can prove the following results for the S2Mix problem.

\medskip

\begin{proposition}[Statistical guarantees for S2Mix]
\label{prop:sketch_Mixture_Blasso}
    Let $\Param\subset\R^d$ be a compact set, and the target $\target\in \Measures(\Param)$ be a sparse measure of $s_0$ spikes $\{\trueparam_1, \ldots, \trueparam_{s_0}\}$, with minimal separation~${\Delta = \min_{k\neq l} \Vert \trueparam_k - \trueparam_l \Vert^2}$. Fix the bandwidth parameter $\bandwidth \leq \bandwidth_{\max}(\target)$ as in \Cref{hyp:Hbandwidth}. Take any template distribution $\templatedensity$ such that it achieves a finite \emph{kernel switch} constant $C_\switch(\bandwidth, \templatedensity) < + \infty$ in~\Cref{hyp:Hphi}. Take any sketching distribution $\SketchDist$ with support containing~${[-1/\bandwidth, 1/\bandwidth]^d}$ and such that $\Vert U_\bandwidth / \SketchDist \Vert_{\infty} < +\infty$. Let $m \in \mathbb{N}^\star$ be the sketch size and consider the following \emph{sketched} BLASSO problem, denoted as sketched Supermix (S2Mix), with sketched forward measurement operator $F_{\textnormal{sketch}}$ and observation (data sketch) $\bsketch$ defined in \Cref{eq:F_sketched,eq:sketch_definition} respectively:
    \begin{equation}
        \min_\mu \left\{ J_{\textnormal{sketch}, \kappa}(\mu) = \frac{1}{2} \Vert \bsketch - F_{\textnormal{sketch}} \mu \Vert_{\mathbb{C}^m}^2 + \kappa \Vert \mu \Vert_{\mathrm{TV}} \right\}. 
        \addtocounter{equation}{1}
	\tag{\theequation--S2Mix}
    \end{equation}
    Suppose the constant $C_\SketchDist \defeq \big \Vert {f^{(4)}_{\bandwidth}}/{\SketchDist} \big\Vert_\infty < +\infty$, where $f^{(4)}_{\bandwidth} = \fourier [\Psi_\bandwidth] /(2\pi)^d$ is the normalized spectral measure of the sinc-4 kernel $\Psi_\bandwidth$ (\Cref{app:sinc4_RKHS}). Let $C_d, C'_\pivot >0$ be constants 
    that may depend polynomially on the dimension $d$. Fix the regularisation%
    \[
    \kappa = \frac{C_{\alpha, m}}{C'_\pivot C_\switch \sqrt{n}} \frac{1}{\sqrt\K}.
    \]
    Take any diverging sequence $\delta_n = o(n^{1/4})$ and set the near and far regions' \emph{effective radius} $r_n = \delta_n n^{-1/4}$. Let $\alpha > 0$ and the sketch size $m$ be such that
    \begin{equation}
        m \geq C_d\,C_\SketchDist 
        \,
        \max(d,\log(\K))\,
        \K\,\log\bigg(\frac{\max(1,|\Param|)\,\K}{\alpha}\bigg)\,,
    \end{equation}
    Then, for any \emph{near optimal} $\mu \in \Measures(\Param)$ such that $J_{\textnormal{sketch}, \kappa}(\mu) \leq J_{\textnormal{sketch}, \kappa}(\target)$, and $n \geq c_d \delta_n^4$, where $c_d$ depends polynomially on $d$, the following controls hold with probability $1 - \alpha$ over the joint draws of the sample $\Obs$ and the sketches~${(\bomega_1, \ldots, \bomega_m)}$. 
    \begin{itemize}
    \item Control of the far region: 
    \begin{equation}
        |\mu|(\farregion(r))\leq  C'_\pivot C_\switch \frac{512}{69} C_{\alpha,m} \Big(\frac{1}{\delta_n^2 }\Big) \sqrt{\K}\,,
    \end{equation}
    \item Control of all the near regions: for all $k\in[\K]$,
    \begin{equation}
        |\mu(\nearregion_k(r))-a^0_k|\leq  C'_\pivot C_\switch   \frac{1536}{69} C_{\alpha,m} \Big(\frac{1}{\delta_n^2}\Big) \sqrt{\K} \,,
    \end{equation}
    \item Detection level: for all Borelian $A\subset\Param$ such that $|\mu|(A)>  C'_\pivot C_\switch \frac{512}{69} C_{\alpha,m} \Big(\frac{1}{\delta_n^2 }\Big) \sqrt{\K}$, there exists $\trueparam_k$ such that
    \begin{equation}
        \min_{\bt\in A}\distance(t,\trueparam_k)\leq  r\,,
    \end{equation} 
\end{itemize}
\end{proposition}
The detailed proof is given in \Cref{proof:sketch_Mixture_Blasso}.

\subsection{Extensions}

An important consequence of the kernel switch is that one can leverage the previous results for any translation invariant mixture model beyond Gaussian, and any smoothing kernel $\lambda_\bandwidth$ beyond the sinus cardinal, as long as the resulting $K_{\model}$ verifies \Cref{hyp:Hphi}. Since the latter is translation invariant, it suffices to show that its Fourier transform is not ``too localized''. For the sinc-4 pivot, this means $\fourier[K_\model]$ should contain low-frequencies and be bounded below so that $C_\switch$ is finite. One could also consider using another translation invariant pivot kernel, for the Gaussian pivot of~\Cref{eq:sketch_kernel_Gaussian} (not treated in this paper) this means $\fourier[K_\model]$ has a slower decay than the Gaussian.

\section{Effective near regions, kernel switch and near-optimal solutions}
\label{sec:BlassoResults}
This section presents the core theoretical contributions of our work, establishing rigorous guarantees for the BLASSO estimator. We begin by detailing the key assumptions underpinning our analysis, notably the local positive curvature ($\LPC$) property for pivot kernels and conditions for sketching. Our main results then demonstrate three significant advancements: first, we introduce the concept of \textit{effective near regions}, showing that the localisation accuracy of BLASSO adapts to the noise level, improving as more data becomes available. Second, we formalize the \textit{kernel switch} principle, which allows leveraging well-characterized pivot kernels to analyze BLASSO problems with different, potentially more complex, model kernels, provided a specific embedding condition holds. This greatly expands the range of problems amenable to theoretical guarantees. Third, we prove that these statistical error bounds apply not only to exact minimizers of the BLASSO objective but also to \textit{near-optimal solutions}, a crucial finding for practical algorithm design where exact optimisation is often infeasible. Finally, we provide a dedicated analysis for the sinc-4 kernel, proving it satisfies the necessary $\LPC$ and sketching assumptions, thereby enabling its use as a robust pivot kernel, particularly for mixture model estimation as discussed in Section~\ref{sec:Mixtures}.
\subsection{Assumptions}
\label{sec:assumptions}
First, we introduce the local positive curvature assumption, which is a key property for BLASSO analysis. This assumption is crucial to establish the statistical guarantees of our results. We also define the Fisher-Rao distance, which is a Riemannian distance associated with a metric tensor denoted $\mathfrak g_\bx$ and a kernel  $K(\cdot,\cdot)$. 

\paragraph{Fisher-Rao distance for general kernels}
Let $K(\cdot,\cdot)$ be a real valued kernel of positive type. Define the metric tensor (identified as a psd matrix) $\mathfrak g_\bx\defeq\nabla_1\nabla_2 K(\bx,\bx)\in\R^{d\times d}$, where~$\nabla_i$ is the gradient with respect to the $i$th variable. %
Define the Fisher-Rao distance by 
\begin{equation}
    \label{eq:Fisher_Rao_distance}
    \mathfrak{d}_{\mathfrak g}(\bs, \bt)\defeq \inf_{p} \int_0^1\sqrt{{p}'(u)^\top\mathfrak g_{p(u)}{p}'(u)}\,\mathrm{d} u\,,
\end{equation}
where the infimum is taken over smooth paths ${p}\,:\,[0,1]\to \Param$ such that ${p}(0)=\bs$ and ${p}(1)=\bt$. The Fisher-Rao distance is thus the infimum of the length of paths connecting two points in the parameter space $\Param$.

\begin{remark}[Metric of translation invariant kernels]
    \label{rem:metric_translation_invariant}
    The Fisher-Rao distance is a Riemannian distance associated with the metric tensor $\mathfrak g_\bx$ and the kernel $K(\cdot,\cdot)$. For translation invariant kernels defined by 
\begin{subequations}
        \label[pluralequation]{eqs:metric_translation_invariant}
    \begin{equation}
        K(\bs, \bt)=\rho(\bs - \bt),
    \end{equation} 
    the metric tensor is constant and the Fisher-Rao distance is the Mahalanobis distance given by a positive definite matrix~$\mathfrak g$ which does not depend on the points $\bx$, it holds
    \begin{align}
            \mathfrak g& =-\nabla^2 \rho(0)\,,\\
            \distance(\bs, \bt)&=\sqrt{(\bs - \bt)^\top\mathfrak g(\bs - \bt)}\,.
    \end{align}
    \end{subequations}
\end{remark}

\paragraph{Local positive curvature assumption}
We now recall the local positive curvature assumption $\LPC$ as presented in \textcite{poon2023geometry}. The $\LPC$ property ensures that the kernel $K(\cdot,\cdot)$ behaves well in terms of its derivatives and curvature properties with respect to some Fisher-Rao distance, which is essential for the convergence of the BLASSO estimator.

\bigskip

\pagebreak[3]

\begin{subequations}\label[pluralequation]{eqs:ass_positive_curvarure}
\begin{assumption}[$\LPC$ with parameters $\K$, $\Delta_0$, $r_0$, $\bar\varepsilon_0$ and $\bar\varepsilon_2$]
	\label{hyp:local_positive_curvature}
	A real valued kernel $K(\cdot,\cdot)$ of positive type satisfies the local positive curvature assumption $($$\LPC$$)$ if the following holds:
	\begin{itemize}
		\item Assume that, for all $0 \leq i,j\leq 2$ such that $i+j\leq 3$,  
		\begin{equation}
		\label{eq:bound_derivatives}
		B_{ij}\defeq\sup_{\bs,\bt\in\Param}\left\Vert K^{(i,j)}(\bs, \bt) \right\Vert_{\bs,\bt}<+\infty\,.
		\end{equation}
		and denote $B_i\defeq 1+ B_{0i}+ B_{1i}$ where $K^{(i,j)}(\bs, \bt)$ is the covariant derivative\footnote{   We refer to \cite[Section 4.1]{poon2023geometry} for further details on covariant derivatives and their operator norm, and to Appendix~\ref{sec:covariant_derivatives_TI} which gives their expression for translation invariant kernels.} of order $i$ with respect to the first variable $\bs$ and of order $j$ with respect to the second variable $\bt$, and $\|\cdot\|_{\bs,\bt}$ denotes the operator norm with respect to the corresponding metric tensor given by $K(\cdot,\cdot)$ on the tangent spaces at $\bs$ and $\bt$.

		\item Assume that there exists $r_0\in(0,1/\sqrt{B_{02}})$ such that $K(\cdot,\cdot)$ has positive 
        \emph{curvature constants defined as }
		\begin{equation}
        \begin{aligned}
		\bar\varepsilon_0&\defeq\frac12\sup_{\varepsilon\geq 0}
			\Big\{
				\varepsilon\,:\,
				K(\bs, \bt)\leq 1-\varepsilon\,,\ \forall \bs,\bt\in\Param \textnormal{ s.t. }\distance(\bs, \bt)\geq r_0
			\Big\}\,,
            \label{eq:positive_curvature_constants}
            \\
		\bar\varepsilon_2&\defeq\frac14\sup_{\varepsilon\geq 0}
		\Big\{
			\varepsilon\,:\, - K^{(0,2)}(\bs, \bt)[\bv,\bv]\geq \varepsilon \|\bv\|_\bt^2\,,\ \forall \bv\in\mathbb T_\bt\,,\ \forall \bs,\bt\in\Param \textnormal{ s.t. }\distance(\bs, \bt)< r_0
		\Big\}\,,
        \end{aligned}
		\end{equation}
		with $\distance$ the Fisher-Rao distance associated with the metric tensor of $K(\cdot,\cdot)$ and $\mathbb T_\bt$ the tangent space at point $\bt$.
		\item Assume that there exists $\K\geq2$ such that the set 
		\begin{equation}
		\label{eq:delta_constant_definition}
		\Big\{\Delta\,:\,32\sum_{l=2}^{\K}\|K^{(i,j)}(\bx_1,\bx_l)\|_{\bx_1,\bx_l}\leq \min\Big(\frac{\bar\varepsilon_0}{B_0},\frac{2\bar\varepsilon_2}{B_2}\Big),\, 
        \forall (i,j)\in\{0,1\}\times\{0,2\}
		, \forall \{\bx_l\}_{l=1}^{\K}\in\mathcal S_\Delta
		\Big\}
		\end{equation}
		is not empty, with $\mathcal S_\Delta\defeq\big\{\{\bx_l\}_{l=1}^{\K}\in\Param^{\K}\,:\,\min_{k\neq l}\distance(\bx_k,\bx_l)\geq \Delta\big\}$. Denote $\Delta_0$ its infimum.
	\end{itemize}
	In this case, the kernel $K(\cdot,\cdot)$ is said to satisfy the $\LPC$ with parameters $\K$, $\Delta_0$, $r_0$, $\bar\varepsilon_0$ and $\bar\varepsilon_2$.
\end{assumption}
\end{subequations}

\paragraph{Pivot kernels}
Next, we introduce the concept of \textit{pivot kernels}, which are a special class of kernels that can be used to analyze the BLASSO problem. Admissible model kernels are those such that the kernel switch condition (\Cref{hyp:pivot_kernel}) defined below holds. This kernel switch allows us to leverage the properties of the pivot kernel to analyze the model kernel, which may be more complex.

The model kernel introduced in Definition \ref{def:kmod_general} is a kernel of positive type associated to a unique Reproducing Kernel Hilbert Space (RKHS) which is denoted by~$\Hilbert_{\model}$. One can prove that $\Hilbert_{\model}$ is separable and isometric to the closure in~$\mathcal{F}$ of $\mathrm{Im}(F)$, the range of $F$ \parencite[Appendix A.1]{decastro2023fastpart}, 
\[
    \big(\Hilbert_{\model},\|\cdot \|_{\Hilbert_{\model}}\big)
    \simeq 
    \big({\overline{\mathrm{Im}(F)}},\|\cdot \|_{\mathcal{F}}\big)
    \subseteq 
    \big(\mathcal{F},\|\cdot \|_{\mathcal{F}}\big)\,,
\]
where $\simeq$ denotes an isometry mapping. 

\begin{subequations}
Let $F^\star\,:\,\mathcal{F} \to (\Measures(\Param),\|\cdot\|_{\mathrm TV})^\star $ be the dual linear map of $F\,:\,\Measures(\Param)\to \mathcal{F}$. One can prove \parencite[Lemma A.2]{decastro2023fastpart} that %
\begin{align}
    F^\star\,:\,
    c\in \mathcal{F} 
    \mapsto 
    \eta_c
    \in
    \big(\mathcal C(\Param),\|\cdot\|_{\infty}\big)
\end{align}
is a continuous linear operator, and 
\begin{align}
 \langle c,F\nu \rangle_{\mathcal{F}}
 = \int_\Param \langle c,F\delta_\bt\rangle_{\mathcal{F}}\diff\nu(t)
 = \langle F^\star c, \nu\rangle_{\mathcal C(\Param),\Measures(\Param)}
 = \int_\Param \eta_c\diff\nu(t)
 \quad\textnormal{for all } c\in\mathcal{F},\,\nu\in\Measures(\Param)\,.
\end{align}
Now, define
\begin{equation}
    \label{def:eta_c}
    \forall c\in\mathcal{F}\,,\quad 
\eta_c\,:\,t\in\Param\mapsto \langle c,F\delta_\bt \rangle_{\mathcal{F}}\in\R\,,
\end{equation}
\end{subequations}
we have the following useful proposition.

\medskip

\begin{proposition}
\label{prop:model_description}
    If \eqref{hyp:continuous_model_kernel} holds then for all $\eta\in\Hilbert_\model$, there exists a unique $c\in {\overline{\mathrm{Im}(F)}}$ such that $\|\eta\|_{\Hilbert_\model}=\|c\|_{\mathcal{F}}$ and $\eta=\eta_c$. Furthermore, it holds that $\eta=F^\star c$.
\end{proposition}

\begin{proof}
    See Appendix~\ref{proof:prop_model_description}
\end{proof}

We can now turn to the the definition of a pivot kernel $K_\pivot(\cdot,\cdot)$
with respect to the model kernel $K_\model(\cdot,\cdot)$. %

\begin{subequations}\label[pluralequation]{eqs:pivot_assumptions}
\begin{assumption}[Pivot kernels]
\label{hyp:pivot_kernel}
A real valued kernel $K_\pivot(\cdot,\cdot)$ of positive type is a pivot kernel with respect to $K_\model(\cdot,\cdot)$ if a) its values at points $(\bt,\bt)$ are normalized, namely
\begin{equation}
    \label{eq:pivot_normalisation}
  \forall \bt\in\Param\,,\quad K_{\pivot}(\bt,\bt)=1\,;  
\end{equation}
and b) its RKHS, denoted by~$\Hilbert_{\pivot}$, is such that $\Hilbert_{\pivot}\subseteq \Hilbert_{\model}$, and it holds that the following identity map~${\mathrm{Id}_{\pivot,\model}\,:\,(\Hilbert_{\pivot},\|\cdot \|_{\Hilbert_{\pivot}})\to (\Hilbert_{\model},\|\cdot \|_{\Hilbert_{\model}})}$ is continuous, namely
\begin{equation}
    \label{eq:constant_switch}
C_{\switch}(K_\pivot, K_\model) \defeq\big\|\mathrm{Id}_{\pivot,\model}\big\|_{\mathrm{op}}
=
\sup_{\eta\in \Hilbert_{\pivot}\setminus\{0\}}
\frac{\|\eta\|_{\Hilbert_\model}}{\|\eta\|_{\Hilbert_\pivot}}<\infty\,.
\end{equation}
\end{assumption}
\end{subequations}

\paragraph{Tail bound assumption on the sketching function derivatives.}
The third and final assumption is a tail bound on the derivatives of the sketching function $\psi_\bomega$ defined in \Cref{eq:pivot_sketch_kernel}, which also depends on the sketching distribution $\SketchDist$. It ensures that the derivatives of the sketching function do not grow too large with respect to the Fisher-Rao distance, which is important to control the sketching error. The tail bounds are defined in terms of a so-called survival function $\bar{\mathbb{P}}_j(t)$, which quantifies the probability that the derivatives of the sketching function exceed a certain threshold $t$. 

\begin{subequations}\label[pluralequation]{eqs:ass_sketch}
\begin{assumption}[Tail bound on the sketching function derivatives]
\label{hyp:sketch_function_bounded}
The function~$\psi_\bomega$ and the law $\bomega\sim\SketchDist$ satisfy the derivatives tail bounds assumption if: a) the kernel $K(\cdot,\cdot)$ defined by~\eqref{eq:pivot_sketch_kernel} satisfies the local positive curvature assumption $($Assumption~\ref{hyp:local_positive_curvature}$)$ with some parameters $\K$, $\Delta_0$, $r_0$, $\bar\varepsilon_0$ and $\bar\varepsilon_2$; and b) there exists a constant~${C>0}$ which depends on $K(\cdot,\cdot)$ and polynomially on $d$ such that, for all~$\alpha\in(0,1)$, there exist positive reals $L_0,L_1,L_2,L_3>0$ which may depend on $\alpha$ verifying:
\begin{align}
    \sum_{j=0}^3\bar{\mathbb{P}}_j(L_j)
    &\leq\frac{\min(\bar\varepsilon_0,\bar\varepsilon_2,\alpha)}{m_0}    
    \\
    \max_{j=0}^3\Big(L_j^2\sum_{i=0}^3\bar{\mathbb{P}}_i(L_i)
    +6\int_{L_j}^\infty t\, \bar{\mathbb{P}}_j(t)\,\mathrm dt
    \Big)
    &\leq\frac{\min(\bar\varepsilon_0,\bar\varepsilon_2)}{m_0}    
    \\
    m_0
    &\defeq C \K
    \bigg(
    C_1\log(\K)\log\Big(\frac{\K}{\alpha}\Big)
    +C_2\log\bigg(\frac{(\K N)^d}{\alpha}\bigg)
    \bigg)
\end{align}
where 
\begin{align}
    \label{eq:survival_function}
    \bar{\mathbb{P}}_j(t)
    &\defeq1-\mathbb P_{\bomega\sim\SketchDist}\bigg(\sup_{\bx\in\Param}\|\psi_\bomega^{(j)}(\bx)\|_\bx\leq t\bigg)\\ \notag
    N
    &\defeq \frac{|\Param|L_1}{\bar\varepsilon_0}+\frac{r_0 L_3L_0+L_2}{\bar\varepsilon_2}
    \\ \notag
    C_1
    &\defeq (L_0^2+L_1^2)(B_0/\bar\varepsilon_0^2+B_2/\bar\varepsilon_2^2)
    \\ \notag
    C_2
    &\defeq L_0^2/\bar\varepsilon_0^2+L_{01}L_0/\bar\varepsilon_0
    +L_2^2/\bar\varepsilon_2^2+L_{01}L_2/\bar\varepsilon_2
    \\ \notag
    L_{ij}
    &\defeq \sqrt{L_i^2+L_j^2}
    \\ \notag
    |\Param|
    &\defeq \sup_{\bs,\bt\in\Param}\distance(\bs, \bt)
\end{align}
with $\psi_\bomega^{(j)}(\bx)$ the covariant derivative of order $j$ at $\bx$, $\distance(\bs, \bt)$ the Fisher-Rao distance, $\|\cdot\|_{\bx}$ the norm of the corresponding metric tensor given by $K(\cdot,\cdot)$ on the tangent space at $\bx$.
\end{assumption}
\end{subequations}

\subsection{Main theorem}

\begin{subequations}
\label[pluralequation]{eqs:main_thm}
\begin{theorem}[Estimation error bounds for BLASSO]
\label{thm:main_thm}
Let $\Param$ be a compact set of $\R^d$, let $\mathcal{F}$ be a separable Hilbert space, and let $K_\pivot(\cdot,\cdot)$ be a kernel of positive type which satisfies Assumption~\ref{hyp:local_positive_curvature} with parameters $\K$, $\Delta_0$, $r_0$, $\bar\varepsilon_0$ and $\bar\varepsilon_2$, and which satisfies Assumption~\ref{hyp:pivot_kernel} with constant $C_{\switch}<+\infty$. 

Let $F\,:\,\Measures(\Param)\to \mathcal{F}$ be a linear map and let $\target\in\mathcal M(\Param)$, $\by \in\mathcal{F}$ and $\gamma\geq 0$ be such that
\begin{equation}
\notag
    \target\in\Model_{\K, \Delta_0, \distance}
    \quad\textnormal{and}\quad
    \gamma\geq\|\by-F\target\|_{\mathcal{F}}\,,
\end{equation}
with $\distance$ the Fisher-Rao distance associated to $K_\pivot(\cdot,\cdot)$. 

 Consider $J_\kappa$ given by~\eqref{eq:blasso} where $\kappa=c_\kappa\gamma/\sqrt{\K}$ with $c_\kappa>0$ some constant, and let $\mu\in\mathcal{M}(\Param)$ be such that 
\begin{equation}
\notag
    J_\kappa(\mu)\leq J_\kappa(\target)\,.
\end{equation}
Then for each $r>0$ such that $r<\min\big(r_0,\sqrt{\bar\varepsilon_0/\bar\varepsilon_2}\big)$, we have:
\begin{itemize}
    \item Control of the far region: 
    \begin{equation}
    \label{eq:control_far_thm}
        |\mu|(\farregion(r))\leq \bar c_\kappa\Big(\frac\gamma{\bar\varepsilon_2 r^2}\Big) \sqrt{\K}\,,
    \end{equation}
    \item Control of all the near regions: for all $k\in[\K]$,
    \begin{equation}
    \label{eq:control_near_thm}
        |\mu(\nearregion_k(r))-a^0_k|\leq \tilde c_\kappa \Big(\frac\gamma{\bar\varepsilon_2r^2}\Big) \sqrt{\K}+\hat c_\kappa \gamma\,,
    \end{equation}
    \item Detection level: for all Borelian $A\subset\Param$ such that $|\mu|(A)>\bar c_\kappa\big(\frac\gamma{\bar\varepsilon_2 r^2}\big)\sqrt{\K}$, there exists $\trueparam_k$ such that
    \begin{equation}
    \label{eq:detection_near_thm}
        \min_{\bt\in A}\distance(t,\trueparam_k)\leq  r\,,
    \end{equation} 
\end{itemize}
where the far and near regions are defined by Definition~\ref{def:Far_and_Near_regions} and
\begin{align}
    \bar c_\kappa&\defeq\frac{(1+\sqrt 2\,C_{\switch}c_\kappa)^2}{2c_\kappa}\geq2\sqrt 2\,C_{\switch}\\
    \tilde c_\kappa&=\bar c_\kappa\max(1,\bar\varepsilon_0)\\
    \hat c_\kappa&=2\sqrt 2\,C_{\switch}(1+\sqrt 2C_{\switch}c_\kappa)%
\end{align}
\end{theorem}
\end{subequations}
The proof of is given in Appendix~\ref{proof:main_thm}. In the statement above, one can see that the choice of the regularisation parameter is $\kappa \propto \gamma / \sqrt\K$ up to a free constant $c_\kappa$. The selection of $c_\kappa$ influences the constants in the error bounds and their dependency on $\K$. We discuss several choices below:
\begin{remark}[Tuning the regularisation parameter $\kappa = c_\kappa \gamma / \sqrt{\K}$]
    \label{rem:tuning_kappa_main_thm} 
    The constant $c_\kappa$ in the definition of the regularisation parameter $\kappa$ affects the pre-factors in the error bounds \eqref{eq:control_far_thm}, \eqref{eq:control_near_thm}, and \eqref{eq:detection_near_thm}.
    \begin{enumerate}
        \item \textbf{Optimal scaling with $\sqrt{\K}$}: Choosing $c_\kappa = 1 / (\sqrt{2} C_\switch)$ minimizes the constant $\bar{c}_\kappa$ to its lower bound $2 \sqrt{2} C_\switch$. This yields error bounds that scale as $\mathcal{O}(\gamma \sqrt{\K} / (\bar\varepsilon_2 r^2))$ for the far region control and the leading term in the near region control. This is generally the preferred choice for optimal statistical rates.

        \item \textbf{Alternative scaling with $\K$}: If one chooses $c_\kappa$ such that $\kappa$ becomes independent of ${\K}$ ({\it e.g.,} by setting~$c_\kappa = \sqrt{\K}$), then $\kappa \propto \gamma$. This leads to $\bar{c}_\kappa \propto \K$, degrading the error bounds to scale as $\mathcal{O}(\gamma \K / (\bar\varepsilon_2 r^2))$. While this might simplify the expression for $\kappa$, it thus results in suboptimal statistical guarantees.

        \item \textbf{Caution on noise level dependency}: It is crucial that $\kappa$ maintains its proportionality to the noise level~$\gamma$. For instance, attempting to make $\kappa$ independent of $\gamma$ (e.g., by setting $c_\kappa \propto 1/\gamma$) would lead to error bounds with a constant bias term, such as $\mathcal{O}((\text{constant} + C_\switch)\sqrt{\K})$, which does not vanish as $\gamma \to 0$. This is undesirable in statistical applications where $\gamma$ typically decreases with increasing sample size (e.g., mixture model estimation, see Section~\ref{sec:Mixtures} and Proposition~\ref{prop:supermix_thm10}).
    \end{enumerate}
\end{remark}

\subsection{New guarantees for sketched BLASSO}
We now consider a sketched BLASSO given by a forward operator $F$ as in \eqref{eq:sketch_forward} with sketching function $\varphi_\bomega(\cdot)$, whose model kernel is given by~\eqref{eq:population_model_kernel}.

\begin{theorem}[Estimation error bounds for sketched BLASSO]
\label{thm:main_thm_sketch}
Let $\Param$ be a compact set of $\R^d$. Let $\varphi_\bomega,\psi_\bomega\,:\,\Param\to\mathbb{C}$ be functions and let $\SketchDist$ be a law on $\mathbb R^d$. Consider 
\begin{align}
    {K}_{\model}(\bs, \bt)&= \Expectation_{\bomega\sim \SketchDist}\big[\varphi_\bomega(s)\overline{\varphi}_\bomega(t)\big]\,,
    \notag
    \\
    {K}_{\pivot}(\bs, \bt)&= \Expectation_{\bomega\sim \SketchDist}\big[\psi_\bomega(s)\overline{\psi}_\bomega(t)\big]\,.
    \tag{$\mathbf{H}_{\psi_\bomega,\SketchDist}$}
\end{align}
Assume that $K_\pivot(\cdot,\cdot)$ satisfies Assumption~\ref{hyp:local_positive_curvature} $($with parameters $\K$, $\Delta_0$, $r_0$, $\bar\varepsilon_0$ and $\bar\varepsilon_2)$ and Assumption~\ref{hyp:pivot_kernel} with constant $C_{\switch}>0$. Assume that the sketching function $\psi_\bomega$ satisfies Assumption~\ref{hyp:sketch_function_bounded} and that its derivatives are $\SketchDist$-almost surely bounded up to order $3$. Denote $C_{\textnormal{sketch}} \defeq 2C\max(C_1, C_2) > 0$. 

Let $\alpha\in(0,1)$ and $m$ be such that
\begin{equation}
    \label{eq:m0_thm}
    m\geq C_{\textnormal{sketch}}\,\max(d,\log(\K))\,\K\,\log\Big(\frac{\max(1,|\Param|)\,\K}{\alpha}\Big)\,,
\end{equation}
where $\displaystyle|\Param|= \sup_{\bs,\bt\in\Param}\distance(\bs, \bt)$ is the diameter of $\Param$ with respect to the Fisher-Rao distance $\distance$ given by $K_\pivot(\cdot,\cdot)$. 

\begin{subequations}
Let $F\,:\,\Measures(\Param)\to \mathbb C^m$ be given by 
    \begin{equation*}
        (F\mu)_i= \frac{1}{\sqrt{m}}\int_\Param {\varphi}_{\bomega_i}(t) \mathrm{d}\mu(t)\,,\quad i=1,\ldots,m\,.
    \end{equation*}
where the sequence of \textit{i.i.d.\!\!} random vectors $(\bomega_i)_i$ is drawn with respect to $\SketchDist$. Let $\target\in\mathcal M(\Param)$, $\bsketch\in\mathbb C^m$ and $\gamma\geq 0$ be such that
\begin{equation}
\notag
    \target\in\Model_{\K, \Delta_0, \distance}
    \quad\textnormal{and}\quad
    \gamma\geq\|\bsketch-F\target\|_{\mathbb C^m}\,,
\end{equation}

 Consider $J_\kappa$ given by~\eqref{eq:blasso} where $\kappa=c_\kappa\gamma/\sqrt{\K}$ with $c_\kappa>0$ some constant, and let
 $\mu\in\mathcal{M}(\Param)$ be such that 
\begin{equation}
\notag
    J_\kappa(\mu)\leq J_\kappa(\target)\,.
\end{equation}

Then, there exists a constant $C'_{\pivot}>0$ which depends only on the kernel $K_\pivot$ such that, with probability at least $1-\alpha$ on the draw of $(\bomega_1, \ldots, \bomega_m)$, for any $r>0$ such that $r<\min\big(r_0,\sqrt{\bar\varepsilon_0/6\bar\varepsilon_2}\big)$, we have:
\begin{itemize}
    \item Control of the far region: 
    \begin{equation}
    \label{eq:control_far_thm_sketch}
        |\mu|(\farregion(r))\leq \bar c_\kappa\Big(\frac{2\gamma}{3\bar\varepsilon_2 r^2}\Big) \sqrt{\K}\,,
    \end{equation}
    \item Control of all the near regions: for all $k\in[\K]$,
    \begin{equation}
    \label{eq:control_near_thm_sketch}
        |\mu(\nearregion_k(r))-a^0_k|\leq \tilde c_\kappa \Big(\frac{2\gamma}{3\bar\varepsilon_2r^2}\Big) \sqrt{\K}+\hat c_\kappa \gamma\,,
    \end{equation}
    \item Detection level: for all borelian $A\subset\Param$ such that $|\mu|(A)>\bar c_\kappa\big(\frac{2\gamma}{3\bar\varepsilon_2 r^2}\big)\sqrt{\K}$, there exists $\trueparam_k$ such that
    \begin{equation}
    \label{eq:detection_near_thm_sketch}
        \min_{\bt\in A}\distance(t,\trueparam_k)\leq  r\,,
    \end{equation} 
\end{itemize}
where the far and near regions are defined by Definition~\ref{def:Far_and_Near_regions} and
\begin{align}
    \bar c_\kappa&\defeq\frac{(1+{C'_\pivot}C_{\switch}c_\kappa)^2}{2c_\kappa}\geq2{C'_\pivot}C_{\switch}\\
    \tilde c_\kappa&=\bar c_\kappa\max(1,\bar\varepsilon_0/4)\\
    \hat c_\kappa&=2{C'_\pivot}C_{\switch}(1+{C'_\pivot}C_{\switch}c_\kappa)%
\end{align}
\end{subequations}
\end{theorem}

The proof of is given in Appendix~\ref{proof:main_thm_sketch}.

\begin{remark}
    The sketch size \eqref{eq:m0_thm} presented in the statement of Theorem~\ref{thm:main_thm_sketch} is given for pivot sketching functions $\psi_\bomega$ with $\SketchDist$-almost surely bounded derivatives up to the order $3$. One could weaken this hypothesis and only ask that the survival function $\bar{\mathbb P}_j(t)$ defined in \Cref{hyp:sketch_function_bounded} decays exponentially, which would incur additional logarithmic terms in $\K$ and $d$ in the sketch size \eqref{eq:m0_thm}.
\end{remark}

\begin{remark}
    For any $m\geq 1$, one can set 
    \[
    \alpha= \max(1,|\Param|)\,\K\,\exp\bigg(-\frac{m}{C_{\textnormal{sketch}}\,\max(d,\log(\K))\,\K}\bigg)
    \]
    and the result of Theorem~\ref{thm:main_thm_sketch} holds with probability $1-\min(1,\alpha)$. We uncover that $m=\mathcal{O}({\K})$ (up to multiplicative log factors in $\K$ and multiplicative polynomial factors in $d$) are sufficient to get accurate statistical estimation error bounds. %
\end{remark}

\subsection{The sinc-4 pivot: LPC and sketching}
We consider the sinc-4 pivot $\Psi_{\bandwidth}(\bx - \by) = \sinc^4(\frac{\bx-\by}{4 \bandwidth})$ defined in Equation~\eqref{eq:sketch_kernel_sinc}. This kernel has been studied in \cite{decastro2019sparse} without proving that the $\LPC$ holds. Since the theoretical analysis of the sketched BLASSO is built upon $\LPC$ \parencite{poon2023geometry}, the statistical error bounds of sketched BLASSO with the sinc-4 pivot was an open question. The next theorem shows that the sinc-4 kernel $\Psi_\bandwidth$ satisfies Assumption~\ref{hyp:local_positive_curvature}. Note that its Fisher-Rao distance is 
\begin{subequations}
\begin{equation}
  \distancegeneric_{\metric, \bandwidth}(\bx, \by) = \frac{1}{2\sqrt{3}\, \bandwidth} \Vert \bx - \by \Vert_2\,,  
\end{equation}
which scales as $\bandwidth^{-1}$.
\begin{theorem}
    \label{thm:sinc4_lca}
    Given any $\bandwidth>0$ and $\K\geq1$, the kernel $\Psi_{\bandwidth}(\bx - \by) = \sinc^4(\frac{\bx-\by}{4 \bandwidth})$ on $\mathcal X=\R^d$ satisfies Assumption~\ref{hyp:local_positive_curvature} with parameters $\K$, $\Delta_0=42.66\, \K^{1/4}d^{7/4}$, $r_0=1/(4\,{d})$, $\bar\varepsilon_0 \geq 1/(32\,d^3)$ and $\bar\varepsilon_2 \geq 23/128$.
\end{theorem}

The proof is given in Appendix~\ref{proof:sinc4_lca}. Given a set of distinct points $\{\bx_1^0,\ldots,\bx_{\K}^0\}\subset \R^d$, if 
\begin{equation}
\label{eq:condition_bandwidth}
\bandwidth\leq  \frac{\min_{k,l}\|\bx_k^0-\bx_l^0\|_2}{147.77\,\K^{1/4}d^{7/4}}
\end{equation}
then any measure supported on these points is a sparse target measure \eqref{eq:class_model} for which $\Delta_0\leq42.66\, \K^{1/4}d^{7/4}$ and hence the $\LPC$ holds with $r_0=1/(4\,{d})$, $\bar\varepsilon_0 \geq 1/(32\,d^3)$ and $\bar\varepsilon_2 \geq 23/128$. Condition~\eqref{eq:condition_bandwidth} is exactly~\eqref{hyp:Hbandwidth}.

\medskip
We now turn to the proof of \Cref{hyp:sketch_function_bounded} for the sinc-4 kernel.  As a translation invariant kernel, we denote $f^{(4)}_{\bandwidth}$ its spectral measure, detailed in~\Cref{app:sinc4_RKHS} and such that
\[ 
\Psi_\bandwidth(\bx - \by)= \int e^{+ \imath \omega^\top (\bx - \by)} f^{(4)}_{\bandwidth}(\bomega) \diff \bomega.
\]
Since $\Psi_\bandwidth(\bm{0})=1$ and $\Psi_\bandwidth$ is a positive type kernel, $f^{(4)}_{\bandwidth}$ is a probability density function. Thus, for any sketching distribution $\SketchDist$, we may write the latter as an expectation in the form of~\Cref{eq:pivot_sketch_kernel} 
\[\Psi_\bandwidth(\bx - \by) =  \Expectation_{\bomega\sim \SketchDist}\big[\psi_\bomega(\bx)\overline{\psi}_\bomega(\by)\big], \quad \textnormal{with sketching function } \psi_\bomega(\bt) = e^{\imath \bomega^\top \bt} \sqrt{\frac{f^{(4)}_{\bandwidth}}{\SketchDist}} (\bomega).
\]
In order to prove~\Cref{hyp:sketch_function_bounded} we need to exhibit $L_0, L_1, L_2, L_3 >0$ such that~\Cref{eqs:ass_sketch} is verified for the $\LPC$ constants of \Cref{thm:sinc4_lca}. This is the purpose of the following theorem.

\begin{theorem}
    \label{thm:sinc4_ass2}
    For any bandwith $\bandwidth$, the sinc-4 kernel $\Psi_\bandwidth(\bx - \by) = \sinc^4(\frac{\bx - \by}{4 \bandwidth})$ and its sketching function $\psi_\bandwidth$ verify \Cref{hyp:sketch_function_bounded} with:
    \begin{align}
    m_0 & = C\K
    \Bigg(
    C_1\log(\K)\log\Big(\frac{\K}{\alpha}\Big)
    +C_2\log\bigg(\frac{(\K N)^d}{\alpha}\bigg)
    \Bigg) \\
    N
    &\defeq |\Param| 32\sqrt{12} \sqrt{C_\SketchDist} d^{7/2} + \frac{128}{23}( 12 \sqrt{12} C_\SketchDist d^{1/2} + \sqrt{C_\SketchDist} 12 d) 
    \\ \notag
    C_1
    &\defeq (1 + 12 d) C_\SketchDist \left(1024 d^6 (2 + \sqrt{12d}) + \frac{128^2}{23^2} (1 + \sqrt{12d} + 12d)\right) = \mathcal{O}\left(d^{15/2} C_\SketchDist\right)
    \\ \notag
    C_2
    &\defeq C_\SketchDist \left(1024  d^6 + 32d^3 \sqrt{1 + 12d} + 
    + \frac{128^2}{23^2} 12^2  d^2 + \frac{128}{23} 12 d  \sqrt{1 + 12d} \right) = \mathcal{O}\left(d^{6} C_\SketchDist\right)
    \\ \notag
    C_\SketchDist 
    &\defeq \sup_{\bomega \in [-\frac{1}{\bandwidth}, \frac{1}{\bandwidth}]^d} \frac{f^{(4)}_{\bandwidth}}{\SketchDist} (\bomega) 
    \end{align}
    where $|\Param| = \sup_{\bs,\bt\in\Param}\distance(\bs, \bt)$ and $C$ depends on $\Psi_\bandwidth$ and polynomially on $d$. Moreover, derivatives of $\psi_\bandwidth$ are almost surely bounded at any order.
\end{theorem}
The proof is given in Appendix~\ref{proof:sinc4_ass2}.
\end{subequations}

Lastly, the following proposition highlights the flexibility of the sinc-4 kernel as pivot for the analysis of BLASSO problems involving translation-invariant kernels with sufficiently slow-decaying spectral measure.

\begin{proposition}[The sinc-4 pivot for translation-invariant kernels]
\label{prop:sinc4_cswitch}
Let $K_\model(\bs, \bt) = \rho_\model(\bs - \bt)$ be a translation-invariant reproducing kernel with a non-zero Fourier transform on low frequencies $\mathbb{B}_\bandwidth = [-1/\bandwidth, 1/\bandwidth]^d$. Then, we have that
\[
C_\switch(\Psi_\bandwidth, \rho) = \essentialsup_{\bomega \in \mathbb{B}_\bandwidth }\sqrt{
\frac{\fourier[\Psi_\bandwidth]}{\fourier[\rho_\model]}(\bomega)} \leq  \frac{
\essentialsup_{\bomega \in \mathbb{B}_\bandwidth } \sqrt{\fourier[\Psi_\bandwidth](\bomega)}
}{\essentialinf_{\bomega \in \mathbb{B}_\bandwidth  } \sqrt{\fourier[\rho_\model](\bomega)}}.
\]
and if $\fourier[\rho_\model]$ is bounded away from $0$ almost everywhere on low-frequencies, \textit{i.e.} $\essentialinf_{\mathbb{B}_\bandwidth} \fourier[\rho_\model] > 0$,
then $C_\switch(\Psi_\bandwidth, \rho) < +\infty$ and $\Hilbert_{\textnormal{sinc-4}}$ continuously embeds into $\Hilbert_\rho$.
\end{proposition}

The proof is given in~\Cref{app:proof_sinc4_csiwtch} and is a consequence of~\textcite[Corollary 3.2]{zhang2013inclusion}. For the smoothed mixture models considered in~\Cref{sec:Mixtures}, the model kernel is $\rho_\templatedensity = \lambda_\bandwidth \star \templatedensity \star \check{\templatedensity}$ which Fourier transform $\fourier[\rho_\templatedensity] = \fourier[\lambda_\bandwidth] | \fourier[\templatedensity] |^2$ contains low-frequencies. This yields the kernel switch constant of~\Cref{hyp:Hphi}
\begin{equation}
    C_\switch(\Psi_\bandwidth, \rho_\templatedensity) = C_\switch(\bandwidth, \templatedensity)
    \leq  \frac{\essentialsup_{\bomega \in \mathbb{B}_\bandwidth } \sqrt{\fourier[\Psi_\bandwidth](\bomega)}}{\essentialinf_{\bomega \in \mathbb{B}_\bandwidth  } \sqrt{\fourier[\lambda_\bandwidth] \big | \fourier[\templatedensity]\big|^2(\bomega)}},
\end{equation}
 which is finite when $\essentialinf \fourier[\rho_\templatedensity] > 0$. In this work, we use the rescaled sinus cardinal as $\lambda_\bandwidth$, which has a \emph{constant} Fourier transform $\fourier[\lambda_\bandwidth] \propto \indicator_{\mathbb{B}_\bandwidth }$. Thus, any template distribution $\templatedensity$ such that $\vert \fourier[\templatedensity] \vert$ is bounded away from 0 on $\mathbb{B}_\bandwidth $ is amenable to the kernel switch in~\Cref{prop:supermix_thm10}, and the scaling of the constant $C_\switch(\bandwidth, \templatedensity)$ depends on the smoothness of $\fourier[\templatedensity]$.
 
\stoptoc

\section{Conclusion and perspectives}

This paper advances the theory of continuous sparse regularisation on measures, providing novel theoretical guarantees for sketched problems, including mixture models. We introduce a key embedding constant, $C_\switch$, which we explicitly characterize for translation-invariant kernels via a ratio of Fourier transforms for translation invariant kernels. Our analysis further reveals how the support of the solution belongs to so-called effective near regions that scale with the noise level and establishes formal error bounds for near-optimal solutions, validating them as reliable estimators in pratice.

These contributions open promising avenues for future research. Theoretically, a natural extension is to exhibit the kernel switch constant for non-translation invariant models, such as Gaussian mixtures with unknown covariance. Methodologically, our work motivates the development of practical conic particle gradient descent algorithms for the BLASSO. This pursuit will also require addressing the crucial challenge of selecting the regularisation parameter~$\kappa$ and the sketch size in a data-driven manner, especially when the true number of components is unknown.

Finally, we note that our results are not limited to mixture models but can be applied to any continuous regression problem where some pivot kernel satisfies the $\LPC$. This includes a wide range of applications in machine learning and statistics, such as density estimation, clustering, and regression tasks on continuous data ({\it e.g.,} tensor regression \parencite{azais2024second}, shallow neural networks \parencite{chizat2019sparse, azais2024second}).

\paragraph{Acknowledgments} The authors would like to thank \emph{Romane Giard} for precious comments and discussions on draft versions of this article. 

\paragraph{Funding} The work of R. Gribonval was partially supported by the AllegroAssai ANR project ANR-19-CHIA-0009 and the SHARP ANR Project ANR-23-PEIA-0008 of the PEPR IA, funded in the framework of the France 2030 program.

\begin{table}[!hb]
    \newcommand{\ttile}[1]{\textsf{\textbf{#1}}}
    \centering
    {\small
    \begin{tabular}{ll}
		\multicolumn{2}{c}{\textit{General notation}} \\
		\midrule
        $[n]$                                   & Set of integers $\{1,\ldots,n\}$;\\
    &   \\
		\multicolumn{2}{c}{Continuous regression} \\
		\midrule
        $\Param$                                & Compact  set of $\R^d$;\\
        $(\mathcal{C}(\Param),\|\cdot\|_\infty)$        & Continuous functions on $\Param$;\\
        $(\mathcal{M}(\Param),,\|\cdot\|_{\mathrm{TV}})$& Radon measures on $\Param$;\\
        $(\mathcal{F},\langle\cdot,\cdot\rangle_{\mathcal{F}})$       
                                                & Separable Hilbert space;\\
        $\kappa$                                 & Regularisation parameter of \eqref{eq:blasso};\\
        $\Model_{\K, \Delta_0, \distance}$                 & Class of models, Eq.~\eqref{eq:class_model};\\
        $\Gamma$ (resp. $\gamma$)               & Noise (resp. noise level), Eq.~\eqref{eqs:noise_level_term};\\
        $\nearregion_k(r)$ (resp. $\farregion(r)$) &Near (resp. far) regions, Eq.~\eqref{eqs:near_far_def};\\
            &   \\
		\multicolumn{2}{c}{Linear operators and kernels} \\
		\midrule
        $K_\model$                              & Model kernel, Eq.~\eqref{def:kmod_general};\\
        $(\Hilbert_{\model},\|\cdot\|_{\Hilbert_{\model}})$ & RKHS associated to $K_\model$, Eq.~\eqref{eq:hmod_representation};\\
        $K_\pivot$                              & Pivot kernel, Eq.~\eqref{eqs:pivot_assumptions};\\
        $(\Hilbert_{\pivot},\|\cdot\|_{\Hilbert_{\pivot}})$ & RKHS associated to $K_\pivot$;\\
        $\varphi_\bomega (\textnormal{resp. } \psi_\bomega)$         & Model (resp. pivot) sketching functions;\\
        $\Theta_\Omega$                         & Radial basis function kernel with covariance $\Omega$, Eq.~\eqref{eq:sketch_kernel_Gaussian};\\
	$\Psi_\bandwidth$                             & Sinus cardinal kernel to the four, Eq.~\eqref{eq:sketch_kernel_sinc};\\
        $\SketchDist$                           & Sketching distribution, Eq.~\eqref{eq:sketch_distribution};\\
        $F$                                     & Forward operator, Eq.~\eqref{eq:sketch_forward};\\
        $F^\star$                               & Adjoint operator of $F$;\\
        $\fourier[\cdot]$                      & Fourier transform;\\
        $\fourier^{-1}[\cdot]$                 & Inverse Fourier transform;\\
    &   \\
		\multicolumn{2}{c}{Differential geometry} \\
		\midrule
        $\nabla_1 K(\bs, \bt)$                       & Gradient at point $\bs$ of $\bu\mapsto K(\bu, \bt)$; \\
        $\nabla_2 K(\bs, \bt)$                       & Gradient at point $\bt$ of $\bu\mapsto K(\bs, \bu)$; \\
        $\mathfrak g_\bx$                         & Metric tensor $\nabla_1\nabla_2 K(\bx,\bx)$ at point $\bx$;\\
        $\distance(\bs, \bt)$                        & Fisher-Rao distance between $\bs$ and $\bt$ with respect to a kernel $K$ with metric $\mathfrak{g}$, Eq.~\eqref{eq:Fisher_Rao_distance};\\
        $\mathfrak{d}_{\metric, \bandwidth}(\bs, \bt)$ & Fisher-Rao distance between $\bs$ and $\bt$ with respect to the sinc-4 kernel $\Psi_\bandwidth$;\\
    &   \\
  \end{tabular}
  }
\caption{List of notations.}
\end{table}
\stoptoc

\newpage
\printbibliography

\newpage

\appendix

\section{Proofs of the main results}

\subsection{Existence of pivot non-degenerate dual certificates}

\begin{figure}[!ht]
	\centering
	\includegraphics[width=.49\linewidth]{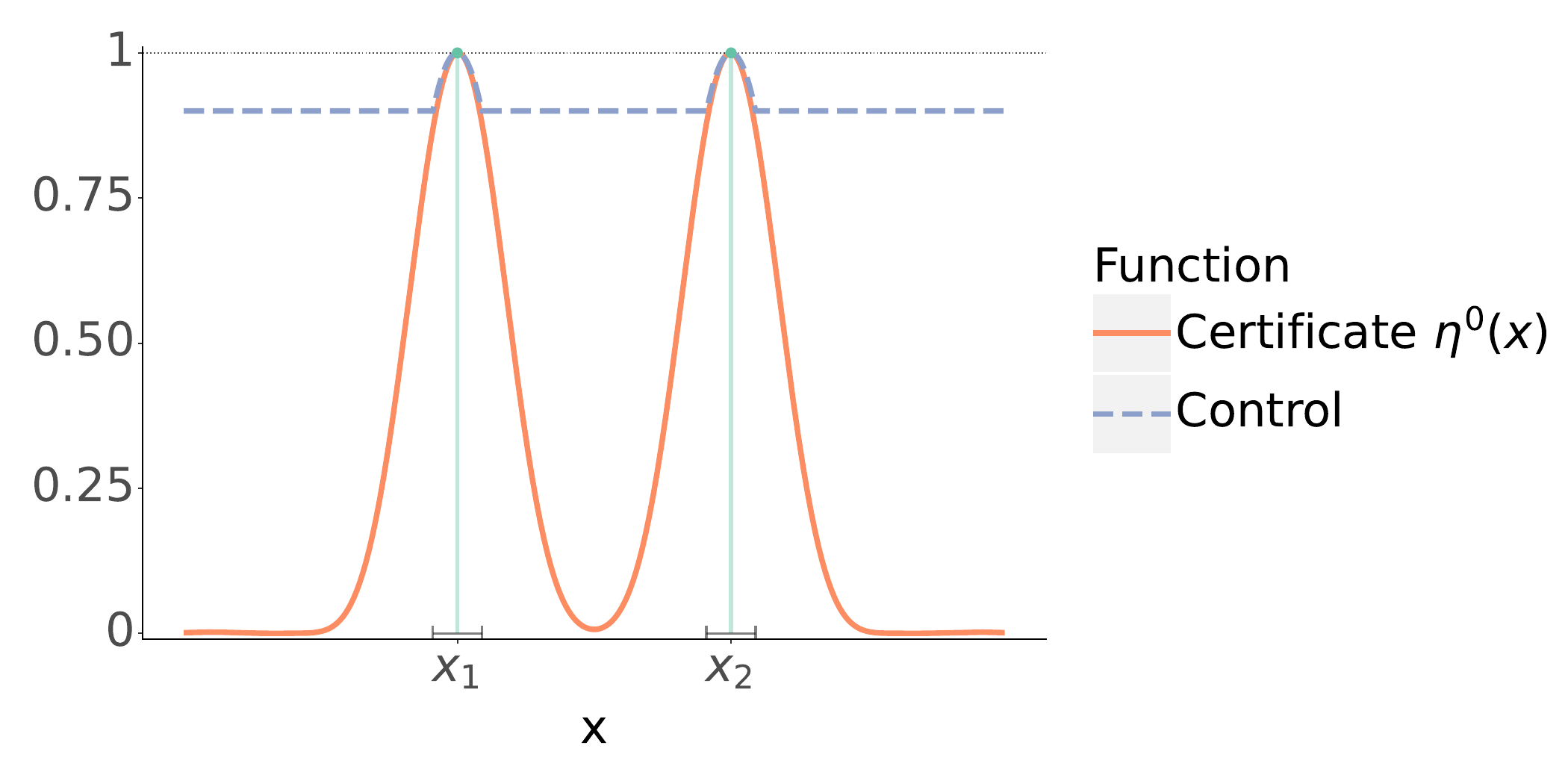} \hfill
	\includegraphics[width=.49\linewidth]{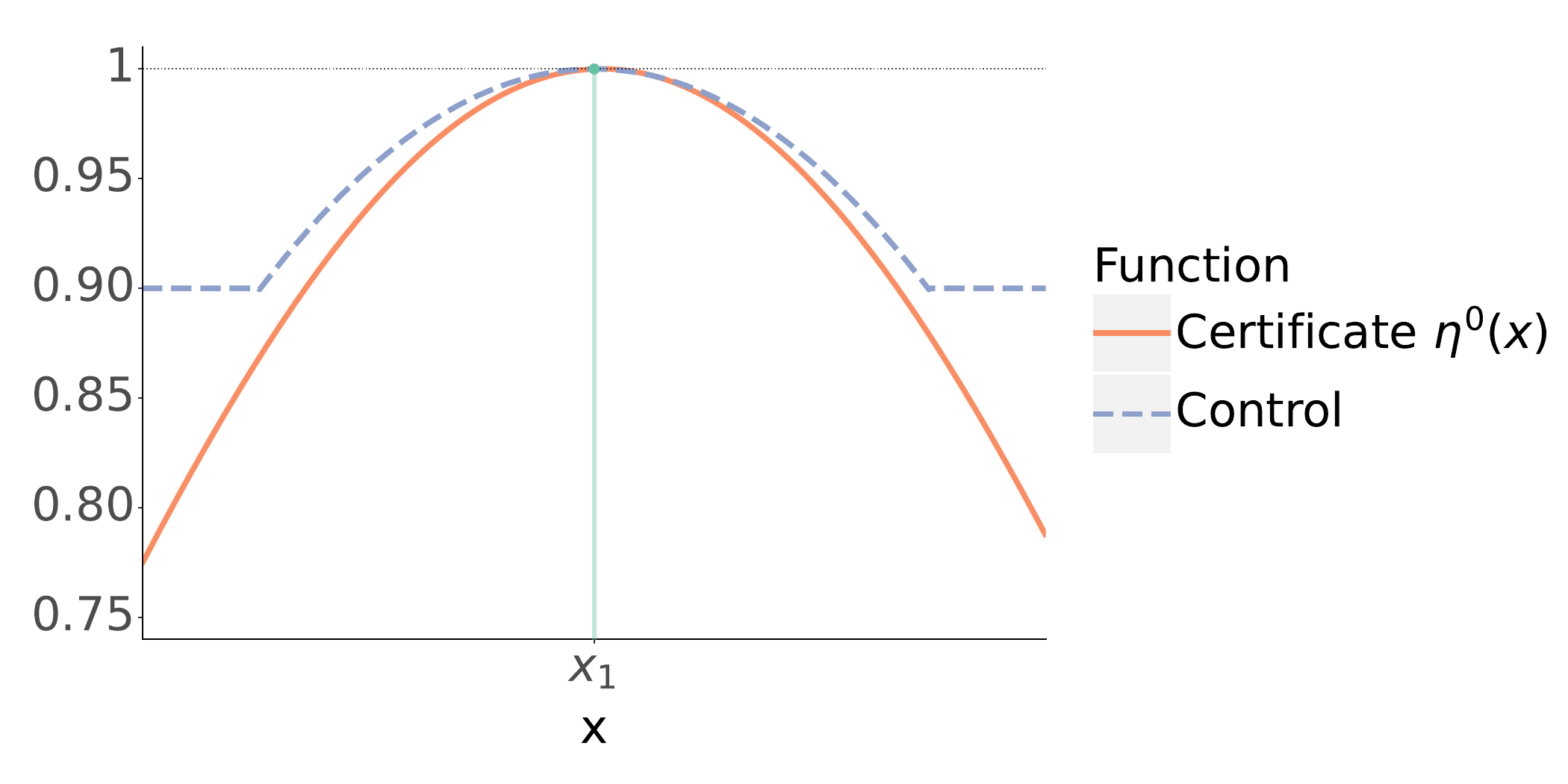} \hfill
	\caption{On the left, a one-dimensional illustration of the non-degenerate dual certificate $\Certificate$ in the near regions and far regions for some target $\target$ with two spikes $x_1$ and $x_2$, shown in light green. On the right, a zoomed location of the curvature control in the near region $\nearregion_1(r)$ around $t_1$. The certificates are drawn in orange while blue dotted lines illustrates the required control.}
	\label{fig:DualCerticates}
\end{figure}

The crux of the statistical analysis lies in the construction of a \textit{non-degenerate dual certificate} $\Certificate$, \textit{i.e.} a sub-gradient of the TV-norm interpolating $\textnormal{sign}(a_k^0) = \pm 1$ at the support point $t^k$ of $\target$, and with controlled behaviors in the near and far regions, see Figure~\ref{fig:DualCerticates}. Both existence and construction of such objects have been the focus of a lot of attention in the continuous sparse recovery literature, see for example \parencite{candes2014towards,poon2023geometry,decastro2019sparse}.

\begin{definition}[Pivot non-degenerate certificates]
Let $K_{\pivot}(\cdot,\cdot)\,:\,\Param\times\Param\to\R$ be a pivot kernel and let $\target\in\Model_{\K, \Delta_0, \distance}$ where the Fisher-Rao distance is given by $K_\pivot(\cdot,\cdot)$.  

We say that $\eta^0$ is an $(\varepsilon_0,\varepsilon_2,r_0)$-pivot non-degenerate dual certificate if
\begin{itemize}
    \item There exists $c^0\in\mathcal{F}$ such that $\eta^0=F^\star c^0$,
    \item For all $i\in[\K]$, $\eta^0(t^0_i)=\mathrm{sign}(a^0_i)$,
    \item For all $x\in\farregion(r_0)$, $|\eta^0(x)|\leq1-\varepsilon_0$,
    \item For all $i\in[\K]$, for all $x\in\nearregion_i(r_0)$, $|\eta^0(x)|\leq1-\varepsilon_2\distance(x,t_i^0)^2$,
\end{itemize}
where $\distance$ is the Fisher-Rao distance associated to $K_{\pivot}(\cdot,\cdot)$, and the near and far regions are those of~$(t_i^0)_{i=1}^\K$. 

For all $i\in[\K]$, we say that $\eta^0_i$ is an $(\varepsilon_0,\varepsilon_2,r_0)$-pivot non-degenerate localizing certificate at point $t_i^0$ if
\begin{itemize}
    \item There exists $c^0_i\in\mathcal{F}$ such that $\eta^0_i=F^\star c^0_i$,
    \item For all $x\in\farregion(r_0)$, $|\eta_i^0(x)|\leq1-\varepsilon_0$,
    \item For all $j\in[\K]$, for all $x\in\nearregion_j(r_0)$, $|\delta_{ij}-\eta_i^0(x)|\leq{\varepsilon_2}\distance(x,t_j^0)^2$,
\end{itemize}
where $\delta_{ij} = 1$ if $i=j$ and $0$ otherwise.
\end{definition}

\medskip

\begin{theorem}[Existence of pivot non-degenerate dual certificates]
\label{thm:existence_certificates}
Let $\Param$ be a compact set of~$\R^d$, let $\mathcal{F}$ be a separable Hilbert space, and let $K_\pivot(\cdot,\cdot)$ be a positive type kernel that satisfies Assumption~\ref{hyp:local_positive_curvature} with parameters $\K$, $\Delta_0$, $r_0$, $\bar\varepsilon_0$ and $\bar\varepsilon_2$, and Assumption~\ref{hyp:pivot_kernel} with constant $C_{\switch}>0$ with respect to some model kernel $K_\model(\cdot,\cdot)$ which satisfies~\eqref{hyp:continuous_model_kernel}.

Then there exist an $({\bar\varepsilon_0},\bar\varepsilon_2,r_0)$-pivot non-degenerate dual certificate such that $\|c^0\|_{\mathcal{F}}\leq \sqrt{2}C_{\switch}\sqrt{\K}$ and $(\bar\varepsilon_0,\bar\varepsilon_2,r_0)$-pivot non-degenerate localizing certificates at points $t_i^0$ such that $\|c^0_i\|_{\mathcal{F}}\leq \sqrt 2C_{\switch}$.
\end{theorem}

\newpage 

\begin{proof}
    This proof is an adaptation of the proof of \cite[Theorem~2, Page 266]{poon2023geometry} for the pivot kernel. Denote by $f_\pivot(x)\defeq K_{\pivot}(x,\cdot)$ the canonical feature map of $\mathbb H_\pivot$. 
    
    Consider some coefficients $(\alpha_{1},\alpha_{2})\in\mathbb R^{\K}\times\mathbb R^{d\K}$ such that
    \begin{equation}
    \label{eq:certificate_proof}
        \eta^0(x)\defeq \sum_{j=1}^\K \alpha_{1,j}K_{\pivot}(\trueparam_j,x)+\langle\alpha_{2,j},\nabla_1K_{\pivot}(\trueparam_j,x)\rangle_{\mathbb R^d}
    \end{equation}
    satisfies $\nabla \eta^0(\trueparam_j)=0$ and $\eta^0(\trueparam_j)=\mathrm{sign}(a_j^0)$, for $j\in[\K]$. These $\K(d+1)$ constraints can be written as a linear system 
    \[
    \Upsilon
    \begin{pmatrix}
            \alpha_1 \\ 
            \alpha_2
    \end{pmatrix}
    =
    \begin{pmatrix}
        (\mathrm{sign}(a_j^0))_{j=1}^\K \\ 
        \mathbf{0}_{\K d}
    \end{pmatrix}=:\mathbf{u},
    \]
    where $\Upsilon\in\mathbb R^{d\K \times d\K}$ is a real symmetric matrix defined as
    \begin{align*}
    \Upsilon
        &=
    \begin{pmatrix}
            (K_{\pivot}(\trueparam_i,\trueparam_j))_{1\leq i,j\leq \K} 
                & (\nabla_2 K_{\pivot}(\trueparam_i,\trueparam_j)^\top)_{1\leq i,j\leq \K} 
            \\ 
            (\nabla_1 K_{\pivot}(\trueparam_i,\trueparam_j))_{1\leq i,j\leq \K} 
                & (\nabla_1\nabla_2 K_{\pivot}(\trueparam_i,\trueparam_j))_{1\leq i,j\leq \K}
    \end{pmatrix}.
    \end{align*}
    
    Then, define
    \[
        \mathbf{f}(x)\defeq
        \Big(
            (K_{\pivot}(\trueparam_i,x))_{i=1}^\K,
            \big(\nabla_1K_{\pivot}(\trueparam_i,x)^\top\big)_{i=1}^\K
        \Big)^\top\in\mathbb{R}^{\K(d+1)}
    \]
    Since the mixed partial derivative $\partial_i \partial_j K_\pivot$ exists and are continuous, \textcite[Lemma 4.34]{steinwart2008support} ensures that the feature map $f_\pivot: \mathbb{R}^d \to \Hilbert_\pivot$ is continuously differentiable and we denote by~${\nabla f_\pivot(\bs) = (\partial_i f_\pivot(\bs))_{i=1}^d \in \Hilbert_{\pivot}^{d}}$ the vector of partial derivatives, the latter verifying:~${\<\partial_i f_\pivot(\bs),  \partial_j f_\pivot(\bt) \>_{\Hilbert_\pivot} = \partial_i \partial_j K_\pivot(\bs, \bt)}$.

    We can thus define $\mathbf{f}=
        \Big(
            (f_{\pivot}(\trueparam_i))_{i=1}^\K,
            \big(\nabla f_{\pivot}(\trueparam_i)^\top\big)_{i=1}^\K
        \Big)^\top\in (\mathbb H_\pivot)^{\K(d+1)}$ and we introduce the following notation:
          \[
            \forall p\geq 1,\,\forall q\geq 1,\, \forall (a_i)_{i=1}^p\in (\mathbb H_\pivot)^p,\, \forall (b_j)_{j=1}^q\in (\mathbb H_\pivot)^q\,,\quad 
            \langle a\otimes b\rangle_{\mathbb H_\pivot}\defeq (\langle a_i,b_j\rangle_{\mathbb H_\pivot})_{i\in[p],j\in[q]}\,.
        \]      
    Note that
        \begin{align*}
    \Upsilon
        &=
    \begin{pmatrix}
            (\langle f_\pivot(\trueparam_i),f_\pivot(\trueparam_j) \rangle_{\mathbb H_\pivot})_{1\leq i,j\leq \K} 
                & (\langle  f_\pivot(\trueparam_i)\otimes\nabla f_\pivot(\trueparam_j) \rangle_{\mathbb H_\pivot})_{1\leq i,j\leq \K} 
            \\ 
            (\langle \nabla f_\pivot(\trueparam_i)\otimes f_\pivot(\trueparam_j) \rangle_{\mathbb H_\pivot})_{1\leq i,j\leq \K} 
                & (\langle \nabla f_\pivot(\trueparam_i)\otimes\nabla f_\pivot(\trueparam_j) \rangle_{\mathbb H_\pivot})_{1\leq i,j\leq \K}
    \end{pmatrix}\\
        \mathbf{f}(x)&=\langle \mathbf{f},f_{\pivot}(x)\rangle_{\mathbb H_\pivot}\,.
    \end{align*}

    We deduce that $\Upsilon$ is the Gram matrix of $\mathbf{f}$, with respect to the dot product of $\mathbb H_\pivot$, namely $\Upsilon=\langle \mathbf{f}\otimes\mathbf{f}\rangle_{\mathbb H_\pivot}$.

    Assuming that $\Upsilon$ is invertible, we can rewrite~\eqref{eq:certificate_proof} as
    \begin{subequations}
            \begin{align}
                \eta^0(x)     &= \big(\Upsilon^{-1}\mathbf{u}\big)^{\top}\mathbf{f}(x)\\
                \eta^0(x)        &= \langle \mathbf{f}^\top \Upsilon^{-1}\mathbf{u},f_{\pivot}(x)\rangle_{\mathbb H_\pivot}
            \end{align}
    and we deduce that $\eta=\mathbf{f}^\top \Upsilon^{-1}\mathbf{u}$ where the equality holds in $\mathbb H_\pivot$, by the reproducing property. 

    Hence
    \begin{align*}
        \|\eta^0\|_{\mathbb H_\pivot}^2&=\langle\mathbf{f}^\top \Upsilon^{-1}\mathbf{u},\mathbf{f}^\top \Upsilon^{-1}\mathbf{u} \rangle_{\mathbb H_\pivot}=\mathbf{u}^\top \Upsilon^{-1} \mathbf{u}\,,
    \end{align*}
    using that $\langle \mathbf{f}\otimes\mathbf{f}\rangle_{\mathbb H_\pivot}=\Upsilon$.

    We also define the block diagonal normalisation matrix  $D_{\mathfrak g}\in \R^{\K(d+1)\times \K(d+1)}$ as
\begin{equation}
    \label{eq:matrix-D-block}
    D_{\mathfrak g} \defeq 
        \begin{pmatrix}
            \Id_{s} \\
                & {\mathfrak g}_{\bx_1}^{-\frac12} \\
            &&\ddots\\
            &&& {\mathfrak g}_{\bx_s}^{-\frac12}
\end{pmatrix}
\end{equation}
and $\tilde \Upsilon \defeq D_{\mathfrak g} \Upsilon D_{\mathfrak g}$ which has constant value $1$ along its diagonal. 

We deduce that 
\begin{equation}
    \label{eq:bound_et0_etaMat}
    \|\eta^0\|_{\mathbb H_\pivot}^2\leq 
        \big\|\tilde \Upsilon^{-1}\big\|
        \|D_{\mathfrak g} \mathbf{u}\|_2^2\leq 2\K\,.
\end{equation}
using that $\|\tilde \Upsilon^{-1}\|\leq 2$ by \cite[Lemma 3]{poon2023geometry}.

Invoke Proposition~\ref{prop:model_description}, \eqref{eq:bound_et0_etaMat} and \eqref{eq:constant_switch} to get that
\begin{equation}
    \label{eq:bound_switch}
    \|c_0\|_{\mathcal{F}}=\|\eta^0\|_{\mathbb H_\model}\leq C_{\switch}\|\eta^0\|_{\mathbb H_\pivot}\leq C_{\switch}\sqrt{2\K}
\end{equation}
\end{subequations}

By \cite[Theorem~2]{poon2023geometry}, $\eta^0$ is a $({\varepsilon_0},\varepsilon_2,r_0)$-pivot non-degenerate dual certificate (note that we slightly changed the values of ${\varepsilon_0},\varepsilon_2$ of \cite[Assumption]{poon2023geometry} in Assumption~\ref{hyp:local_positive_curvature}).

The existence of $(\varepsilon_0,\varepsilon_2,r_0)$-pivot non-degenerate localizing certificates follow the same line. In this case $\mathbf{u}^\top=((\mathrm{1}_{\{i=j\}}\mathrm{sign}(a_j^0))_{j\in[\K]}^\top,\mathbf{0}_{\K d}^T)$ and hence $\|D_{\mathfrak g} \mathbf{u}\|_2=1$.
\end{proof}

\subsection{Existence of sketch pivot non-degenerate dual certificates}

\begin{definition}[Sketch pivot non-degenerate certificates]
Let $\psi_\bomega\,:\,\Param\to\mathbb{C}$ be a function and let $\SketchDist$ be a law on~$\mathbb R^d$. Consider 
\begin{align}
    {K}_{\pivot}(\bs, \bt)&= \Expectation_{\bomega\sim \SketchDist}\big[\psi_\bomega(s)\overline{\psi}_\bomega(t)\big]\,,
    \tag{$\mathcal H_{\psi_\bomega,\SketchDist}$}
\end{align}
and let $\target\in\Model_{\K, \Delta_0, \distance}$ where the Fisher-Rao distance is given by $K_\pivot(\cdot,\cdot)$. Let $m\geq1$ and let $(\bomega_i)_i$ be {\rm i.i.d.} vectors drawn with respect to $\SketchDist$.

We say that $\eta^0$ is an $(\varepsilon_0,\varepsilon_2,r_0)$-sketch pivot non-degenerate dual certificate if
\begin{itemize}
    \item There exists $c^0\in\mathbb{C}^m$ such that $\displaystyle\eta^0=\frac{1}{\sqrt{m}} \sum_{j=1}^mc^0_j \varphi_{\bomega_j}$,
    \item For all $i\in[\K]$, $\eta^0(t^0_i)=\mathrm{sign}(a^0_i)$,
    \item For all $x\in\farregion(r_0)$, $|\eta^0(x)|\leq1-\varepsilon_0$,
    \item For all $i\in[\K]$, for all $x\in\nearregion_i(r_0)$, $|\eta^0(x)|\leq1-\varepsilon_2\distance(x,t_i^0)^2$,
\end{itemize}
where $\distance$ is the Fisher-Rao distance associated to $K_{\pivot}(\cdot,\cdot)$, and the near and far regions are those of~$(t_i^0)_{i=1}^\K$. 

For all $i\in[\K]$, we say that $\eta^0_i$ is an $(\varepsilon_0,\varepsilon_2,r_0)$-pivot non-degenerate localizing certificate at point $t_i^0$ if
\begin{itemize}
    \item There exists $c_i^0\in\mathbb{C}^m$ such that $\displaystyle\eta^0=\frac{1}{\sqrt{m}} \sum_{j=1}^mc^0_{i,j} \varphi_{\bomega_j}$,
    \item For all $x\in\farregion(r_0)$, $|\eta_i^0(x)|\leq1-\varepsilon_0$,
    \item For all $j\in[\K]$, for all $x\in\nearregion_j(r_0)$, $|\mathrm{1}_{\{i=j\}}-\eta_i^0(x)|\leq{\varepsilon_2}\distance(x,t_j^0)^2$,
\end{itemize}
where $\mathrm{1}_{\{i=j\}}$ equals one if $i=j$ and $0$ otherwise.
\end{definition}

\medskip
\begin{theorem}[Existence of sketch pivot non-degenerate dual certificates]
\label{thm:existence_sketch_certificates}
Let $\Param$ be a compact set of $\R^d$. Let $\varphi_\bomega,\psi_\bomega\,:\,\Param\to\mathbb{C}$ be functions and let $\SketchDist$ be a law on $\mathbb R^d$. Consider 
\begin{align}
    {K}_{\model}(\bs, \bt)&= \Expectation_{\bomega\sim \SketchDist}\big[\varphi_\bomega(s)\overline{\varphi}_\bomega(t)\big]\,,
    \notag
    \\
    {K}_{\pivot}(\bs, \bt)&= \Expectation_{\bomega\sim \SketchDist}\big[\psi_\bomega(s)\overline{\psi}_\bomega(t)\big]\,.
    \tag{$\mathcal H_{\psi_\bomega,\SketchDist}$}
\end{align}
Assume that $\psi_\bomega$ satisfies Assumption~\ref{hyp:sketch_function_bounded}, and assume that $K_\pivot(\cdot,\cdot)$ satisfies Assumption~\ref{hyp:local_positive_curvature} $($with parameters $\K$, $\Delta_0$, $r_0$, $\bar\varepsilon_0$ and $\bar\varepsilon_2)$ and Assumption~\ref{hyp:pivot_kernel} with constant $C_{\switch}>0$.  

Let $\alpha\in(0,1)$ and $m$ be such that
\begin{equation}
    m\geq C_{\textnormal{sketch}}\,\max(d,\log(\K))\,\K\,\log\bigg(\frac{\max(1,|\Param|)\,\K}{\alpha}\bigg)\,,
\end{equation}
where $C_{\textnormal{sketch}}>0$ depends on ${K}_{\pivot}$ and  polynomially on $d$, $\displaystyle|\Param|= \sup_{\bs,\bt\in\Param}\distance(\bs, \bt)$ and the Fisher-Rao distance $\distance$ is given by $K_\pivot(\cdot,\cdot)$. 

Then, there exists a constant $C'_{\pivot}>0$ which depends only on the kernel $K_\pivot$ such that, with probability at least $1-\alpha$, for all $m\geq m_0$, there exist a $({\bar\varepsilon_0/4},3\bar\varepsilon_2/2,r_0)$-sketch pivot non-degenerate dual certificate such that $\|c^0\|_{\mathbb C^m}\leq {C'_{\pivot}}C_{\switch}\sqrt{\K}$ and $(\bar\varepsilon_0/4,3\bar\varepsilon_2/2,r_0)$-sketch pivot non-degenerate localizing certificates at points~$t_i^0$ such that $\|c^0_i\|_{\mathbb C^m}\leq C'_{\pivot}C_{\switch}$.
\end{theorem}

\begin{proof}
    This result is a consequence of  \cite[Theorem~4, Page 274]{poon2023geometry} and Theorem~\ref{thm:existence_certificates}. The only subtlety relies in \cite[Section 6.5, Step~3, Page 284]{poon2023geometry} where one has to invoke Proposition~\ref{prop:model_description} and \eqref{eq:constant_switch} to get the bound $\|c^0\|_{\mathbb C^m}\leq {C'_{\pivot}}C_{\switch}\sqrt{\K}$. Except this latter point, the reasoning is identical.
\end{proof}

\subsection{Proof of Theorem~\ref{thm:sinc4_lca}}
\label{proof:sinc4_lca}

Consider the sinc-4 kernel $\Psi_\bandwidth$ as defined by~\eqref{eq:sketch_kernel_sinc}. We want to prove that this kernel satisfies the local positive curvature assumption, Assumption~\ref{hyp:local_positive_curvature}. We split the proof into $5$ parts.

\subsubsection{Useful geometric object}
For $\bandwidth > 0$ and $\bx, \by \in \mathbb{R}^d$, we define:
	\begin{itemize}
        \item The random features $\psi_\bomega(\bx) = \sqrt{\frac{f^{(4)}_{\bandwidth}}{\SketchDist}(\bomega)} e^{-\imath \bomega^\top \bx}$ with the spectral measure $f^{(4)}_{\bandwidth} = (2\pi)^{-d} \fourier[\Psi_\bandwidth]$ detailed in~\Cref{app:sinc4_RKHS}
		\item The kernel $K_\bandwidth(\bx, \by) = \Expectation_{\bomega \sim \SketchDist_\bandwidth} \left[ \overline{\psi_\bomega(\bx)} \psi_\bomega(\by) \right] = \Psi_\bandwidth(\bx - \by ) = \Psi_1(\frac{\bx - \by}{\bandwidth}) = K_1(\frac{\bx}{\bandwidth}, \frac{\by}{\bandwidth})$.
		\item The metric $\metric_{\bx, \bandwidth} = \nabla_{1} \nabla_{2} K_\bandwidth(\bx, \bx) = - \nabla^2 \Psi_{\bandwidth}(\bm{0}) =  \frac{1}{12 \bandwidth^2} \Id $. It is independent of the point $\bx$ in the case of translation invariant kernel.
        \item The norm at point $\bx$: $\Vert \bz \Vert_{\bx, \bandwidth}^2 = \bz^\top \metric_{\bx, \bandwidth} \bz = \frac{1}{12 \bandwidth^2} \Vert \bz \Vert_{2}^2 $
        \item The geodesic distance from $\bx$ to $\by$ which is proportional to the euclidean norm since the kernel is translation invariant: 
        \[
        \mathfrak{d}_{\metric, \bandwidth}(\bx, \by) \defeq \inf \left\{ \int_0^1 \Vert g'(t) \Vert_{g(t)} \; : \; g(0) = \bx, \; g(1) = \by \right\} = \frac{1}{2\sqrt{3} \bandwidth} \Vert \bx - \by \Vert_2 
        \]
	\end{itemize}
	Moreover, the metric have the following re-scaling properties
	\begin{align}
	    \metric_{\bx, \bandwidth} & = \frac{1}{\bandwidth^2} \metric_{\bu, 1}, \; \forall \bx, \bu \in \R^d, \\
        \mathfrak{d}_{\metric, \bandwidth}(\bx, \by) & = \mathfrak{d}_{\metric, 1}\left(\frac{\bx}{\bandwidth}, \frac{\by}{\bandwidth}\right).
	\end{align}
	Finally, notice that this choice of $K_\bandwidth$ defines a flat euclidean geometry so that euclidean and Riemannian gradients (and Hessians) coincide. With a slight abuse of notations, we will use $\Psi_\bandwidth(\bx, \by)$ in place of $K_\bandwidth(\bx, \by) = \Psi_\bandwidth(\bx - \by)$.

\subsubsection{Reduction to bandwidth one}
\textbf{ $\bullet$ Covariant derivatives:} First, we begin by proving that it suffices to prove the result for $\bandwidth=1$. The controls we look after, see Equations~\eqref{eqs:ass_positive_curvarure}, involve the covariant derivatives of the kernel of order $(i,j)$ seen as multi-linear maps $\Psi^{(ij)}(\bx, \by) : \left( \mathbb{C}^d\right)^{i+j} \to \mathbb{C}$:
\begin{subequations}
\begin{equation}
    \label{eq:kernelcovariantderivative}
    \Psi^{(ij)}(\bx, \by)[\bq_1, \ldots, \bq_{i+j}] \defeq \Expectation_{\bomega} \left[ \overline{D_i [\psi_\bomega] (\bx) [\bq_1, \ldots, \bq_i]} D_j [\psi_\bomega] (\by) [\bq_{i+1}, \ldots, \bq_{i+j}] \right]
\end{equation}
In the particular case of $\psi_\bomega$ yielding the $\Psi_\bandwidth$ kernel, we have simple expression of the partial derivatives of order $r$ as $\partial_{i_1, \ldots, i_r} \psi_\bomega(\bx) = (-\imath)^{r} \omega_{i_1} \ldots \omega_{i_r} \psi_\bomega(\bx)$ so that the differential of order $r$ can be written as a rank-$1$ tensor in the canonical basis of $(\mathbb{C}^{d})^r$
\begin{equation}
 D_r[\psi_\bomega](\bx)[\cdot] = (- \imath)^r  \psi_\bomega(\bx) (\bomega \otimes \ldots \otimes \bomega)[\cdot] = (- \imath)^r  \psi_\bomega(\bx) \, \bomega^{\otimes r} [\cdot] .   
\end{equation}
This leads to the following simplification of \Cref{eq:kernelcovariantderivative} putting $Q = (\bq_1, \ldots, \bq_{i+j})$
\begin{equation}
    \Psi^{(ij)}(\bx, \by)[\bq_1, \ldots, \bq_{i+j}] =  (- 1)^{j} \mathbb{E}_\bomega \left[ \overline{D_{i+j}[\psi_\bomega](\bx)[Q]} \psi_\bomega(\by) \right] 
\end{equation}

The operator norm at $\bx, \by$ is then defined as \textcite[][Eq. (26)]{poon2023geometry} 
\begin{equation}
    \label{eq:OperatorNorm}
    \Vert \Psi^{(ij)}(\bx, \by)\Vert_{\bx, \by} \coloneqq \sup_{\substack{\Vert \bq_l \Vert_{\bx} \leq 1, \, l=1,\ldots,i \\\Vert \bq_{l+i} \Vert_{\by} \leq 1, \, l=1,\ldots,j}} \left \vert \mathbb{E}_\bomega \left[ \overline{D_i[\psi_\bomega](\bx)[\bq_1, \ldots, \bq_i]} D_j[\psi_\bomega](\by)[\bq_{i+1}, \ldots, \bq_{i+j}] \right] \right\vert .
\end{equation}
\end{subequations}

The first fact is that $\Psi_\bandwidth$ and $\Psi_1$ (and their covariant derivatives) have the same operator norm provided we rescale the coordinate by $1 / \bandwidth$, as stated by the following lemma.
\begin{lemma}[Operator norm invariance by scaling $\bx$ onto $\bx / \bandwidth$] 
\label{lemma:norm_tau_equal_norm_one} 
Let $\bandwidth > 0$ and $\bx, \by \in \mathcal{X}$, then, putting $\bu = \bx / \bandwidth $ and $\bv = \by / \bandwidth$ for all $0 \leq i, j \leq 2$ we have
    \[
        \Vert K_\bandwidth^{(ij)}(\bx, \by)\Vert_{\bx, \by} = \Vert K_1^{(ij)}(\bu, \bv)\Vert_{\bu, \bv}
    \]
\end{lemma}

\begin{proof} 

We begin by noting that the norm $\Vert \bz \Vert_{\bx} = \Vert \bz \Vert_{\by} \propto \frac{1}{\bandwidth} \Vert \bz \Vert_2$ is a rescaled euclidean metric, and is thus independent of the points $\bx$ or $\by$. Thus, the operator norm defined in \Cref{eq:OperatorNorm} can be written as
\begin{align*}
    \Vert K_\bandwidth^{(ij)}(\bx, \by)\Vert_{\bx, \by} =& \sup_{\substack{\Vert \bq_l \Vert_\bandwidth \leq 1, \\ l=1,\ldots,i+j}} \left \vert \mathbb{E}_\bomega \left[ \overline{D_i[\psi_\bomega](\bx)[\bq_1, \ldots, \bq_i]} D_j[\psi_\bomega](\by)[\bq_{i+1}, \ldots, \bq_{i+j}] \right] \right\vert ,\\
    =&\sup_{\substack{\Vert \bq_l \Vert_\bandwidth \leq 1, \\ l=1,\ldots,i+j}}  \left \vert (- 1)^{j} \mathbb{E}_\bomega \left[ \overline{D_{i+j}[\psi_\bomega](\bx)[Q]} \psi_\bomega(\by) \right] \right\vert 
\end{align*}

Using this symmetry, it suffice to prove the desired identity for $j=0$ and $i=0, \ldots, 4$
    \begin{itemize}
        \item  $i=j=0$,  $\Vert K_\bandwidth^{(00)}(\bx, \by)\Vert_{\bx} =  K_\bandwidth(\bx, \by) = K_1(\frac{\bx}{\bandwidth}, \frac{\by}{\bandwidth}) = K_1(\bu, \bv)$.
        \item $i=1, j=0$: using that $\bandwidth^2 \metric_{\bx, \bandwidth} = \metric_{\bu, 1}$ when $\bu = \bx / \bandwidth$  
        \begin{align*}
            \Vert K_\bandwidth^{(10)}(\bx, \by)\Vert_{\bx} &=  \Vert \metric_{\bx, \bandwidth}^{-1/2} \nabla_1 K_\bandwidth(\bx, \by) \Vert_2\\
            &= \Vert \metric_{\bx, \bandwidth}^{-1/2} \frac{1}{\bandwidth} \nabla_1 K_1(\bu, \bv) \Vert_2\\ 
            &= \Vert \metric_{\bu, 1}^{-1/2}\nabla_1 K_1(\bu, \bv) \Vert_2\\
            &= \Vert K_1^{(10)}(\bu, \bv)\Vert_{\bu}\,.
        \end{align*}
        \item $i=2, j=0$: similarly
        \begin{align*}
           \Vert K_\bandwidth^{(20)}(\bx, \by)\Vert_{\bx} &=  \Vert \metric_{\bx, \bandwidth}^{-1/2} \nabla^2 \Psi_\bandwidth(\bx - \by) \metric_{\bx, \bandwidth}^{-1/2} \Vert_2 ,\\
            &= \Vert \metric_{\bx, \bandwidth}^{-1/2} \frac{1}{\bandwidth^2}  \nabla^2 \Psi_1(\bu - \bv) \metric_{\bx, \bandwidth}^{-1/2}  \Vert_2 , \\
            &=  \Vert \metric_{\bu, 1}^{-1/2} \nabla^2 \Psi_1(\bu -\bv) \metric_{\bu, 1}^{-1/2} \Vert_2 , \\
            &=  \Vert K_1^{(20)}(\bu, \bv)\Vert_{\bu}.
        \end{align*}
    \end{itemize}
    The same follows for higher order derivatives.  
\end{proof}

\noindent
\textbf{ $\bullet$ Local curvature:} \textcite{poon2023geometry} define geometric quantities related to the local curvature of the kernel~$K_\bandwidth$ (and hence of $\Certificate$) on near and far regions. 

\begin{definition}[$\bar{\varepsilon}_0$ and $\bar{\varepsilon_2}$] 
Let $r \geq 0$,
    \begin{align*}
        \bar{\varepsilon}_0(r, \bandwidth) &\defeq  \sup \left\{\varepsilon \, ; \, K_\bandwidth(\bx, \by) \leq 1 - \varepsilon, \; \forall \bx, \by \textnormal{ s.t. } \mathfrak{d}_{\metric, \bandwidth}(\bx, \by) \geq r \right\}, \\
        \bar{\varepsilon}_2(r, \bandwidth) &\defeq  \sup \left\{\varepsilon \, ; \, \Psi^{(02)}_\bandwidth(\bx, \by)[\bz, \bz] \geq \varepsilon \Vert \bz \Vert_{\bx, \bandwidth}^2 , \; \forall z \in \R^d, \; \forall \bx, \by \textnormal{ s.t. } \mathfrak{d}_{\metric, \bandwidth}(\bx, \by) < r \right\}\,.
    \end{align*}
\end{definition}
These controls depend both on the kernel and the radius $r = r_{near}$ of the regions considered. However, the following Lemmas show that this controls do not depend on the bandwidth $\tau$.

\begin{lemma}[Isometry between $K_\bandwidth$ and $K_1$]
	\label{lemma:isometry}
	Let $\bandwidth > 0$ and
	\[\mathcal{H}_\bandwidth \defeq \left\{ g \in L^2(\mathbb{R}^d) \, : \Vert g \Vert_{\mathcal{H}_\bandwidth}^2 \defeq \frac{1}{(2\pi)^d} \int \frac{\left\lvert \fourier[g] \right\rvert^2}{ \fourier[\Psi_{\bandwidth}]} < + \infty \right\},
	\]
	then the rescaling map $\textrm{iso}_\bandwidth : g \in \mathcal{H}_1 \to  g (\frac{\cdot}{\bandwidth}) \in \mathcal{H}_{\bandwidth}$ is an isometry between $(\mathcal{H}_{1}, \Vert \cdot \Vert_{\bandwidth = 1} )$ and $(\mathcal{H}_{\bandwidth}, \Vert \cdot \Vert_{\bandwidth})$.
	Moreover, for any $\bx, \by \in \mathbb{R}^d$, we have:
	\begin{enumerate}[label={\sf $($\roman*$)$}]
		\item $ K_1(\frac{\bx}{\bandwidth}, \frac{\bx}{\bandwidth}) = K_\bandwidth(\bx, \by) $
		\item $\nabla_1 K_1(\frac{\bx}{\bandwidth}, \frac{\by}{\bandwidth}) = \bandwidth \nabla_1 K_\bandwidth(\bx, \by)$
		\item $\nabla_{12} K_1(\frac{\bx}{\bandwidth}, \frac{\by}{\bandwidth}) = \bandwidth^2 \nabla_{12} K_\bandwidth(\bx, \by)$
	\end{enumerate}
\end{lemma}
\begin{proof}
	Let $\bandwidth > 0$,
	\begin{enumerate}[label={\sf $($\roman*$)$}]
	\item By definition of $K_{\bandwidth}$ and $\Psi_{\bandwidth}$.
	\item $\nabla_1 K_1(\frac{\bx}{\bandwidth}, \frac{\by}{\bandwidth}) = \nabla \Psi_1(\frac{\bx - \by}{\bandwidth}) = \bandwidth  \frac{1}{\bandwidth} \nabla \Psi_1(\frac{\bx - \by}{\bandwidth}) = \bandwidth \nabla \Psi_\bandwidth(\bx - \by) = \bandwidth \nabla_1 K_\bandwidth(\bx, \by)$
	\item $\nabla_{12} K_1(\frac{\bx}{\bandwidth}, \frac{\by}{\bandwidth}) = - \nabla^2 \Psi_1(\frac{\bx - \by}{\bandwidth}) =  - \bandwidth^2 \frac{1}{\bandwidth^2}\nabla^2 \Psi_1(\frac{\bx - \by}{\bandwidth})= - \bandwidth^2 \nabla^2 \Psi_\bandwidth(\bx - \by)  = + \bandwidth^2 \nabla_{12} K_\bandwidth(\bx, \by)$
	\end{enumerate}

	For the isometry result, we simply use the change of variable $\bomega' = \bomega / \bandwidth$ in the RKHS norm definition, with diagonal Jacobian yielding $\diff \bomega' = \bandwidth^{-d} \diff \bomega$. 
	\begin{align*}
		\Vert g \Vert_{\mathcal{H}_1} & = \frac{1}{(2 \pi)^d}  \int \frac{\left\lvert \fourier[g] (\bomega) \right\rvert^2}{ \fourier[\Psi_{1}](\bomega)}  \diff \bomega , \\
		&=  \frac{1}{(2 \pi)^d}  \int  \frac{\left\lvert \fourier[g] (\bandwidth \bomega')\right\rvert^2}{ \fourier[\Psi_{1}](\bandwidth \bomega')}  \bandwidth^d \diff \bomega'  , & \textrm{(change of variables : } \bomega = \bandwidth \bomega' \textrm{)}\\
		& =  \frac{1}{(2 \pi)^d}  \int \frac{ \bandwidth^{-2d}  \left\lvert \fourier[g(\cdot / \bandwidth)] ( \bomega') \right\rvert^2}{\fourier[\Psi_{1}](\bomega')}  \bandwidth^d  \diff \bomega'  , & \textrm{(} \fourier[g](\bandwidth \bomega') = \bandwidth^{-d} \fourier[g(\cdot / \bandwidth)](\bomega') \textrm{)}, \\
		& =  \frac{1}{(2 \pi)^d}  \int\frac{\left\lvert \fourier[ g(\cdot / \bandwidth)] ( \bomega') \right\rvert^2}{ \bandwidth^{d} \fourier[\Psi_{1}](\bandwidth \bomega')}  \diff \bomega'  , &  \\
		& =  \frac{1}{(2 \pi)^d}  \int\frac{\left\lvert \fourier[ g(\cdot / \bandwidth)] ( \bomega') \right\rvert^2}{ \fourier[\Psi_{\bandwidth}](\bomega')}  \diff \bomega'  , & \textrm{(same argument on } g = \Psi_1 \textrm{)}, \\
		&= \Vert \textrm{iso}_\bandwidth(g) \Vert_{\mathcal{H}_\bandwidth}\,,
	\end{align*}
which gives the result.
\end{proof}

\begin{lemma}[$\bar{\varepsilon}_0$ and $\bar{\varepsilon}_2$]
	Let $\bandwidth > 0$ and $r > 0$, then for
    \begin{align*}
        j \in \{0, 2 \}, \quad \bar{\varepsilon}_j\left(r, \bandwidth \right) = \bar{\varepsilon}_j(r, 1)
    \end{align*}
\end{lemma}

\begin{proof} Let $r > 0$. First, notice that due to the scaled euclidean metric induced by $K_\bandwidth$,
\[
\mathfrak{d}_{\metric, \bandwidth}(\bx, \by) \geq r  \iff \mathfrak{d}_{\metric, 1}\left(\frac{\bx}{\bandwidth}, \frac{\by}{\bandwidth}\right) \geq r 
\]
For $j=0$, it suffices to use \Cref{lemma:isometry} \texttt{(i)} together with the above equivalence to show that
\begin{align*}
            \bar{\varepsilon}_0\left(r, \bandwidth \right) = & \sup \left\{ \varepsilon \; : \;  K_\bandwidth(\bx, \by) \leq 1 - \varepsilon, \; \forall \bx, \by \textnormal{ s.t. } \mathfrak{d}_{\metric, \bandwidth}(\bx, \by) \geq r \right\} , \\ 
        = &  \sup \left\{ \varepsilon \; : \;K_1\left(\frac{\bx}{\bandwidth}, \frac{\by}{\bandwidth}\right)\leq 1 - \varepsilon, \; \forall \bx, \by \textnormal{ s.t. } \mathfrak{d}_{\metric, 1}\left(\frac{\bx}{\bandwidth}, \frac{\by}{\bandwidth}\right) \geq r \right\}, \\
        = &  \sup \left\{ \varepsilon \; : \;K_1(\bu, \bv) \leq 1 - \varepsilon, \; \forall \bu, \bv \textnormal{ s.t. } \mathfrak{d}_{\metric, 1}\left(\bu, \bv \right) \geq r \right\}, \\
        = &  \bar{\varepsilon}_0\left(r, 1 \right) .
    \end{align*}

For $j=2$, we begin by noting that $\psi_\bomega(\bx) = \sqrt{f^{(4)}_{\bandwidth}}(\bomega) e^{- \imath \bomega^\top \bx}$ so that, in the same lines as \Cref{lemma:isometry} \texttt{(iii)}:
\begin{align*}
    K^{(02)}_\bandwidth (\bx, \by) [\bz, \bz] &= \Expectation_\bomega \left[ \overline{\psi_\bomega(\frac{\bx}{\bandwidth})} \bz^\top \nabla^2 \psi_\bomega(\frac{\by}{\bandwidth}) \bz \right], \\
    &= \Expectation_\bomega \left[ \overline{\psi_\bomega(\frac{\bx}{\bandwidth})} \psi_\bomega(\frac{\by}{\bandwidth})  \bz^\top (-\imath)^2 \bomega \bomega^\top \bz  \right] ,\\
    &= - \bz^\top \Expectation_\bomega \left[  \bomega \bomega^\top  \overline{\psi_\bomega(\frac{\bx}{\bandwidth})} \psi_\bomega(\frac{\by}{\bandwidth}) \right] \bz,\\
    &= - \bz^\top \nabla^2 K_\bandwidth(\bx - \by) \bz, \\
    &=  \bz^\top \frac{1}{\bandwidth^2} (-\bandwidth^2)\nabla^2 \Psi_\bandwidth(\bx - \by) \bz, \\
    &=   \frac{1}{\bandwidth^2}  \bz^\top \left[- \nabla^2 \Psi_1\left(\frac{\bx - \by}{\bandwidth}\right)\right] \bz, \\
    &=  \frac{1}{\bandwidth^2} \Psi^{(02)}_1 \left(\frac{\bx}{\bandwidth}, \frac{\by}{\bandwidth}\right) [\bz, \bz].
\end{align*}
Moreover, for all $\bz$ and $\bx$, $\Vert \bz \Vert_{\bx,\bandwidth}^2 = \frac{1}{\bandwidth^2}\Vert \bz \Vert_{\bx/\bandwidth,1}^2$. Thus, we again have equality of the $\varepsilon$-sets over which the supremum is taken.
\begin{align*}
    \bar{\varepsilon}_2\left(r, \bandwidth \right) & =  \sup \left\{\varepsilon \, ; \, K^{(02)}_\bandwidth(\bx, \by)[\bz, \bz] \geq \varepsilon \Vert \bz \Vert_{\bx, \bandwidth}^2 , \; \forall z \in \R^d, \; \forall \bx, \by \textnormal{ s.t. } \mathfrak{d}_{\metric, \bandwidth}(\bx, \by) < r \right\}, \\
    &= \sup \left\{\varepsilon \, ; \, \frac{1}{\bandwidth^2} K^{(02)}_1\left(\frac{\bx}{\bandwidth}, \frac{\by}{\bandwidth}\right)[\bz, \bz] \geq \varepsilon \frac{1}{\bandwidth^2}  \Vert \bz \Vert_{\bx/\bandwidth, 1}^2 , \; \forall z \in \R^d, \; \forall \bx, \by \textnormal{ s.t. } \mathfrak{d}_{\metric, 1}\left(\frac{\bx}{\bandwidth}, \frac{\by}{\bandwidth}\right) < r  \right\}, \\
    &= \sup \left\{\varepsilon \, ; \, K^{(02)}_1\left(\frac{\bx}{\bandwidth}, \frac{\by}{\bandwidth}\right)[\bz, \bz] \geq \varepsilon \Vert \bz \Vert_{\bx/\bandwidth, 1}^2 , \; \forall z \in \R^d, \; \forall \bx, \by \textnormal{ s.t. } \mathfrak{d}_{\metric, 1}\left(\frac{\bx}{\bandwidth}, \frac{\by}{\bandwidth}\right) < r  \right\}, \\
    &= \sup \left\{\varepsilon \, ; \, K^{(02)}_1\left(\bu, \bv\right)[\bz, \bz] \geq \varepsilon \Vert \bz \Vert_{\bu, 1}^2 , \; \forall z \in \R^d, \; \forall \bu, \bv \textnormal{ s.t. }\mathfrak{d}_{\metric, 1}\left(\bu, \bv \right) < r \right\}, \\
    &= \bar{\varepsilon}_2\left(r, 1 \right) .
\end{align*}
\end{proof}

\textbf{ $\bullet$ Separation:}
Then, for a number of components $s \in \mathbb{N}$ and $h > 0$, the kernel width of $K_\bandwidth$ is defined by \textcite{poon2023geometry} as
\begin{equation}
    \Delta(h, s) \defeq \inf \left\{ \Delta\, ;\, \sum_{k=2}^s \left \Vert K_\bandwidth^{(ij)}(\trueparam_1, \trueparam_k) \right \Vert_{\trueparam_1, \trueparam_k} \leq h, \quad (i,j) \in \{0,1\} \times \{0, 2\}, \quad \{\trueparam_k \}_{k=1}^s \in \mathcal{S}_{\Delta, \bandwidth} \right\},
\end{equation}
where $\mathcal{S}_{\Delta, \bandwidth} \defeq \{ (\bx_k)_{k=1}^s \in \Param \, ; \, \mathfrak{d}_\bandwidth(\bx_k, \bx_l) \geq \Delta, \forall k \neq l \}$ is the set of $\bs$-tuples of $\Delta$-separated points. This notion is somewhat related to the coherence of the kernel $K$ meaning that the interference is small enough so that functions $K_\bandwidth(\trueparam_k, \cdot)$ can interpolate between parameters. 

\begin{remark}
   The particular form of the metric $\mathfrak{d}_\bandwidth$ gives $S_{\Delta, \bandwidth} = \bandwidth \cdot S_{\Delta, 1}$. Together with \Cref{lemma:norm_tau_equal_norm_one}, this implies that $\Delta(h,s)$ does not depend on $\bandwidth$. 
\end{remark}

\subsubsection{First part of $\LPC$}

This lead to tracking uniform bound on the Kernel and its derivative 
\begin{align}
    B_{ij} \defeq \sup_{\bx, \by} \Vert \Psi_\bandwidth^{(ij)}(\bx - \by) \Vert_{\bx, \by}  & & B_i \defeq 1+ B_{0i} + B_{1i},
\end{align}
which regimes w.r.t. to the dimension $d$ are given by the following lemma. 
\begin{proposition}%
One has
    \label{prop:Bij}
    \[
    \forall i,j, \quad B_{ij} \leq (12d)^{\frac{i+j}{2}},
    \]
and hence $1/\sqrt{B_{02}}\geq 1/\sqrt{12d}$ and $B_2 \leq 1+ \sqrt{12 d} + 12d$.
\end{proposition}

\begin{proof}
    Again, the rate does not depend on $\bandwidth$ as \Cref{lemma:norm_tau_equal_norm_one} gives 
    \[\sup_{\bx, \by} \Vert K_\bandwidth^{(ij)}(\bx, \by) \Vert_{\bx, \by} = \sup_{\bu, \bv} \Vert K_1^{(ij)}(\bu, \bv) \Vert_{\bu, \bv}\,.
    \]
    For the rest of the proof, $\bandwidth$ is fixed equal to $1$.

    Trivially, for $i=j=0$, a triangle inequality yield a uniform control on $(\bu, \bv)$ as:
     \begin{equation}
         \sup_{\bu, \bv} \vert K_1(\bu, \bv) \vert  \leq   \sup_{\bu, \bv}  \int \left \vert e^{\imath \bomega^\top (\bu - \bv)} \right \vert  f^{(4)}_{\bandwidth}(\bomega) \diff \bomega = \int f^{(4)}_{\bandwidth}(\bomega) \diff \bomega = 1 \,.%
     \end{equation}
     For higher order derivative, the proof relies on~\Cref{app:translation_invariant_derivative} and the fact that the covariant derivatives of the kernel involves differential $\nabla^{r} \psi_\bomega(\bu) : (\mathbb{C}^d)^{r} \to \mathbb{C}$ for which we can control the supremum uniformly over all $ \bu \in \mathcal{X}$:
     \begin{equation}
        \sup_{\Vert \bq_l \Vert_2^2 \leq 12, \, l=1, \ldots, r}  \left \vert \nabla^{r}  \psi_\bomega(\bu)[\bq_1, \ldots, \bq_r] \right \vert =  \sqrt{\frac{f^{(4)}_{\bandwidth}}{\SketchDist}(\bomega)} \times \sup_{\Vert \bq_l \Vert_2^2 \leq 12, \, l=1, \ldots, r}  \prod_{l=1}^r \left \vert \left\langle\bomega, \bq_l \right \rangle \right \vert.
     \end{equation}
     Hence, the covariant derivative of the kernel can be uniformly upper bounded for all $\bu, \bv \in \mathcal{X}$ by
    \begin{align*}
        \left \vert K_1^{(ij)}(\bu, \bv)[\bq_1, \ldots, \bq_{i+j}] \right \vert & \leq  \int_{[-1, 1]^d} \left(\prod_{l=1}^{i+j} \sup_{\Vert \bq_l \Vert_2^2 \leq 12} \vert <\bomega, \bq_l > \vert \right) f^{(4)}_{\bandwidth}(\bomega) \vert e^{-\imath \bomega^T (\bu - \bv)} \vert \diff \bomega, \\
        & \leq \sup_{\bomega \in [-1, 1]^d} \prod_{l=1}^{i+j} \sup_{\Vert \bq_l \Vert^2_2 \leq 12}  \vert <\bomega, \bq_l > \vert , \tag{$f^{(4)}_{\bandwidth}$ p.d.f on the unit hypercube} \\
        &  \leq  \sup_{\bomega \in [-1, 1]^d} (2\sqrt{3}  \Vert \bomega \Vert_2)^{i+j}, \tag{Cauchy-Schwarz}\\
        & \leq \sqrt{12d}^{i+j}\,.%
	\end{align*}
 Taking the sup of the latter on $\bu, \bv$ terminates the proof.
\end{proof}

\subsubsection{Second part of $\LPC$}

Then, we focus on the triplets $r_0$, $\bar\varepsilon_0$ and $\bar\varepsilon_2$ defined in \Cref{eq:positive_curvature_constants} which, in turn, define the minimal separation~$\Delta_0$ required between the parameters. Specifically, we need to exhibit some radius $r_0\in(0,1/\sqrt{B_{02}})$ for the near regions such that the kernel $K_\bandwidth$ has positive curvature constants $\bar\varepsilon_0$ and $\bar\varepsilon_2$.
\begin{subequations}
    \label[pluralequation]{sub_eqs:second_part_LPCA_sinc}
\begin{align}
\bar\varepsilon_0&\defeq\frac12\sup_{\varepsilon\geq 0}
    \Big\{
        \varepsilon\,:\,
        K_1(\bs, \bt)\leq 1-\varepsilon\,,\ \forall \bs,\bt\in\Param \textnormal{ s.t. }\distance(\bs, \bt)\geq r_0
    \Big\}\,,\\
\bar\varepsilon_2&\defeq\frac14\sup_{\varepsilon\geq 0}
\Big\{
    \varepsilon\,:\, - K_1^{(0,2)}(\bs, \bt)[\bv,\bv]\geq \varepsilon \|\bv\|_\bt^2\,,\ \forall \bv\in\mathbb T_\bt\,,\ \forall \bs,\bt\in\Param \textnormal{ s.t. }\distance(\bs, \bt)< r_0 \Big\}\,,\\
\Delta_0 &\defeq \inf \Big\{\Delta\,:\,32\sum_{l=2}^{\K}\|K_1^{(i,j)}(\bx_1,\bx_l)\|_{\bx_1,\bx_l}\leq \min\left(\frac{\bar\varepsilon_0}{B_0},\frac{2\bar\varepsilon_2}{B_2}\right)\,,\ (i,j)\in\{0,1\}\times\{0,2\} \,,\ \{\bx_l\}_{l=1}^{\K}\in\mathcal S_\Delta \Big\}.
\end{align}
\end{subequations}
This is the purpose of the following Lemma, which also exhibits the scaling of the last quantities with respect to the dimension $d$. %

\medskip

\begin{proposition}%
 Let $r_0$ be such that $r_0=1/(4\,{d})$ then $\bar\varepsilon_0 \geq 1/(32\,d^3)$ and $\bar\varepsilon_2 \geq 23/128$ in \Cref{sub_eqs:second_part_LPCA_sinc}.
\end{proposition}

\begin{proof}
\begin{subequations}
    One can check the following calculus:
    \begin{align}
    \label{eq:control_sinc_inuvariate}
        \forall z\in\R\,,\quad |\sinc z|\leq (1-\frac{z^2}{12})\indicator_{|z|\leq 2}+\frac12\indicator_{|z|> 2}\,.
    \end{align}
    Now, let $c>0$ and $\bt\in\R^d$ such that $(1/\sqrt{12})\|\bt\|_2\geq c/d$, then there exists a coordinate $z$ of the vector $\bt$ such that $|z|\geq c\sqrt{12}d^{-3/2}$ and we deduce by \eqref{eq:control_sinc_inuvariate} that $\bar\varepsilon_0\geq c^2d^{-3}/2$. The result is given for $c=1/4$.

    We turn to $\bar\varepsilon_2$. Let $\bs,\bt,\bv\in\R^d$ be fixed. Writing the exact Taylor expansion of $\bu \mapsto- K_1^{(0,2)}(\bs, \bu)[\bv,\bv]$ at $\bs$ and using \eqref{eq:Hessian_TVkernel}, there exists $\bxi\in[\bs,\bt]$ (which is the geodesic between these points) such that
    \begin{align*}
        - K_1^{(0,2)}(\bs, \bt)[\bv,\bv] 
        &= - K_1^{(0,2)}(\bs,\bs)[\bv,\bv]  - K_1^{(0,3)}(\bs,\bs)[\bv,\bv,\bs - \bt] - \frac12K_1^{(0,4)}(\bs,\bxi)[\bv,\bv,\bs - \bt,\bs - \bt] 
    \\
    &= \|v\|^2_{\bt}  +0 + \frac{1}{2} \bv^\top\big(A(\bs,\bxi)[\bs - \bt,\bs - \bt]\big) \bv 
    \\
    &\geq (1-6\rho(A(\bs,\bxi)[\bs - \bt,\bs - \bt]))\,\|\bv\|^2_{\bt}
    \end{align*}
    where the quadratic form $K_1^{(0, 4)}(\bs, \bxi)[\bs - \bt, \bs - \bt]$ is identified to the $d \times d$ matrix in the canonical basis:
    \begin{equation}
        \label{def:A_matrix_supplement_SMIX}
        \forall i_1,i_2\in[d]\,,\quad (A(\bs,\bxi)[\bs - \bt, \bs - \bt])_{i_1,i_2}
        =\sum_{i_3=1}^d\sum_{i_4=1}^d
        (s -t)_{i_3}(s - t)_{i_4}\partial_{i_1,i_2,i_3,i_4}{\Psi_1(\bs-\bxi)}\,,
    \end{equation}
    and $\rho(A)$ denotes its spectral radius.

    Then, notice that $f=\sinc$ has all derivatives upper bounded by $1/2$ in absolute value so that
    \begin{align*}
    \frac{d^{4}}{d x^{4}} f^4{\left(x \right)}&=
        \displaystyle 4 f^{3}{\left(x \right)} \frac{d^{4}}{d x^{4}} f{\left(x \right)} + 48 f^{2}{\left(x \right)} \frac{d}{d x} f{\left(x \right)} \frac{d^{3}}{d x^{3}} f{\left(x \right)} \\
        &+
        36 f^{2}{\left(x \right)} \left(\frac{d^{2}}{d x^{2}} f{\left(x \right)}\right)^{2} + 144 f{\left(x \right)} \left(\frac{d}{d x} f{\left(x \right)}\right)^{2} \frac{d^{2}}{d x^{2}} f{\left(x \right)} + 24 \left(\frac{d}{d x} f{\left(x \right)}\right)^{4}\\
        &\leq \frac{256}{16}=16.
    \end{align*}
    Remembering that $\Psi_1=\sinc^4(\cdot/4)$, we deduce from the previous calculus that 
    \begin{align}
        \forall i_1,i_2,i_3,i_4\in[d]\,,\quad |\partial_{i_1,i_2,i_3,i_4}{\Psi_1(\bs-\bxi)}|
        &\leq 16\times \frac{1}{4^4}=\frac{1}{16}\,,\\
        \forall i_1,i_2\in[d]\,,\quad |(A(\bs,\bxi)[\bs - \bt, \bs - \bt])_{i_1,i_2}|
        &\leq \frac{d\|\bs - \bt\|_2^2}{16}\,.
    \end{align}
    Finally, using that $\rho(A)\leq{d}\sup_{i_1,i_2}|A_{i_1,i_2}|$, one has 
    \[
    6\rho(A(\bxi,\bt))\leq \frac{6d^2\|\bs - \bt\|_2^2}{16} = \frac{9d^2}{2} \,{\distance(\bs, \bt)^2}\,,
    \]
    Hence for $r_0=1/(4{d})$ one has $- K_1^{(0,2)}(\bs, \bt)[\bv,\bv] \geq (1-9/32)\|\bv\|^2_{\bt}$.
\end{subequations}
\end{proof}

\begin{remark}
    Note that $r_0< B_{02}^{-1/2}$ since  $r_0 = 1/(4\,{d})$ and $B_{02}^{-1/2}\geq 1/{\sqrt{12d}}$.
\end{remark}

\subsubsection{Last part of $\LPC$}
Note that $B_0\leq 1 +\sqrt{12d}+1 = 2+ \sqrt{12d}$ and $\bar\varepsilon_0/B_0 \geq 1/(32\,d^3(2+ \sqrt{12d}))\geq 1/(175\,d^{7/2})$. Also, it holds 
$B_2\leq 12d +({12d})^{3/2}+1$ and $2\bar\varepsilon_2/B_2 \geq 23/(64+768d +64\,({12d})^{3/2})\geq 23/(1088\,d^{3/2})$. Hence
\[
\frac1{32}\min\left(\frac{\bar\varepsilon_0}{B_0},\frac{2\bar\varepsilon_2}{B_2}\right)\geq \frac{1}{5600\,d^{7/2}}\,.
\]
Using \Cref{lem:bound_sinc4_derivatives} below, we deduce that
\[
\sum_{k=2}^s \Vert K_1^{(ij)} (\bu_k, \bu_1) \Vert_{\bu_{k}, \bu_1}\leq \Big(\frac{4}{3}\Big)^2\sqrt{48d}^{i+j} \frac{sd^2}{\Delta^4}\leq \frac{1}{5600\,d^{7/2}}\leq \frac1{32}\min\left(\frac{\bar\varepsilon_0}{B_0},\frac{2\bar\varepsilon_2}{B_2}\right) 
\]
as soon as $\Delta\geq\Delta_0:=42.66\, s^{1/4}d^{7/4}$.

\begin{lemma}
\label{lem:bound_sinc4_derivatives}
It holds that
   \[
\Vert K_1^{(ij)} (\bs, \bt) \Vert_{\bs, \bt}\leq \Big(\frac{4}{3}\Big)^2\sqrt{48d}^{i+j} \frac{d^2}{\distance(\bs,\bt)^4}\,.
\]
\end{lemma}

\begin{proof}
    We bound each entry of $\Psi_1^{(i+j)} (\bs- \bt)$ using that $|\sinc^{(i)}(x/4)|\leq (1+i)\,4/|x|$ for all $x\neq 0$. One has
    \[
    \Psi_1(\bs- \bt)=\prod_{\ell=1}^d\sinc^4(\frac{s_\ell- t_\ell}{4})\leq \frac{4^4}{(s_{\ell_0}- t_{\ell_0})^4}\prod_{\ell\neq\ell_0}^d\sinc^4(\frac{s_\ell- t_\ell}{4})\leq \frac{4^4d^2}{\|\bs- \bt\|_2^4}
    =\Big(\frac{4}{3}\Big)^2\frac{d^2}{\distance(\bs,\bt)^4}
    \,,
    \]
    because one coordinate $\ell_0$ exists such that $(s_{\ell_0}- t_{\ell_0})^2\geq \|\bs- \bt\|_2^2/d$. Let $h:=\bs- \bt$ and $\partial_u$ the derivative with respect to coordinate $u$, one has
\[
\partial_u \Psi_1(h)=\sinc^{(1)}\big(\frac{h_u}{4}\big)\sinc^{3}\big(\frac{h_u}{4}\big)\prod_{\ell\neq u}\sinc^4\big(\frac{h_\ell}4\big)
\leq \frac{2\,4^4}{(h_{\ell_0})^4}
\leq 2\,\Big(\frac{4}{3}\Big)^2\frac{d^2}{\distance(\bs,\bt)^4}
    \,.
\]
because one coordinate $\ell_0$ exists such that $h_{\ell_0}^2\geq \|h\|_2^2/d$. For the second order derivatives, if $u=v$ one has
\begin{align*}
   \partial^2_{u,v}\Psi_1(h)
   &=\big(\frac14\big)\sinc^{(2)}\big(\frac{h_u}{4}\big)\sinc^{3}\big(\frac{h_u}{4}\big)\prod_{\ell\neq u}\sinc^4\big(\frac{h_\ell}4\big) 
   + \big(\frac34\big)\big(\sinc^{(1)}\big(\frac{h_u}{4}\big)\big)^2\sinc^{2}\big(\frac{h_u}{4}\big)\prod_{\ell\neq u}\sinc^4\big(\frac{h_\ell}4\big) \\
   &\leq 4 \,\Big(\frac{4}{3}\Big)^2\frac{d^2}{\distance(\bs,\bt)^4}\,,
\end{align*}
and otherwise if $u \neq v$
\[
\partial^2_{u,v}\Psi_1(h)=\sinc^{(1)}\big(\frac{h_u}{4}\big)\sinc^{3}\big(\frac{h_u}{4}\big)\sinc^{(1)}\big(\frac{h_v}{4}\big)\sinc^{3}\big(\frac{h_v}{4}\big)\prod_{\ell\neq u,v}\sinc^4\big(\frac{h_\ell}4\big)
\leq 4 \,\Big(\frac{4}{3}\Big)^2\frac{d^2}{\distance(\bs,\bt)^4}\,.
\]
Pursuing the reasoning, one gets
\[
|\partial^{i+j}\Psi_1(h)|\leq 2^{i+j}\,\Big(\frac{4}{3}\Big)^2\frac{d^2}{\distance(\bs,\bt)^4}\,.
\]
Now, observe that
\[
\Vert K_1^{(ij)} (\bs, \bt) \Vert_{\bs, \bt}=\sqrt{12}^{i+j} \Vert K_1^{(ij)} (\bs, \bt) \Vert_2\leq ({12d})^{\frac{i+j}{2}}\max |\partial^{i+j}\Psi_1|
\leq (48{d})^{\frac{i+j}{2}}\,\Big(\frac{4}{3}\Big)^2\frac{d^2}{\distance(\bs,\bt)^4}\,,
\]
which gives the result.
\end{proof}

\subsection{Proof of Theorem~\ref{thm:sinc4_ass2}}
\label{proof:sinc4_ass2}
We refer the reader to the definitions of $N$, $C_1$, $C_2$, $C_3$ and $L_{ij}$ in \Cref{hyp:sketch_function_bounded}. Since $\bar\varepsilon_0$ and $\bar\varepsilon_2$ are already treated in~\Cref{thm:sinc4_lca}, we only need to find the scaling of $L_1$, $L_2$, and $L_3$. Concretely, for $r \in \{0, \ldots, 3\}$, we need to study the tail of the distribution of the random variable
\begin{equation}
    Y_r(\bomega) = \sup_{\bx \in \mathcal{X}} \Vert D_r[\psi_\bomega](\bx) \Vert_{\bx},
\end{equation}
where $\psi_\bomega(\bx) = e^{\imath \bomega^\top \bx} \sqrt{\frac{f^{(4)}_\bandwidth}{\SketchDist}(\bomega)} $ and the convention $\Vert D_0[\psi_\bomega(\bx) \Vert_{\bx} = \vert \psi_\bomega(\bx) \vert$. 
Specifically, we want to exhibits $L_r \in \mathbb{R}$ such that
\[
 \sum_{r=0}^3 1 - \mathbb{P}_{\bomega \sim \SketchDist}(Y_r(\bomega) \leq L_r) \leq \frac{\min(\bar\varepsilon_0, \bar\varepsilon_2, \alpha)}{m}.
\]
In the sinc-4 case, this is given by the following Lemma, independently of $\alpha$, $r_0, \bar\varepsilon_0, \bar\varepsilon_2$ and $m$
\begin{lemma}
    For any $\bandwidth > 0$, and any $r \in \{0, \ldots, 3\}$, the covariant derivatives are almost surely bounded and we have 
    \begin{equation}
        \label{eq:sinc4_Lr_scaling}
        L_r = (\sqrt{12d})^r \sup_{\bomega} \sqrt{\frac{f^{(4)}_{\bandwidth}}{\SketchDist}}(\bomega).
    \end{equation}
\end{lemma}
\begin{proof}
Since the support of $\SketchDist$ is bounded and $\psi_\bomega(\cdot)$ is a weighted Fourier feature, we don't need refined probabilistic control on the tail bound of $L_r(\bomega)$. Indeed, using~\Cref{app:translation_invariant_derivative} we have
\begin{align}
    Y_0(\bomega) &= \sup_\bx \vert \psi_\bomega(\bx) \vert = \sqrt{\frac{f^{(4)}_{\bandwidth}}{\SketchDist}}(\bomega), \\
    Y_1(\bomega) &= \sup_\bx \Vert \metric_{\bx, \tau}^{-1/2} \nabla \psi_\bomega(\bx) \Vert_2 =  \sqrt{\frac{f^{(4)}_{\bandwidth}}{\SketchDist}}(\bomega) \sqrt{12}\bandwidth\Vert \bomega \Vert_2, \\
    Y_2(\bomega) &= \sup_\bx \Vert \metric_{\bx, \tau}^{-1/2} \nabla^2 \psi_\bomega(\bx) \metric_{\bx, \tau}^{-1/2} \Vert_{op} =  \sqrt{\frac{f^{(4)}_{\bandwidth}}{\SketchDist}}(\bomega) (\sqrt{12}\bandwidth)^2 \Vert \bomega \Vert_2^2, \\ 
    Y_3(\bomega) &= \sup_\bx \Vert \nabla^3 \psi_\bomega(\bx)[\metric_{\bx, \tau}^{-1/2} \cdot, \metric_{\bx, \tau}^{-1/2} \cdot, \metric_{\bx, \tau}^{-1/2} \cdot]\Vert_{op} =  \sqrt{\frac{f^{(4)}_{\bandwidth}}{\SketchDist}}(\bomega)(\sqrt{12}\bandwidth)^3 \Vert \bomega \Vert_2^3.  
\end{align}
Thus, since $\SketchDist$ has a bounded support on $\Omega=[-1/\bandwidth, 1/\bandwidth]^d$, and using $\Vert \bomega \Vert_2 \leq \sqrt{d} \Vert \bomega \Vert_{\infty}$, the random variables $Y_r(\bomega)$ are a.s. bounded by
\[
Y_r(\bomega) \leq \sup_{\bomega \in \Omega} \sqrt{\frac{f^{(4)}_{\bandwidth}}{\SketchDist}}(\bomega)(\sqrt{12} \bandwidth)^r \bigg(\frac{\sqrt{d}}{\bandwidth}\bigg)^r = \sup_{\bomega \in \Omega} \sqrt{\frac{f^{(4)}_{\bandwidth}}{\SketchDist}}(\bomega)(\sqrt{12} \sqrt{d} )^r.
\]
Denoting $C_\SketchDist = \sup_{\bomega \in \Omega} \frac{f^{(4)}_{\bandwidth}}{\SketchDist}(\bomega)$, we get that 
\[
\mathbb{P}_{\bomega \sim \SketchDist}\Big( Y_r(\bomega) \leq \sqrt{C_\SketchDist} (\sqrt{12} \sqrt{d} )^r \Big) = 1,
\]
so that $\psi_\bomega$ satisfies \Cref{hyp:sketch_function_bounded} with $L_r \defeq \sqrt{C_\SketchDist} (\sqrt{12} \sqrt{d} )^r$. The proof ends by plugging these values in $N, C_1, C_2$ and $L_{ij}$.
\end{proof}

\subsection{Proof of Theorem~\ref{thm:main_thm}}
\label{proof:main_thm}
By Theorem~\ref{thm:existence_certificates}, we known that there exists pivot non-degenerate certificates. Let $\eta^0$ be an $({\bar\varepsilon_0},\bar\varepsilon_2,r_0)$-pivot non-degenerate dual certificate. Note that
\begin{align}
    \frac{1}{2} \big\Vert \by - \Forward \mu \big\Vert_{\mathcal{F}}^2 + \kappa \Vert \mu \Vert_{\mathrm TV}
    \leq
    \frac{1}{2} \big\Vert \by - \Forward \target \big\Vert_{\mathcal{F}}^2 + \kappa \Vert \target \Vert_{\mathrm TV}
    \leq
    \frac{\gamma^2}{2} + \kappa \Vert \mu \Vert_{\mathrm TV}
\end{align}
Using $\eta^0=F^\star c^0$ and \eqref{eqs:f_star_properties}, we obtain
\begin{align}
    \kappa\mathcal D_{\eta^0}(\mu\,||\,\target) + \kappa\langle\eta^0,\mu-\target\rangle_{\mathcal C(\Param)\times\mathcal M(\Param)}+\frac{1}{2} \big\Vert \by - \Forward \mu \big\Vert_{\mathcal{F}}^2
    &\leq 
    \frac{\gamma^2}{2}\notag\\
    \kappa\mathcal D_{\eta^0}(\mu\,||\,\target) + \frac{1}{2} \big\Vert \kappa c^0 + \Forward \mu - \by\big\Vert_{\mathcal{F}}^2
    &\leq 
    \frac{\gamma^2}{2} + \frac{\kappa^2\|c^0\|_{\mathcal{F}}^2}{2}-\langle \kappa c^0,\Gamma\rangle_{\mathcal{F}}\label{eq:intermediate_main_proof}\\
    \mathcal D_{\eta^0}(\mu\,||\,\target) 
    &\leq \frac{(\gamma+\kappa\|c^0\|_{\mathcal{F}})^2}{2\kappa}\notag\\
    \notag
    &\leq \frac{(\gamma+C_{\switch}\kappa\sqrt{2\K})^2}{2\kappa}
\end{align}
where we recall that
\[
\mathcal D_{\eta^0}(\mu\,||\,\target)\defeq\|\mu\|_{\mathrm{TV}}-\|\target\|_{\mathrm{TV}}-\langle\eta^0,\mu-\target\rangle_{\mathcal C(\Param)\times\mathcal M(\Param)}\,,
\]
see \eqref{eq:def_bregman}, and $\Gamma\defeq\by-\Forward\target$. Choosing $\kappa=c_\kappa\gamma/\sqrt{\K}$ one gets
\begin{equation}
\label{eq:control_divergence_general}
    \mathcal D_{\eta^0}(\mu\,||\,\target)\leq \bar c_\kappa\gamma \sqrt{\K}
    \quad\textnormal{where}\quad\bar c_\kappa\defeq\frac{(1+\sqrt 2\,C_{\switch}c_\kappa)^2}{2c_\kappa}\geq2\sqrt{2}\,C_{\switch}\,.
\end{equation}
On the other hand, using that $\eta^0$ is an $({\bar\varepsilon_0},\bar\varepsilon_2,r_0)$-pivot non-degenerate dual certificate, for $r\in(0,r_0]$,
\begin{align}
    \mathcal D_{\eta^0}(\mu\,||\,\target)
    &=
    \|\mu\|_{\mathrm{TV}}-\langle\eta^0,\mu\rangle_{\mathcal C(\Param)\times\mathcal M(\Param)}\notag\\
    &\geq 
    \|\mu\|_{\mathrm{TV}}-\sum_{l=1}^{\K}\int_{\nearregion_l(r)}|\eta^0|\,\mathrm{d}|\mu|-\int_{\farregion_l(r)}|\eta^0|\,\mathrm{d}|\mu|
    \notag
\end{align}
Now observe that, for $r\in(0,r_0]$,
\begin{itemize}
    \item For all $i\in[\K]$, $\eta^0(t^0_i)=\mathrm{sign}(a^0_i)$,
    \item For all $x\in\farregion(r)$, 
    \begin{equation}
    \label{eq:new_control_certificate}
        |\eta^0(x)|\leq1-\min(\bar\varepsilon_0,\bar\varepsilon_2 r^2)=1-\bar\varepsilon_2 r^2\,,
    \end{equation}
    \item For all $i\in[\K]$, for all $x\in\nearregion_i(r)$, $|\eta^0(x)|\leq1-\bar\varepsilon_2\distance(x,t_i^0)^2$.
\end{itemize}
Hence,
\begin{align}
    \label{eq:lb_divergence_general}
    \mathcal D_{\eta^0}(\mu\,||\,\target)
    &\geq 
    \bar\varepsilon_2 r^2|\mu|(\farregion(r))+\bar\varepsilon_2\sum_{l=1}^{\K}\int_{\nearregion_l(r)}\distance(x,t_i^0)^2\,\mathrm{d}|\mu|(x)\,.
\end{align}
Using~\eqref{eq:control_divergence_general} and~\eqref{eq:lb_divergence_general}, one gets~\eqref{eq:control_far_thm} and, by way of contradiction,~\eqref{eq:detection_near_thm}.

It remains to prove~\eqref{eq:control_near_thm}. Observe that, using~\eqref{eq:control_divergence_general} and~\eqref{eq:new_control_certificate} and~\eqref{eq:lb_divergence_general},
\begin{align}
    \big|\mu(\nearregion_j(r))-a^0_j\big|&=\Big|\int_{\nearregion_j(r)}\mathrm{d}(\mu-\target)\Big| 
    \notag\\
    &=
        \Big|
            \int_{\Param}\eta^0_j\mathrm{d}(\mu-\target)
            + \int_{\nearregion_j(r)}\!\!\!(1-\eta^0_j)\mathrm{d}(\mu-\target)\notag\\
    &\quad\quad\quad\quad\quad\quad\quad\quad\  
            - \sum_{l\neq j}\int_{\nearregion_l(r)}\!\!\!\eta^0_j\mathrm{d}(\mu-\target)
            - \int_{\farregion(r)}\!\!\!\eta^0_j\mathrm{d}(\mu-\target)
        \Big| 
    \notag\\
    &\leq
        \Big|
        \int_{\Param}\eta^0_j\mathrm{d}(\mu-\target)
        \Big|
        +\bar\varepsilon_2\sum_{l=1}^\K \int_{\nearregion_l(r)}\distance(x,t_j^0)^2\mathrm{d}|\mu|(x)+(1-\bar\varepsilon_2r^2)|\mu|(\farregion(r))
    \notag\\
    & \leq 
        \Big|
        \int_{\Param}\eta^0_j\mathrm{d}(\mu-\target)
        \Big|
        + \frac{c_r}{r^2}
        \Big(
            \bar\varepsilon_2 r^2|\mu|(\farregion(r))+\bar\varepsilon_2\sum_{l=1}^{\K}\int_{\nearregion_l(r)}\distance(x,t_i^0)^2\,\mathrm{d}|\mu|(x)
        \Big)
    \notag\\
    & \leq 
        \Big|
        \int_{\Param}\eta^0_j\mathrm{d}(\mu-\target)
        \Big|
        + c_r\bar c_\kappa\Big(\frac{\gamma}{r^2}\Big) \sqrt{\K}
    \label{eq:final_proof_main}
\end{align}
where 
\[
c_r\defeq \max(r^2,\frac{1-\bar\varepsilon_2r^2}{\bar\varepsilon_2})\leq\frac{\max(1,\bar\varepsilon_0)}{\bar\varepsilon_2}\,.
\]
Now, note that
\begin{align}
    \Big|
        \int_{\Param}\eta^0_j\mathrm{d} (\mu - \targetdensity)
        \Big|&=\Big|
        \langle\eta^0_j,\mu-\target\rangle_{\mathcal C(\Param)\times\mathcal M(\Param)}
        \Big|\notag\\
        &=\Big|
        \langle c^0_j,F\mu-F\target\rangle_{\mathcal{F}}
        \Big|\notag\\
        &\leq \|c^0_j\|_{\mathcal{F}}(\gamma + \|\by -F\mu\|_{\mathcal{F}})\label{eq:intermediate_3_main_proof}
\end{align}
and, using \eqref{eq:intermediate_main_proof}, one has
\[
\big\Vert \kappa c^0 + \Forward \mu - \by \big\Vert_{\mathcal{F}}
\leq \gamma+\kappa\|c^0\|_{\mathcal{F}}\,,
\]
which leads to
\begin{equation}
    \label{eq:intermediate_2_main_proof}
    \big\Vert\Forward \mu - \by \big\Vert_{\mathcal{F}}\leq \gamma+2\kappa\|c^0\|_{\mathcal{F}}\,.
\end{equation}
By~\eqref{eq:intermediate_3_main_proof} and \eqref{eq:intermediate_2_main_proof}, one has
\begin{align}
\label{eq:intermediate_4_main_proof}
    \Big|
        \int_{\Param}\eta^0_j\mathrm{d}(\mu-\target)
    \Big|
    \leq 
    2\|c^0_j\|_{\mathcal{F}}(\gamma+\kappa\|c^0\|_{\mathcal{F}})
    \leq 2\sqrt 2\,C_{\switch}(1+\sqrt 2\,C_{\switch}c_\kappa)\gamma
\end{align}
using Theorem~\ref{thm:existence_certificates}. Combining~\eqref{eq:final_proof_main} and~\eqref{eq:intermediate_4_main_proof}, one gets the result.

\subsection{Proof of Theorem~\ref{thm:main_thm_sketch}}
\label{proof:main_thm_sketch}
The choice $C_{\textnormal{sketch}} = 2C\max(C_1, C_2) $ together with \Cref{eq:m0_thm} ensures that $m \geq m_0$ from \Cref{hyp:sketch_function_bounded}. First invoke Theorem~\ref{thm:existence_sketch_certificates} to get that there exists a constant $C'_{\pivot}>0$ which depends only on the kernel $K_\pivot$ such that, with probability at least $1-\alpha$, for all $m\geq m_0$, there exist a $({\bar\varepsilon_0/4},3\bar\varepsilon_2/2,r_0)$-sketch pivot non-degenerate dual certificate such that $\|c^0\|_{\mathbb C^m}\leq {C'_{\pivot}}C_{\switch}\sqrt{\K}$ and $(\bar\varepsilon_0/4,3\bar\varepsilon_2/2,r_0)$-sketch pivot non-degenerate localizing certificates at points~$t_i^0$ such that $\|c^0_i\|_{\mathbb C^m}\leq C'_{\pivot}C_{\switch}$. Then the proof goes the same line as in Appendix~\ref{proof:main_thm} substituting $\bar\varepsilon_0$ by ${\bar\varepsilon_0/4}$ and $\bar\varepsilon_2$ by $3\bar\varepsilon_2/2$. 

\subsection{Proof of \Cref{prop:supermix_thm10}}
\label{proof:supermix_thm10}

The proof is an application of~\Cref{thm:main_thm} using the sinc-4 as pivot kernel. Let $\mu_0 \in \Measures(\Param)$ be a sparse measure of $s_0$ spikes $\{\trueparam_1, \ldots, \trueparam_{s_0}\}$. First, recall~\Cref{prop:control_population_noise_level_whp} gives a control with probability $1 - \alpha$
\begin{equation}
    \gamma_n := C_\alpha \bigg(\frac{1}{\bandwidth}\bigg)^{\frac d2}{\frac{1}{\sqrt n}}\,,
\end{equation}
where $C_\alpha \defeq 2\sqrt{1 + C_1 \log({C_2}/{\alpha})}$ is a constant only depending on $\alpha$, while $C_1$ and $C_2$ are universal constants.

By~\Cref{thm:sinc4_lca}, the sinc-4 pivot kernel $K_\pivot(\bs, \bt) = \Psi_\bandwidth(\bs - \bt)$ verifies the $\LPC$ assumption with radius $r_0=1/(4\,{d})$, $\bar\varepsilon_0 \geq 1/(32\,d^3)$, $\bar\varepsilon_2 \geq 23/128$ and $\Delta_0 = 147.77\,\K^{1/4}d^{7/4}$, which are independent of the bandwidth $\bandwidth$. Moreover, its associated Fisher-Rao distance is a rescaled euclidean $\distance_\bandwidth(\bs, \bt) = (2\sqrt{3} \bandwidth)^{-1} \Vert \bs - \bt \Vert^2$. Thus, as discussed in \Cref{eq:condition_bandwidth}, taking $\bandwidth$ verifying~\Cref{hyp:Hbandwidth} implies that the target belongs to the model set $\Model_{s_0, \Delta_0, \distancegeneric_\bandwidth}$. Then, we may invoke \Cref{thm:main_thm} with regularisation 
\[
\kappa = c_\kappa \gamma/\sqrt{\K} = \frac{1}{\sqrt{2}C_\switch}C_\alpha {\bandwidth}^{-\frac{d}2}{\frac{1}{\sqrt n}} \frac{1}{\sqrt{\K}},
\]
where we used the optimal proportionality constant $c_\kappa = 1 / \sqrt{2}C_\switch$ as discussed in~\Cref{rem:tuning_kappa_main_thm}. Injecting this in~\Cref{eq:control_far_thm,eq:control_near_thm,eq:detection_near_thm}, one gets for any radius $r >0$ such that $r<\min\big(r_0,\sqrt{\bar\varepsilon_0/\bar\varepsilon_2}\big) = \mathcal{O}(d^{-3/2}) =: c'_d$,
\begin{itemize}
    \item Control of the far region: 
    \begin{equation}
        |\mu|(\farregion(r))\leq 2 \sqrt{2} C_\switch \frac{128}{23} C_\alpha {\bandwidth}^{-\frac{d}2} \Big(\frac{1}{r^2 \sqrt{n} }\Big) \sqrt{\K}\,,
    \end{equation}
    \item Control of all the near regions: for all $k\in[\K]$,
    \begin{equation}
        \label{eq:control_far_proof_prop21}
        |\mu(\nearregion_k(r))-a^0_k|\leq 2 \sqrt{2} C_\switch \frac{128}{23} C_\alpha {\bandwidth}^{-\frac{d}2} \Big(\frac{1}{r^2 \sqrt{n} }\Big) \sqrt{\K}+ 4 \sqrt{2} C_\switch  C_\alpha {\bandwidth}^{-\frac{d}2} \frac{1}{\sqrt{n}}\,,
    \end{equation}
    \item Detection level: for all Borelian $A\subset\Param$ such that $|\mu|(A)>2 \sqrt{2} C_\switch \frac{128}{23} C_\alpha {\bandwidth}^{-\frac{d}2} \big(\frac{1}{r^2 \sqrt{n} }\big) \sqrt{\K}$, there exists $\trueparam_k$ such that
    \begin{equation}
        \min_{\bt\in A}\distance(t,\trueparam_k)\leq  r\,,
    \end{equation} 
\end{itemize}
The last part of the proof consists in realizing that we may take the radius $r = r_n$ to depend on $n$. The only constraint is that $r_n^2 \sqrt{n} $ goes to $+\infty$ as $n$ grows so that the controls asymptotically make sense, \textit{i.e.} the BLASSO estimator concentrates around the true support. Thus, we may consider a slowly diverging sequence $\delta_n$ such that~$r_n^2 \sqrt{n} = \delta_n^2$, or equivalently $r_n = n^{-1/4} \delta_n$. For comments over the choice of $\delta_n$, we refer to \Cref{rem:effective_near_region}. Note that, for large $n$, $\delta^2_n = r_n^2 \sqrt{n} < \sqrt{n}$ so that the first term on the right-hand side of~\Cref{eq:control_far_proof_prop21} is the dominant one. Hence, we obtain the final control in~\Cref{eq:control_near_mixture}
\[
|\mu(\nearregion_k(r))-a^0_k|\leq 6 \sqrt{2} C_\switch \frac{128}{23} C_\alpha {\bandwidth}^{-\frac{d}2} \Big(\frac{1}{r^2 \sqrt{n} }\Big) \sqrt{\K}.
\]
Finally, the condition $r_n < c'_d$ can be restated implicitly in $n$ as $n \geq (c'_d)^{-4} \delta_n^4$, and the constant $c_d \defeq (c'_d)^{-4} = \mathcal{O}(d^6)$ depends polynomially on the dimension. This concludes the proof.

\subsection{Proof of \Cref{prop:sketch_Mixture_Blasso}}
\label{proof:sketch_Mixture_Blasso}
The proof is an application of~\Cref{thm:main_thm_sketch} which give the controls for sketched BLASSO problems. We first prove the control on the sketch noise level and then turn on proving the actual proposition.

\subsubsection{Proof of \Cref{lem:noise_level_sketched}}
\label{proof:lemma_noise_level_sketch}

First, we recall and prove~\Cref{lem:noise_level_sketched} giving the control over the sketched noise level with probability $1- \alpha$ w.r.t. the joint draws of the sample $\Obs$ and the sketch. Precisely, we prove that with probability $1- \alpha$ it holds that 
    \[
    \Vert \Gamma_{\textnormal{sketch}} \Vert_{\mathbb{C}^m} \leq C_{\alpha,m } \frac{1}{\sqrt{n}},
    \]
    with 
    \[
    C_{\alpha,m} = 2\sqrt{\bigg[ \Big(\frac1\bandwidth\Big)^{d} + \frac{1}{2 \sqrt{m}}\left\Vert \frac{U_\bandwidth}{\SketchDist} \right\Vert_\infty\! \log \Big( \frac{2}{\alpha} \Big) \bigg]  \bigg[ 1 + C_1 \log \big( \frac{2C_2}{\alpha} \Big) \bigg]},
    \]
    and $C_1, C_2$ universal constants.

    \pagebreak[3]

    We follow the proof of \textcite[Lemma 3]{decastro2019sparse} and adapt it with a union bound argument to take into account the extra-randomness from the sketch. First, introduce the i.i.d.~random variables $(\bY_j)_{j=1}^n$
    \[
    \bY_j \defeq \mathcal{A} \delta_{\bZ_j} - \Expectation_{\bZ} \mathcal{A} \delta_\bZ,
    \] 
    where $\mathcal{A} : \nu \in \Measures(\mathbb{R}^d) \mapsto \frac{1}{\sqrt{m}} \left( W(\bomega_i) \fourier[\nu](\bomega_i) \right)_{i=1}^m \in \mathbb{C}^m$ is the so-called \textit{sketching} operator yielding the random Fourier features of a Radon measure $\nu$. Then, we have
    \begin{equation}
    \label{eq:noise_as_U_process}
    \Big\Vert \Gamma_{\textnormal{sketch}} \Big\Vert_{\mathbb{C}^m}^2 
    = \Big \Vert \frac{1}{n} \sum_{j=1}^n \bY_j \Big \Vert_{\mathbb{C}^m}^2 =  \frac{1}{n^2} \sum_{j=1}^n \left \Vert \bY_j \right \Vert_{\mathbb{C}^m}^2 + \frac{1}{n^2} \sum_{j\neq l}^n \< \bY_j, \bY_l \>_{\mathbb{C}^m} ,
    \end{equation}
    Fix $\alpha \in (0,1)$, and let us assume that there exists a constant $M$ such that $ \Vert \bY_j \Vert^2_{\mathbb{C}^m} \leq M$. Then, using the same arguments as in the proof of \textcite[Lemma 3]{decastro2019sparse}, the first term of \Cref{eq:noise_as_U_process} is bounded by~$M/n$ while the second term is a canonical U-process. Thus, \textcite[Proposition 2.3]{arcones1993limit} ensures that there exists two universal constants $C_1, C_2 >0$ such that 
    \begin{equation}
        \label{eq:Arcones}
        \Big\Vert \Gamma_{\textnormal{sketch}} \Big\Vert_{\mathbb{C}^m}^2 \leq M \Big(1 + C_1\log\bigg(\frac{2C_2}{\alpha}\bigg) \Big) \frac{1}{n}
    \end{equation}
    with probability at least $1 - \frac{\alpha}{2}$.
    
    We then proceed to show that $\Vert \bY_j \Vert^2_{\mathbb{C}^m}$ is bounded with probability $1-\alpha/2$ using Hoeffding's inequality. Indeed, for any $j \in [n]$ we may write
    \[
    \Vert \bY_j \Vert_{\mathbb{C}^m}^2 = \frac{1}{\nsketch} \sum_{i=1}^m \left\vert W(\bomega_i) \left( e^{-\imath \bomega_i^\top \bZ_j } - \fourier[\targetdensity](\bomega_i) \right) \right\vert^2 = \frac{1}{\nsketch} \sum_{i=1}^m X_i^j,
    \]
    with $X_i^j = \left\vert W(\bomega_i) \left( e^{-\imath \bomega_i^\top \bZ_j } - \fourier[\targetdensity](\bomega_i) \right) \right \vert^2$. Moreover, the expectation with respect to the draw of $\bomega_{1:m}$ is
    \[
    \Expectation_{\bomega_{1:m}} \Vert \bY_j \Vert_{\mathbb{C}^m}^2 = \Vert L_\bandwidth \delta_{\bZ_j} - L_\bandwidth \targetdensity \Vert_{\mathcal{F}_\bandwidth}^2,
    \]
    where we recall the operator $L_\bandwidth \nu = \lambda_\bandwidth \star \nu$ is the smoothing kernel operator, which can be seen as the kernel mean embedding of $\nu$ in the RKHS $\mathcal{F}_\bandwidth$. The $(X_i^j)_{i=1}^m$ are independent and bounded random variables with
    \[
    0 \leq X_i^j \leq  W^2(\bomega_i) \left(  \vert e^{-\imath \bomega_i^\top \bZ_j}\vert + \vert \fourier[\targetdensity](\bomega_i)\vert  \right )^2 \leq 4 \Vert W^2 \Vert_\infty\,,
    \]
    using Jensen's inequality to get that 
    \[
    \vert \fourier[\targetdensity](\bomega_i)\vert
    =
    \vert \Expectation_\bZ e^{-\imath \bomega_i^\top \bZ}\vert
    \leq 
    \Expectation_\bZ \vert e^{-\imath \bomega_i^\top \bZ}\vert=1\,.
    \]
    Applying Hoeffding's inequality yields that with probability greater than $1 - \alpha/2$
    \[
    \Vert \bY_j \Vert_{\mathbb{C}^m}^2 \leq \Vert L_\bandwidth \delta_{\bZ_j} - L_\bandwidth \targetdensity \Vert_{\mathcal{F}_\bandwidth}^2 + \frac{2  \Vert W^2 \Vert_\infty}{\sqrt{m}} \log \bigg( \frac{2}{\alpha} \bigg)
    \]
    
    Now observe that 
    \begin{align*}
        \Vert L_\bandwidth \delta_{\bZ_j} \Vert_{\mathcal{F}_\tau}^2
        &=\lambda_\bandwidth(\bZ_j-\bZ_j)=\lambda_\bandwidth(\bm{0})\\
        \Vert  L_\bandwidth \targetdensity \Vert_{\mathcal{F}_\tau}^2
        &\leq \Expectation_{\bZ}{\Vert L_\bandwidth \delta_{\bZ} \Vert_{\mathcal{F}_\tau}^2}=\lambda_\bandwidth(\bm{0})
    \end{align*}
    by the reproducing property of the RKHS $\mathcal{F}_\tau$ and Jensen's inequality. We deduce that 
    \[\Vert L_\bandwidth \delta_{\bZ_j} - L_\bandwidth \targetdensity \Vert_{\mathcal{F}_\bandwidth}^2  \leq 4 \lambda_\bandwidth(\bm{0})\,,
    \]
    so that we can conclude
    \begin{equation}
        \label{eq:bound_Yij_whp}
          \Vert \bY_j \Vert_{\mathbb{C}^m}^2 \leq M_{\alpha, m} \defeq  4 \lambda_\bandwidth(0) + \frac{2  \Vert W^2\Vert_\infty}{\sqrt{m}} \log \bigg( \frac{2}{\alpha} \bigg)
    \end{equation}
    with probability greater than $1 - \alpha/2$. A union bound bound yields that the joint probability of \Cref{eq:Arcones,eq:bound_Yij_whp} is greater than $1-\alpha$. The proof ends by plugging $\lambda_\bandwidth(\bm{0}) = (\frac{1}{\bandwidth})^d$ and $W^2 = \frac{U_\bandwidth}{\SketchDist}$.

\subsubsection{Proof of the main proposition}
The proof of~\Cref{prop:sketch_Mixture_Blasso} follows the same line as the proof of~\Cref{prop:supermix_thm10} given above in~\Cref{proof:supermix_thm10}. First, remark that \Cref{thm:sinc4_ass2} ensures the sinc-4 kernel verifies~\Cref{hyp:sketch_function_bounded} and is thus amenable to sketching with the following constants:
  \begin{align*}
    m_0 & = C\K
    \Bigg(
    C_1\log(\K)\log\Big(\frac{\K}{\alpha}\Big)
    +C_2\log\bigg(\frac{(\K N)^d}{\alpha}\bigg)
    \Bigg) \\
    N
    &\defeq |\Param| 32\sqrt{12} \sqrt{C_\SketchDist} d^{7/2} + \frac{128}{23}( 12 \sqrt{12} C_\SketchDist d^{1/2} + \sqrt{C_\SketchDist} 12 d) 
    \\ \notag
    C_1
    &=\mathcal{O}\left(d^{15/2} C_\SketchDist\right)
    \\ \notag
    C_2
    &= \mathcal{O}\left(d^{6} C_\SketchDist\right)
    \\ \notag
    C_\SketchDist 
    &\defeq \sup_{\bomega \in [-\frac{1}{\bandwidth}, \frac{1}{\bandwidth}]^d} \frac{f^{(4)}_{\bandwidth}}{\SketchDist} (\bomega) 
    \end{align*}
    Where the diameter $|\Param| = \sup_{\bs,\bt\in\Param}\distance_\bandwidth(\bs, \bt)$, and $C$ depends on $\Psi_\bandwidth$ and polynomially on $d$.

With our hypothesis on the bandwidth $\bandwidth$ \Cref{hyp:Hbandwidth}, the target $\target$ belongs to the model set $\Model_{\K, \Delta_0, \distancegeneric_\bandwidth}$ where $\distancegeneric_\bandwidth$ is the Fisher-Rao metric of the sinc-4 kernel $\Psi_\bandwidth$ (see \Cref{eq:condition_bandwidth}). We can then invoke~\Cref{thm:main_thm_sketch}, which ensures a positive constant $C'_{\pivot}$ exists and yields the sketched BLASSO controls with high-probability. Precisely, take the sketch size $m \geq m_0$, and use the control on the sketched noise level $\gamma_\textnormal{sketch}$ to tune the regularisation
\[
\kappa = \frac{c_\kappa' \gamma_\textnormal{sketch}}{\sqrt{\K}} = \frac{1}{C'_\pivot C_\switch}{\frac{C_{\alpha,m}}{\sqrt n}} \frac{1}{\sqrt{\K}},
\]
where we used the optimal $c'_\kappa = \frac{1}{C'_\pivot C_\switch}$. Injecting this in~\Cref{eq:control_far_thm_sketch,eq:control_near_thm_sketch,eq:detection_near_thm_sketch}, one gets for any radius $r >0$ such that $r<\min\big(r_0,\sqrt{\bar\varepsilon_0/6\bar\varepsilon_2}\big) =: c'_d$,
\begin{itemize}
    \item Control of the far region: 
    \begin{equation}
        |\mu|(\farregion(r))\leq 2 C'_\pivot C_\switch \frac{256}{69} C_{\alpha,m} \Big(\frac{1}{r^2 \sqrt{n} }\Big) \sqrt{\K}\,,
    \end{equation}
    \item Control of all the near regions: for all $k\in[\K]$,
    \begin{equation}
        \label{eq:control_far_proof_prop22}
        |\mu(\nearregion_k(r))-a^0_k|\leq  2 C'_\pivot C_\switch   \frac{256}{69} C_{\alpha,m} \Big(\frac{1}{r^2 \sqrt{n} }\Big) \sqrt{\K} + 4 C'_\pivot C_\switch  \frac{C_{\alpha,m}}{\sqrt{n}}\,,
    \end{equation}
    \item Detection level: for all Borelian $A\subset\Param$ such that $|\mu|(A)>2  C'_\pivot C_\switch \frac{256}{69} C_{\alpha,m} \Big(\frac{1}{r^2 \sqrt{n} }\Big) \sqrt{\K}$, there exists $\trueparam_k$ such that
    \begin{equation}
        \min_{\bt\in A}\distance(t,\trueparam_k)\leq  r\,,
    \end{equation} 
\end{itemize}
Using a similar reasoning as in~\Cref{proof:supermix_thm10}, we finish the proof by taking the radius $r = r_n$ to depend on $n$. Let $\delta_n$ be a slowly diverging sequence such that~$r_n^2 \sqrt{n} = \delta_n^2$, or equivalently $r_n = n^{-1/4} \delta_n$. Again, the first term on the right-hand side of~\Cref{eq:control_far_proof_prop22} is the dominant one. Hence, we obtain the final control in~\Cref{eq:control_near_mixture}
\[
|\mu(\nearregion_k(r))-a^0_k|\leq 6 C'_\pivot C_\switch   \frac{256}{69} C_{\alpha,m} \Big(\frac{1}{r^2 \sqrt{n} }\Big) \sqrt{\K}.
\]
Again, the condition on the radius $r_n < c'_d$ is met for any such that $n \geq (c'_d)^{-4} \delta_n^{4}$, with $c_d \defeq (c'_d)^{-4} = \mathcal{O}(d^6)$ scaling polynomially in the dimension $d$. This concludes the proof.

\section{Technical results and remarks}

\subsection{Covariant derivatives: the translation invariant kernel case}
\label{sec:covariant_derivatives_TI}
    \label{app:translation_invariant_derivative}
    The covariant derivative of the kernel $K(\bs, \bt) = \mathbb{E}_{\bomega \sim \SketchDist}[\psi_\bomega(\bs) \overline{\psi_\bomega(\bt)}]$ is properly introduced by \textcite[Section 4.1]{poon2023geometry} as a ``bi''-multilinear map $K^{(ij)}(\bs, \bt) : (\mathbb{C}^d)^i \times (\mathbb{C}^d)^j \to \mathbb{C}$, and it involves the covariant derivatives of order $i$ and $j$ of the feature map $\psi_\bomega : \mathbb{R}^d \to \mathbb{C}$. In the case of translation invariant kernels $K(\bs, \bt)=\rho(\bs - \bt)$, the Riemannian covariant derivatives $ K^{(i,j)}(\bs, \bt)[U,V]$ and its Riemannian operator norm $\Vert K^{(i,j)}(\bs, \bt)\Vert_{\bs,\bt}$ are given by the standard Euclidean derivatives $(-1)^j\nabla^{i+j}\rho(h)$ of $\rho$ at point $\bm{h}=\bs - \bt$ and the Euclidean operator norm, namely
    \begin{subequations}
        \begin{align}
            \label{eq:translation_invariant_derivative}
            K^{(i,j)}(\bs, \bt)[U,V]&= (-1)^j\mathrm{Trace}(\nabla^{i+j}\rho(\bm{h}) \times U_1\otimes\cdots\otimes U_i\otimes V_1\otimes\cdots\otimes V_j)\,,\\
            \label{eq:translation_invariant_derivative_norm}
            \Vert K^{(i,j)}(\bs, \bt)\Vert_{\bs,\bt}&=\Big\Vert \nabla^{i+j}\rho(\bm{h})\bigtimes_{l=1}^{i+j}\mathfrak{g}^{-\frac{1}{2}}\Big\Vert_{2}\,.
        \end{align}
    where $U\in(\mathbb C^d)^i$ $($resp. $V\in(\mathbb C^d)^j)$ are $i$ $($resp. $j)$ tangent vectors at points $\bs$ $($resp. $\bt)$, $\otimes$ denotes the tensor product and 
    \begin{equation}
    \label{eq:covariant_derivatives_ti_kernels}
    \forall T\in(\mathbb C^d)^{i+j}\,,\quad 
    \big(\nabla^{i+j}\rho(h)\bigtimes_{l=1}^{i+j}\mathfrak{g}^{-\frac{1}{2}}\big)[T]\defeq
    \nabla^{i+j}\rho(h)[\mathfrak{g}^{-\frac{1}{2}}T_1,\ldots,\mathfrak{g}^{-\frac{1}{2}}T_{i+j}]\,.
    \end{equation}

    In particular, one has
    \begin{align}
        T_1^\top \nabla\rho(h)          &=K^{(1,0)}(\bs, \bt)[T_1]=-K^{(0,1)}(\bs, \bt)[T_1]\,,\\
        \label{eq:Hessian_TVkernel}
        T_1^\top \nabla^2\rho(h)T_2     &=K^{(2,0)}(\bs, \bt)[T_1,T_2]=K^{(0,2)}(\bs, \bt)[T_1,T_2]=-K^{(1,1)}(\bs, \bt)[T_1,T_2]\,,\\
        \Big\Vert\mathfrak{g}^{-\frac{1}{2}}\nabla\rho(h)\Big\Vert_{2}
                                        &=\Vert K^{(1,0)}(\bs, \bt)\Vert_{s}\,,\\
        \Big\Vert\mathfrak{g}^{-\frac{1}{2}}\nabla^2\rho(h)\mathfrak{g}^{-\frac{1}{2}}\Big\Vert_{2}
                                        &=\Vert K^{(2,0)}(\bs, \bt)\Vert_{\bs,\bt}\,.
    \end{align}
    \end{subequations}

\subsection{Model RKHS and dual operator}

\begin{proposition}
\label{prop:steinwart}
Assume \eqref{hyp:continuous_model_kernel} holds true.
    \begin{itemize}
        \item The unique RKHS of $K_\model$ is given by 
        \begin{align}
        \label{eq:hmod_representation}
            \Hilbert_{\model}
            =
            \big\{
                \eta\,:\,\Param\to \R\ \big|\ \exists c\in\mathcal{F}\,,\ \eta=\eta_c
            \big\}
            =
            \big\{
                \eta\,:\,\Param\to \R\ \big|\ \exists! c\in{\overline{\mathrm{Im}(F)}}\,,\ \eta=\eta_c
            \big\}\,,
        \end{align}
        and for all $c$ orthogonal to ${\overline{\mathrm{Im}(F)}}$ in $\mathcal{F}$, $\eta_c=0$.
        \item The isometry is given by the mapping $c\in{\overline{\mathrm{Im}(F)}}\mapsto \eta_c\in\Hilbert_{\model}$.
        \item The norms satisfy 
        \[
    \forall c\in{\overline{\mathrm{Im}(F)}}\,,\quad 
    \|\eta_c\|_{\Hilbert_\model}=\|c\|_{\mathcal{F}}\,.
    \]
    \end{itemize}
\end{proposition}
\begin{proof}
    See \cite[Theorem 4.21]{steinwart2008support}.
\end{proof}

\subsection{Proof of Proposition~\ref{prop:model_description}}
\label{proof:prop_model_description}
\begin{subequations}\label[pluralequation]{eqs:f_star_properties}
By Assumption~\eqref{hyp:continuous_model_kernel}, it holds that $\Hilbert_\model\subset\mathcal C(\Param)$. Define $\mathrm{Id}_\model\,:\, \Hilbert_\model\to \mathcal C(\Param)$ the identity mapping. One has
\begin{align}
    F^\star c&=\mathrm{Id}_\model(\eta_c)\,,\textnormal{ for all } c\in{\overline{\mathrm{Im}(F)}}\,,\\
    \langle \mathrm{Id}_\model(\eta_c), \nu\rangle_{\mathcal C(\Param),\Measures(\Param)}
 &= \langle c, F\nu\rangle_{\mathcal{F}}\,,\textnormal{ for all } c\in\mathcal{F},\,\nu\in\Measures(\Param)\,.
\end{align}
\end{subequations}
The result follows from \eqref{eqs:f_star_properties} and Proposition~\ref{prop:steinwart}. For the sake of readability, we denote $\mathrm{Id}_\model(\eta_c)$ by $\eta_c$.

\section{Functional framework}

\subsection{Radon measures}
\label{app:Measures}

\begin{definition}[Set $\left(\Measures(\R^d), \Vert \cdot \Vert_{\mathrm TV} \right)$]
	We work in the vector space $\Measures(\R^d)$ of finite signed Radon measures on $\R^d$, endowed with the total variation norm:
	\begin{equation*}
	\Vert \mu \Vert_{\mathrm TV} \defeq \int_{\R^d} \diff \vert \mu \vert,
	\end{equation*}
	where $\vert \mu \vert = \mu^+ + \mu^-$ and $\mu = \mu^+ - \mu^-$ is the Jordan decomposition of $\mu \in \Measures(\R^d)$.
\end{definition}
By the Riesz-Markov theorem, this space is the dual of real-valued continuous functions on $\R^d$ \textit{vanishing at infinity}\footnote{$f$ vanishes at infinity iff $\forall \varepsilon > 0$, there exists a compact $K_{\varepsilon} \subset \R^d$ such that $ \forall x \in K_{\varepsilon}^c$, $\vert f(x) \vert < \varepsilon$} equipped with the supremum norm, which we denote as $\left(C_0(\R^d, \R), \Vert \cdot \Vert_{\infty}\right)$ or $\left(C_0(\R^d), \Vert \cdot \Vert_{\infty}\right)$. The total variation norm then admit a variational formulation as:
\begin{equation}
	\Vert \mu \Vert_{\mathrm TV} = \sup \left\{ \int \eta \diff \mu; \; \eta \in \mathcal{C}(\R^d), \,\Vert \eta \Vert_{\infty} \leq 1 \right\},
\end{equation}

\subsection{Fourier transform}
\label{app:FourierTransform}

Here, we detail the definition of the Fourier transform and its inverse used throughout the paper. 
\begin{definition}[Fourier transform]
	For any $g \in L^1(\R^d)$, we define its Fourier transform $\fourier\left[g\right]$ and its inverse as 
\begin{subequations}
    \begin{equation}
	\label{eq:FourierTransform}
	\forall \bomega \in \R^d, \quad \fourier\left[g\right](\bomega) = \int_{\R^d} g(\obs) e^{-i \bomega^\top \obs} \diff \obs.
	\end{equation}
    Furthermore, for $g$ continuous and both $g$ and $\fourier\left[g\right]$ in $L^1(\R^d)$, we define the inverse transform
	\begin{equation}
	\label{eq:InverseFourierTransform}
	\forall \obs \in \R^d, \quad \fourier^{-1}\Big[\fourier\big[g\big]\Big](\obs) = \frac{1}{(2\pi)^d} \int_{\R^d}   \fourier\left[g\right](\bomega) e^{i \bomega^\top \obs} \diff \bomega.
	\end{equation}
\end{subequations}
\end{definition}

\subsection{Reproducing kernel Hilbert space of translation invariant kernel}
\label{app:RKHS}
In this subsection, we remind particular results for translation invariant kernels $k(\bs,\bt) = \rho(\bt-\bs)$ and their characterisation in the Fourier domain. We begin by stating Bochner's theorem, a well known result characterizing functions of the positive definite type as the characteristic function of a positive measure, coined the \textit{spectral measure}.
\begin{theorem}[Bochner]
    \label{thm:Bochner}
    Let $k(\bs,\bt) = \rho(\bt-\bs)$. Then $k$ is a reproducing kernel if and only if there exists a positive measure $\nu$ such that
    \begin{equation*}
        \rho(\bt-\bs) = \int e^{+ \imath \bomega^\top (\bt - \bs)} \diff \nu(\bomega).
    \end{equation*}
    Furthermore, if the kernel is normalized, \textit{i.e.} $\rho(\bm{0}) = 1$, then $\nu$ is a probability distribution.
\end{theorem}
This allows to define the RKHS associated to the translation invariant kernel $\rho$ as functions with sufficiently fast decay rate in the Fourier domain.
\begin{proposition}[RKHS of translation invariant kernels]
    \label{prop:rkhs_of_ti_kernels}
    The RKHS of the translation invariant kernel $k(\bs, \bt) = \rho(\bt - \bs)$ is defined as 
    \begin{equation*}
        \mathcal{H}_\rho \defeq \left\{g \in L^2(\Param) \; : \; \Vert g \Vert_{\mathcal{H}_\rho}^2 = \frac{1}{(2\pi)^d} \int_{\Param} \frac{\left\vert \fourier[g](\bomega) \right\vert^2}{\fourier[\rho](\bomega) }\diff \bomega \right\}
    \end{equation*}
    with scalar product
    \begin{equation*}
        \<g_1,g_2\>_{\mathcal{H}_\rho} = \frac{1}{(2 \pi)^d} \int \frac{ \fourier[g_1](\bomega) \overline{\fourier[g_2](\bomega)}}{\fourier[\rho](\bomega) }\diff \bomega 
    \end{equation*}
\end{proposition}
\begin{proof}[Element of proof] 
    A formal proof is given in \textcite[][, Theorem 10.12, p. 139]{wendland2005scattered}, here we focus on proving the reproducing property. Moreover, we refer to \textcite[][Section 7.3.3, p.192]{bach2024learning} for readers interested in an intuitive construction of the RKHS norm. 
    
    First, we prove that $(\mathcal{H}_\rho, \< \cdot, \cdot \>_{\mathcal{H}_\rho})$ contains function of the form $k(\bs,  \cdot)$. To see this, one needs to remember that the Fourier transform is covariant by translation, meaning $\fourier[f(\cdot - \bs)](\bomega) = e^{-\imath \bomega^\top \bs} \fourier[f](\bomega)$. Hence, we have that 
    \begin{align*}
        \Vert k(\bs, \cdot) \Vert_{\mathcal{H}_\rho}  &= \Vert \rho(\cdot - \bs) \Vert_{\mathcal{H}_\rho}  = \frac{1}{(2\pi)^d}\int   \frac{\left\vert \fourier[\rho](\bomega) \right\vert^2}{\fourier[\rho](\bomega)} \vert e^{-\imath \bomega^\top \bs}\vert^2 \diff \bomega = \int   \frac{\fourier[\rho](\bomega)}{(2\pi)^d} \diff \bomega = \rho(\bm{0}) < + \infty .
    \end{align*}
    Thus, $k(\bs, \cdot)$ belongs to the RKHS. Following the same line, we can prove that the reproducing property holds. Let $g \in \mathcal{H}_{\rho}$ and $\bs \in \mathbb{R}^d$
    \begin{align*}
         \<g, \rho(\cdot - \bs) \>_{\mathcal{H}_\rho}  &= \frac{1}{(2 \pi)^d} \int \frac{ \fourier[g](\bomega) \overline{\fourier[\rho](\bomega)}}{\fourier[\rho](\bomega)} e^{+ \imath \bomega^\top \bs}\diff \bomega = \fourier^{-1}[\fourier[g]](\bs) = g(\bs).
    \end{align*}

    Completeness is treated in \textcite[][, Theorem 10.12, p. 139]{wendland2005scattered}.
\end{proof}

\begin{remark}[Link between the Fourier transform of $\rho$ and its spectral measure $\nu$]
    With the choice of convention introduced in~\Cref{app:FourierTransform} for the inverse Fourier transform, the following identity holds
    \begin{equation*}
        \rho = (2 \pi)^d \fourier^{-1}[\nu] \iff \nu = \frac{\fourier[\rho]}{(2\pi)^d}.
    \end{equation*}
    As a consequence we have 
    \begin{align*}
        \Vert \rho \star g \Vert^2_{\mathcal{H}_\rho} & =  \frac{1}{(2\pi)^d} \int_{\Param} \frac{\left\vert \fourier[\rho](\bomega) \right\vert^2 \left\vert \fourier[g](\bomega) \right\vert^2}{\fourier[\rho](\bomega) }\diff \bomega, \\
        &=   \int_{\Param} \frac{\fourier[\rho](\bomega) }{(2\pi)^d} \left\vert \fourier[g](\bomega) \right\vert^2 \diff \bomega, \\
        &= \int_{\Param} \nu(\bomega) \left\vert \fourier[g](\bomega) \right\vert^2 \diff \bomega, \\ 
    \end{align*}
    If the kernel is normalized, \textit{i.e.} $\rho(0) = 1$, then $\nu$ is a probability distribution and the last line can be written as an expectation.
\end{remark}

\subsection{Fourier transform and spectral measure of the sinus cardinal smoothing kernel $\lambda_\bandwidth$}
We recall our choice of smoothing kernel $\lambda_\bandwidth$ which is the rescaled sinus cardinal
\[
\lambda_\bandwidth(\bx) =  \Big(\frac{1}{\bandwidth}\Big)^d \sinc\Big(\frac{\bx}{\bandwidth}\Big),
\]
with its spectral measure from \Cref{thm:Bochner} denoted as $U_\bandwidth$. Denoting~${\sinc_\bandwidth = \sinc( \cdot / \bandwidth)}$, we have with our choice of Fourier conventions that $\Pi_1(\bomega) \defeq \fourier[\sinc_1] = \pi^d \indicator_{[-1, 1]^d}$. Thus, the Fourier transform of its scaled version is the gate function
\[
\Pi_\bandwidth(\bomega) \defeq \fourier[\sinc_\bandwidth](\bomega) = \fourier[\sinc_1(\cdot/\bandwidth)](\bomega) = (\pi \bandwidth)^d \indicator_{[-1/\bandwidth, 1/\bandwidth]^d}(\bomega),
\] 
with a frequency cut-off as $1/\bandwidth$. Finally, by linearity, the Fourier transform of the smoothing kernel is
\begin{equation}
    \fourier[\lambda_\bandwidth](\bomega) = \pi^d \indicator_{[-1/\bandwidth, 1/\bandwidth]^d}(\bomega),
\end{equation}
ans its spectral measure is
\begin{equation}    
U_\bandwidth(\bomega) = \frac{\fourier[\lambda_\bandwidth]}{(2\pi)^d}(\bomega) = \frac{1}{2^d} \indicator_{[-1/\bandwidth, 1/\bandwidth]^d}(\bomega).
\end{equation}

\subsection{Spectral measure of the sinc-4 kernel}
\label{app:sinc4_RKHS}
The sinc-4 kernel $\Psi_{\bandwidth}(\bx - \by) = \sinc(\frac{\bx-\by}{4 \bandwidth})^4 = \sinc_{4\bandwidth}(\bx - \by)^4$ defined in~\Cref{eq:sketch_kernel_sinc} is a translation invariant kernel with associated RKHS given by
\begin{align}
\Hilbert_{\sinc} & \defeq \left\{ g : \R^d \to \R \textnormal{ s.t. } \Vert g \Vert_{\Hilbert_{{\sinc}}}^2  = \frac{1}{(2\pi)^d} 
\int_{\R^d} \frac{\vert \fourier[g] \vert^2}{\fourier[\Psi_\bandwidth]} < + \infty \right\}\,.
\end{align}  
We denote its spectral measure as~${f^{(4)}_{\bandwidth} = \fourier[\Psi_\bandwidth] / (2\pi)^d}$, which verifies by Bochner's theorem (\Cref{thm:Bochner})
\[
\Psi_\bandwidth(\bx) = \int_{\mathbb{R}^d} e^{+ \imath \bx^\top \bomega} f^{(4)}_{\bandwidth}(\bomega) \diff \bomega.
\]
Since $\Psi_\bandwidth(\bm{0})=1$, then $f^{(4)}_{\bandwidth}$ is a normalized probability density function. Moreover, it has an explicit analytical form. Indeed, remember that with our choice of Fourier conventions, the function~${\sinc_{4\bandwidth} = \sinc( \cdot / 4 \bandwidth)}$ is the inverse Fourier transform of the rescaled low-pass filter $\Pi_{4\bandwidth}(\bomega) = (4 \pi \bandwidth)^d \indicator_{[-1/4\bandwidth, 1/4 \bandwidth]^d}(\bomega)$. Thus, $f^{(4)}_{\bandwidth}$ is proportional to $\Pi_{4\bandwidth}^{(\star 4)}$, which we define as the convolution of $\Pi_{4\bandwidth}$ with itself taken four time. The proportionality constant is such that $f^{(4)}_{\bandwidth}$ is a normalized probability density function. Precisely, we know that $f^{(4)}_{\bandwidth}$ is even since $\Psi_\bandwidth$ is even and real-valued, so that
\[
f^{(4)}_{\bandwidth} 
= \check{f}^{(4)}_{\bandwidth} 
= \fourier^{-1}[ \Psi_\bandwidth] 
= \fourier^{-1}[ \sinc_{4\bandwidth} ^4] 
= \fourier^{-1}[ \fourier^{-1}[ \Pi_{4\bandwidth}]^4] 
= (2\pi)^{-4d} \fourier^{-1}[\fourier[ \Pi_{4\bandwidth}](- \cdot)^4] 
= (2\pi)^{-4d} \Pi_{4\bandwidth}^{(\star 4)}
\]
where we used that $g^{(\star n)} = \fourier^{-1}[ \fourier[g]^n]$ and $\fourier[ \Pi_{4\bandwidth}]$ is even.

Note that $f^{(4)}_{\bandwidth}$ may also be interpreted as the density of the sum of independent uniform random variable on~${[-1/4\bandwidth, 1/4\bandwidth]^d}$, also known as the (rescaled) Irwin-Hall distribution of order 4. It is even by symmetricity around $\bm{0}$ and its expression is
\begin{equation}
    \forall \bomega\in[1/\bandwidth, 1/\bandwidth]^d, \quad f^{(4)}_{\bandwidth}(\bomega) = {(2\bandwidth)^d} \prod_{j=1}^d \frac{1}{6} \sum_{k=1}^{\lfloor 2\bandwidth \omega_j + 2 \rfloor} (-1)^k \binom{4}{k} \big( 2\bandwidth \omega_j + 2 - k \big)^4
\end{equation}
Moreover, its mode is in $\bomega=\bm{0}$ and its maximum value $(\frac{2}{3})^d (2 \bandwidth)^d$ so that 
\[
\sup_{\bomega \in [-1/\bandwidth, 1/\bandwidth]^d} \fourier[\Psi_\bandwidth](\bomega) = (2\pi)^d\sup_{\bomega \in [-1/\bandwidth, 1/\bandwidth]^d} f^{(4)}_\bandwidth(\bomega) = (2\pi)^d \Big(\frac{2}{3}\Big)^d (2 \bandwidth)^d= \mathcal{O}(\bandwidth^d)
\]

\subsection{Kernel switch between translation invariant RKHS}
\label{app:cswitch_ti_kernel}

A direct consequence of Bochner's theorem and \Cref{prop:rkhs_of_ti_kernels} is that we can explicitly give the embedding constant $C_\switch$ between the RKHS of two translation invariant kernels. This is discussed \textit{e.g.} in \textcite[Corollary 3.2]{zhang2013inclusion} and the constant is precisely the essential supremum of the ratio of their Fourier transform as introduced in~\Cref{hyp:Hphi}.  

Considering two continuous and translation-invariant kernels.
\begin{align*}
    K_\pivot(\bs, \bt) &= \rho_\pivot(\bt - \bs) = \frac{1}{(2\pi)^d} \int e^{+\imath \bomega^\top (\bt -\bs)} \fourier[\rho_\pivot](\bomega) \diff \bomega,\\
   K_\model(\bs, \bt) &= \rho_\model(\bt - \bs) = \frac{1}{(2\pi)^d} \int e^{+\imath \bomega^\top (\bt -\bs)} \fourier[\rho_\model](\bomega)  \diff \bomega
\end{align*}
Then, assuming the support of $\fourier[\rho_\model]$ contains that of $\fourier[\rho_\pivot]$, we have for any $\eta \in \Hilbert_\pivot$ that
\begin{equation*}
    \Vert \eta \Vert^2_{\Hilbert_\model}\!\!\! = \frac{1}{(2\pi)^d} \int \frac{\big|\fourier[\eta]\big|^2}{\fourier[\rho_\model]}(\bomega) \diff \bomega = \frac{1}{(2\pi)^d} \int \frac{\fourier[\rho_\pivot]}{\fourier[\rho_\model]}(\bomega) \frac{\big|\fourier[\eta]\big|^2}{\fourier[\rho_\pivot]}(\bomega) \diff \bomega \leq\!\!\!\!\!\! \essentialsup\limits_{\bomega \in \Supp \fourier[\rho_\pivot]} \frac{\fourier[\rho_\pivot]}{\fourier[\rho_\model]}(\bomega) \cdot \Vert \eta \Vert^2_{\Hilbert_\pivot}. %
\end{equation*}
Thus, 
\begin{equation}
    \label{eq:cswitch_ti_kernels}
    \frac{\Vert \eta \Vert_{\Hilbert_\model}}{\Vert \eta \Vert_{\Hilbert_\pivot}} \leq \essentialsup\limits_{\bomega \in \Supp \fourier[\rho_\pivot]}  \sqrt{\frac{\fourier[\rho_\pivot]}{\fourier[\rho_\model]}(\bomega)} =: C_\switch .
\end{equation}

\subsection{Proof of~\Cref{prop:sinc4_cswitch}}
\label{app:proof_sinc4_csiwtch}
The proposition is a direct consequence of the preceding section above. Plugging $K_\pivot = \Psi_\bandwidth$ the sinc-4 kernel, and letting $K_\model = \rho$ be any translation-invariant kernel such that $\Supp \fourier[\Psi_\bandwidth] = [-1/\bandwidth, 1/\bandwidth]^d \subset \Supp \fourier[\rho_\model]$. Then,~\Cref{eq:cswitch_ti_kernels} becomes
\begin{equation}
    C_\switch(\Psi_\bandwidth, \rho) = \essentialsup_{\bomega \in [-1/\bandwidth, 1/\bandwidth]^d}\sqrt{
\frac{\fourier[\Psi_\bandwidth]}{\fourier[\rho_\model]}(\bomega)}.
\end{equation}
In addition, we have the straightforward upper-bound:
\[
C_\switch(\Psi_\bandwidth, \rho) \leq  \frac{
\essentialsup_{\bomega \in [-1/\bandwidth, 1/\bandwidth]^d} \sqrt{\fourier[\Psi_\bandwidth](\bomega)}
}{\essentialinf_{\bomega \in [-1/\bandwidth, 1/\bandwidth]^d } \sqrt{\fourier[\rho_\model](\bomega)}}.
\]
We know that $\fourier[\Psi_\bandwidth]$ is upper-bounded with a mode $\max \fourier[\Psi_\bandwidth] = (\frac{8\pi}{3} \bandwidth)^d$. Hence, the numerator scales as $\mathcal{O}(\bandwidth^{d/2})$ and a sufficient condition for $C_\switch(\Psi_\bandwidth, \rho)$ to be finite is that $\fourier[\rho_\model]$ is essentially bounded away from zero on  $[-1/\bandwidth, 1/\bandwidth]^d$.

\end{document}